\documentclass[aos,preprint]{imsart}

\RequirePackage[OT1]{fontenc}
\RequirePackage{amsthm,amsmath,natbib}
\RequirePackage[colorlinks,citecolor=blue,urlcolor=blue]{hyperref}

\usepackage{color}  
\definecolor{AliceBlue}{rgb}{0.94,0.97,1.00}
\definecolor{AntiqueWhite1}{rgb}{1.00,0.94,0.86}
\definecolor{AntiqueWhite2}{rgb}{0.93,0.87,0.80}
\definecolor{AntiqueWhite3}{rgb}{0.80,0.75,0.69}
\definecolor{AntiqueWhite4}{rgb}{0.55,0.51,0.47}
\definecolor{AntiqueWhite}{rgb}{0.98,0.92,0.84}
\definecolor{BlanchedAlmond}{rgb}{1.00,0.92,0.80}
\definecolor{BlueViolet}{rgb}{0.54,0.17,0.89}
\definecolor{CadetBlue1}{rgb}{0.60,0.96,1.00}
\definecolor{CadetBlue2}{rgb}{0.56,0.90,0.93}
\definecolor{CadetBlue3}{rgb}{0.48,0.77,0.80}
\definecolor{CadetBlue4}{rgb}{0.33,0.53,0.55}
\definecolor{CadetBlue}{rgb}{0.37,0.62,0.63}
\definecolor{CornflowerBlue}{rgb}{0.39,0.58,0.93}
\definecolor{DarkBlue}{rgb}{0.00,0.00,0.55}
\definecolor{DarkCyan}{rgb}{0.00,0.55,0.55}
\definecolor{DarkGoldenrod1}{rgb}{1.00,0.73,0.06}
\definecolor{DarkGoldenrod2}{rgb}{0.93,0.68,0.05}
\definecolor{DarkGoldenrod3}{rgb}{0.80,0.58,0.05}
\definecolor{DarkGoldenrod4}{rgb}{0.55,0.40,0.03}
\definecolor{DarkGoldenrod}{rgb}{0.72,0.53,0.04}
\definecolor{DarkGray}{rgb}{0.66,0.66,0.66}
\definecolor{DarkGreen}{rgb}{0.00,0.39,0.00}
\definecolor{DarkGrey}{rgb}{0.66,0.66,0.66}
\definecolor{DarkKhaki}{rgb}{0.74,0.72,0.42}
\definecolor{DarkMagenta}{rgb}{0.55,0.00,0.55}
\definecolor{DarkOliveGreen1}{rgb}{0.79,1.00,0.44}
\definecolor{DarkOliveGreen2}{rgb}{0.74,0.93,0.41}
\definecolor{DarkOliveGreen3}{rgb}{0.64,0.80,0.35}
\definecolor{DarkOliveGreen4}{rgb}{0.43,0.55,0.24}
\definecolor{DarkOliveGreen}{rgb}{0.33,0.42,0.18}
\definecolor{DarkOrange1}{rgb}{1.00,0.50,0.00}
\definecolor{DarkOrange2}{rgb}{0.93,0.46,0.00}
\definecolor{DarkOrange3}{rgb}{0.80,0.40,0.00}
\definecolor{DarkOrange4}{rgb}{0.55,0.27,0.00}
\definecolor{DarkOrange}{rgb}{1.00,0.55,0.00}
\definecolor{DarkOrchid1}{rgb}{0.75,0.24,1.00}
\definecolor{DarkOrchid2}{rgb}{0.70,0.23,0.93}
\definecolor{DarkOrchid3}{rgb}{0.60,0.20,0.80}
\definecolor{DarkOrchid4}{rgb}{0.41,0.13,0.55}
\definecolor{DarkOrchid}{rgb}{0.60,0.20,0.80}
\definecolor{DarkRed}{rgb}{0.55,0.00,0.00}
\definecolor{DarkSalmon}{rgb}{0.91,0.59,0.48}
\definecolor{DarkSeaGreen1}{rgb}{0.76,1.00,0.76}
\definecolor{DarkSeaGreen2}{rgb}{0.71,0.93,0.71}
\definecolor{DarkSeaGreen3}{rgb}{0.61,0.80,0.61}
\definecolor{DarkSeaGreen4}{rgb}{0.41,0.55,0.41}
\definecolor{DarkSeaGreen}{rgb}{0.56,0.74,0.56}
\definecolor{DarkSlateBlue}{rgb}{0.28,0.24,0.55}
\definecolor{DarkSlateGray1}{rgb}{0.59,1.00,1.00}
\definecolor{DarkSlateGray2}{rgb}{0.55,0.93,0.93}
\definecolor{DarkSlateGray3}{rgb}{0.47,0.80,0.80}
\definecolor{DarkSlateGray4}{rgb}{0.32,0.55,0.55}
\definecolor{DarkSlateGray}{rgb}{0.18,0.31,0.31}
\definecolor{DarkSlateGrey}{rgb}{0.18,0.31,0.31}
\definecolor{DarkTurquoise}{rgb}{0.00,0.81,0.82}
\definecolor{DarkViolet}{rgb}{0.58,0.00,0.83}
\definecolor{DeepPink1}{rgb}{1.00,0.08,0.58}
\definecolor{DeepPink2}{rgb}{0.93,0.07,0.54}
\definecolor{DeepPink3}{rgb}{0.80,0.06,0.46}
\definecolor{DeepPink4}{rgb}{0.55,0.04,0.31}
\definecolor{DeepPink}{rgb}{1.00,0.08,0.58}
\definecolor{DeepSkyBlue1}{rgb}{0.00,0.75,1.00}
\definecolor{DeepSkyBlue2}{rgb}{0.00,0.70,0.93}
\definecolor{DeepSkyBlue3}{rgb}{0.00,0.60,0.80}
\definecolor{DeepSkyBlue4}{rgb}{0.00,0.41,0.55}
\definecolor{DeepSkyBlue}{rgb}{0.00,0.75,1.00}
\definecolor{DimGray}{rgb}{0.41,0.41,0.41}
\definecolor{DimGrey}{rgb}{0.41,0.41,0.41}
\definecolor{DodgerBlue1}{rgb}{0.12,0.56,1.00}
\definecolor{DodgerBlue2}{rgb}{0.11,0.53,0.93}
\definecolor{DodgerBlue3}{rgb}{0.09,0.45,0.80}
\definecolor{DodgerBlue4}{rgb}{0.06,0.31,0.55}
\definecolor{DodgerBlue}{rgb}{0.12,0.56,1.00}
\definecolor{FloralWhite}{rgb}{1.00,0.98,0.94}
\definecolor{ForestGreen}{rgb}{0.13,0.55,0.13}
\definecolor{GhostWhite}{rgb}{0.97,0.97,1.00}
\definecolor{GreenYellow}{rgb}{0.68,1.00,0.18}
\definecolor{HotPink1}{rgb}{1.00,0.43,0.71}
\definecolor{HotPink2}{rgb}{0.93,0.42,0.65}
\definecolor{HotPink3}{rgb}{0.80,0.38,0.56}
\definecolor{HotPink4}{rgb}{0.55,0.23,0.38}
\definecolor{HotPink}{rgb}{1.00,0.41,0.71}
\definecolor{IndianRed1}{rgb}{1.00,0.42,0.42}
\definecolor{IndianRed2}{rgb}{0.93,0.39,0.39}
\definecolor{IndianRed3}{rgb}{0.80,0.33,0.33}
\definecolor{IndianRed4}{rgb}{0.55,0.23,0.23}
\definecolor{IndianRed}{rgb}{0.80,0.36,0.36}
\definecolor{LavenderBlush1}{rgb}{1.00,0.94,0.96}
\definecolor{LavenderBlush2}{rgb}{0.93,0.88,0.90}
\definecolor{LavenderBlush3}{rgb}{0.80,0.76,0.77}
\definecolor{LavenderBlush4}{rgb}{0.55,0.51,0.53}
\definecolor{LavenderBlush}{rgb}{1.00,0.94,0.96}
\definecolor{LawnGreen}{rgb}{0.49,0.99,0.00}
\definecolor{LemonChiffon1}{rgb}{1.00,0.98,0.80}
\definecolor{LemonChiffon2}{rgb}{0.93,0.91,0.75}
\definecolor{LemonChiffon3}{rgb}{0.80,0.79,0.65}
\definecolor{LemonChiffon4}{rgb}{0.55,0.54,0.44}
\definecolor{LemonChiffon}{rgb}{1.00,0.98,0.80}
\definecolor{LightBlue1}{rgb}{0.75,0.94,1.00}
\definecolor{LightBlue2}{rgb}{0.70,0.87,0.93}
\definecolor{LightBlue3}{rgb}{0.60,0.75,0.80}
\definecolor{LightBlue4}{rgb}{0.41,0.51,0.55}
\definecolor{LightBlue}{rgb}{0.68,0.85,0.90}
\definecolor{LightCoral}{rgb}{0.94,0.50,0.50}
\definecolor{LightCyan1}{rgb}{0.88,1.00,1.00}
\definecolor{LightCyan2}{rgb}{0.82,0.93,0.93}
\definecolor{LightCyan3}{rgb}{0.71,0.80,0.80}
\definecolor{LightCyan4}{rgb}{0.48,0.55,0.55}
\definecolor{LightCyan}{rgb}{0.88,1.00,1.00}
\definecolor{LightGoldenrod1}{rgb}{1.00,0.93,0.55}
\definecolor{LightGoldenrod2}{rgb}{0.93,0.86,0.51}
\definecolor{LightGoldenrod3}{rgb}{0.80,0.75,0.44}
\definecolor{LightGoldenrod4}{rgb}{0.55,0.51,0.30}
\definecolor{LightGoldenrodYellow}{rgb}{0.98,0.98,0.82}
\definecolor{LightGoldenrod}{rgb}{0.93,0.87,0.51}
\definecolor{LightGray}{rgb}{0.83,0.83,0.83}
\definecolor{LightGreen}{rgb}{0.56,0.93,0.56}
\definecolor{LightGrey}{rgb}{0.83,0.83,0.83}
\definecolor{LightPink1}{rgb}{1.00,0.68,0.73}
\definecolor{LightPink2}{rgb}{0.93,0.64,0.68}
\definecolor{LightPink3}{rgb}{0.80,0.55,0.58}
\definecolor{LightPink4}{rgb}{0.55,0.37,0.40}
\definecolor{LightPink}{rgb}{1.00,0.71,0.76}
\definecolor{LightSalmon1}{rgb}{1.00,0.63,0.48}
\definecolor{LightSalmon2}{rgb}{0.93,0.58,0.45}
\definecolor{LightSalmon3}{rgb}{0.80,0.51,0.38}
\definecolor{LightSalmon4}{rgb}{0.55,0.34,0.26}
\definecolor{LightSalmon}{rgb}{1.00,0.63,0.48}
\definecolor{LightSeaGreen}{rgb}{0.13,0.70,0.67}
\definecolor{LightSkyBlue1}{rgb}{0.69,0.89,1.00}
\definecolor{LightSkyBlue2}{rgb}{0.64,0.83,0.93}
\definecolor{LightSkyBlue3}{rgb}{0.55,0.71,0.80}
\definecolor{LightSkyBlue4}{rgb}{0.38,0.48,0.55}
\definecolor{LightSkyBlue}{rgb}{0.53,0.81,0.98}
\definecolor{LightSlateBlue}{rgb}{0.52,0.44,1.00}
\definecolor{LightSlateGray}{rgb}{0.47,0.53,0.60}
\definecolor{LightSlateGrey}{rgb}{0.47,0.53,0.60}
\definecolor{LightSteelBlue1}{rgb}{0.79,0.88,1.00}
\definecolor{LightSteelBlue2}{rgb}{0.74,0.82,0.93}
\definecolor{LightSteelBlue3}{rgb}{0.64,0.71,0.80}
\definecolor{LightSteelBlue4}{rgb}{0.43,0.48,0.55}
\definecolor{LightSteelBlue}{rgb}{0.69,0.77,0.87}
\definecolor{LightYellow1}{rgb}{1.00,1.00,0.88}
\definecolor{LightYellow2}{rgb}{0.93,0.93,0.82}
\definecolor{LightYellow3}{rgb}{0.80,0.80,0.71}
\definecolor{LightYellow4}{rgb}{0.55,0.55,0.48}
\definecolor{LightYellow}{rgb}{1.00,1.00,0.88}
\definecolor{LimeGreen}{rgb}{0.20,0.80,0.20}
\definecolor{MediumAquamarine}{rgb}{0.40,0.80,0.67}
\definecolor{MediumBlue}{rgb}{0.00,0.00,0.80}
\definecolor{MediumOrchid1}{rgb}{0.88,0.40,1.00}
\definecolor{MediumOrchid2}{rgb}{0.82,0.37,0.93}
\definecolor{MediumOrchid3}{rgb}{0.71,0.32,0.80}
\definecolor{MediumOrchid4}{rgb}{0.48,0.22,0.55}
\definecolor{MediumOrchid}{rgb}{0.73,0.33,0.83}
\definecolor{MediumPurple1}{rgb}{0.67,0.51,1.00}
\definecolor{MediumPurple2}{rgb}{0.62,0.47,0.93}
\definecolor{MediumPurple3}{rgb}{0.54,0.41,0.80}
\definecolor{MediumPurple4}{rgb}{0.36,0.28,0.55}
\definecolor{MediumPurple}{rgb}{0.58,0.44,0.86}
\definecolor{MediumSeaGreen}{rgb}{0.24,0.70,0.44}
\definecolor{MediumSlateBlue}{rgb}{0.48,0.41,0.93}
\definecolor{MediumSpringGreen}{rgb}{0.00,0.98,0.60}
\definecolor{MediumTurquoise}{rgb}{0.28,0.82,0.80}
\definecolor{MediumVioletRed}{rgb}{0.78,0.08,0.52}
\definecolor{MidnightBlue}{rgb}{0.10,0.10,0.44}
\definecolor{MintCream}{rgb}{0.96,1.00,0.98}
\definecolor{MistyRose1}{rgb}{1.00,0.89,0.88}
\definecolor{MistyRose2}{rgb}{0.93,0.84,0.82}
\definecolor{MistyRose3}{rgb}{0.80,0.72,0.71}
\definecolor{MistyRose4}{rgb}{0.55,0.49,0.48}
\definecolor{MistyRose}{rgb}{1.00,0.89,0.88}
\definecolor{NavajoWhite1}{rgb}{1.00,0.87,0.68}
\definecolor{NavajoWhite2}{rgb}{0.93,0.81,0.63}
\definecolor{NavajoWhite3}{rgb}{0.80,0.70,0.55}
\definecolor{NavajoWhite4}{rgb}{0.55,0.47,0.37}
\definecolor{NavajoWhite}{rgb}{1.00,0.87,0.68}
\definecolor{NavyBlue}{rgb}{0.00,0.00,0.50}
\definecolor{OldLace}{rgb}{0.99,0.96,0.90}
\definecolor{OliveDrab1}{rgb}{0.75,1.00,0.24}
\definecolor{OliveDrab2}{rgb}{0.70,0.93,0.23}
\definecolor{OliveDrab3}{rgb}{0.60,0.80,0.20}
\definecolor{OliveDrab4}{rgb}{0.41,0.55,0.13}
\definecolor{OliveDrab}{rgb}{0.42,0.56,0.14}
\definecolor{OrangeRed1}{rgb}{1.00,0.27,0.00}
\definecolor{OrangeRed2}{rgb}{0.93,0.25,0.00}
\definecolor{OrangeRed3}{rgb}{0.80,0.22,0.00}
\definecolor{OrangeRed4}{rgb}{0.55,0.15,0.00}
\definecolor{OrangeRed}{rgb}{1.00,0.27,0.00}
\definecolor{PaleGoldenrod}{rgb}{0.93,0.91,0.67}
\definecolor{PaleGreen1}{rgb}{0.60,1.00,0.60}
\definecolor{PaleGreen2}{rgb}{0.56,0.93,0.56}
\definecolor{PaleGreen3}{rgb}{0.49,0.80,0.49}
\definecolor{PaleGreen4}{rgb}{0.33,0.55,0.33}
\definecolor{PaleGreen}{rgb}{0.60,0.98,0.60}
\definecolor{PaleTurquoise1}{rgb}{0.73,1.00,1.00}
\definecolor{PaleTurquoise2}{rgb}{0.68,0.93,0.93}
\definecolor{PaleTurquoise3}{rgb}{0.59,0.80,0.80}
\definecolor{PaleTurquoise4}{rgb}{0.40,0.55,0.55}
\definecolor{PaleTurquoise}{rgb}{0.69,0.93,0.93}
\definecolor{PaleVioletRed1}{rgb}{1.00,0.51,0.67}
\definecolor{PaleVioletRed2}{rgb}{0.93,0.47,0.62}
\definecolor{PaleVioletRed3}{rgb}{0.80,0.41,0.54}
\definecolor{PaleVioletRed4}{rgb}{0.55,0.28,0.36}
\definecolor{PaleVioletRed}{rgb}{0.86,0.44,0.58}
\definecolor{PapayaWhip}{rgb}{1.00,0.94,0.84}
\definecolor{PeachPuff1}{rgb}{1.00,0.85,0.73}
\definecolor{PeachPuff2}{rgb}{0.93,0.80,0.68}
\definecolor{PeachPuff3}{rgb}{0.80,0.69,0.58}
\definecolor{PeachPuff4}{rgb}{0.55,0.47,0.40}
\definecolor{PeachPuff}{rgb}{1.00,0.85,0.73}
\definecolor{PowderBlue}{rgb}{0.69,0.88,0.90}
\definecolor{RosyBrown1}{rgb}{1.00,0.76,0.76}
\definecolor{RosyBrown2}{rgb}{0.93,0.71,0.71}
\definecolor{RosyBrown3}{rgb}{0.80,0.61,0.61}
\definecolor{RosyBrown4}{rgb}{0.55,0.41,0.41}
\definecolor{RosyBrown}{rgb}{0.74,0.56,0.56}
\definecolor{RoyalBlue1}{rgb}{0.28,0.46,1.00}
\definecolor{RoyalBlue2}{rgb}{0.26,0.43,0.93}
\definecolor{RoyalBlue3}{rgb}{0.23,0.37,0.80}
\definecolor{RoyalBlue4}{rgb}{0.15,0.25,0.55}
\definecolor{RoyalBlue}{rgb}{0.25,0.41,0.88}
\definecolor{SaddleBrown}{rgb}{0.55,0.27,0.07}
\definecolor{SandyBrown}{rgb}{0.96,0.64,0.38}
\definecolor{SeaGreen1}{rgb}{0.33,1.00,0.62}
\definecolor{SeaGreen2}{rgb}{0.31,0.93,0.58}
\definecolor{SeaGreen3}{rgb}{0.26,0.80,0.50}
\definecolor{SeaGreen4}{rgb}{0.18,0.55,0.34}
\definecolor{SeaGreen}{rgb}{0.18,0.55,0.34}
\definecolor{SkyBlue1}{rgb}{0.53,0.81,1.00}
\definecolor{SkyBlue2}{rgb}{0.49,0.75,0.93}
\definecolor{SkyBlue3}{rgb}{0.42,0.65,0.80}
\definecolor{SkyBlue4}{rgb}{0.29,0.44,0.55}
\definecolor{SkyBlue}{rgb}{0.53,0.81,0.92}
\definecolor{SlateBlue1}{rgb}{0.51,0.44,1.00}
\definecolor{SlateBlue2}{rgb}{0.48,0.40,0.93}
\definecolor{SlateBlue3}{rgb}{0.41,0.35,0.80}
\definecolor{SlateBlue4}{rgb}{0.28,0.24,0.55}
\definecolor{SlateBlue}{rgb}{0.42,0.35,0.80}
\definecolor{SlateGray1}{rgb}{0.78,0.89,1.00}
\definecolor{SlateGray2}{rgb}{0.73,0.83,0.93}
\definecolor{SlateGray3}{rgb}{0.62,0.71,0.80}
\definecolor{SlateGray4}{rgb}{0.42,0.48,0.55}
\definecolor{SlateGray}{rgb}{0.44,0.50,0.56}
\definecolor{SlateGrey}{rgb}{0.44,0.50,0.56}
\definecolor{SpringGreen1}{rgb}{0.00,1.00,0.50}
\definecolor{SpringGreen2}{rgb}{0.00,0.93,0.46}
\definecolor{SpringGreen3}{rgb}{0.00,0.80,0.40}
\definecolor{SpringGreen4}{rgb}{0.00,0.55,0.27}
\definecolor{SpringGreen}{rgb}{0.00,1.00,0.50}
\definecolor{SteelBlue1}{rgb}{0.39,0.72,1.00}
\definecolor{SteelBlue2}{rgb}{0.36,0.67,0.93}
\definecolor{SteelBlue3}{rgb}{0.31,0.58,0.80}
\definecolor{SteelBlue4}{rgb}{0.21,0.39,0.55}
\definecolor{SteelBlue}{rgb}{0.27,0.51,0.71}
\definecolor{VioletRed1}{rgb}{1.00,0.24,0.59}
\definecolor{VioletRed2}{rgb}{0.93,0.23,0.55}
\definecolor{VioletRed3}{rgb}{0.80,0.20,0.47}
\definecolor{VioletRed4}{rgb}{0.55,0.13,0.32}
\definecolor{VioletRed}{rgb}{0.82,0.13,0.56}
\definecolor{WhiteSmoke}{rgb}{0.96,0.96,0.96}
\definecolor{YellowGreen}{rgb}{0.60,0.80,0.20}
\definecolor{aliceblue}{rgb}{0.94,0.97,1.00}
\definecolor{antiquewhite}{rgb}{0.98,0.92,0.84}
\definecolor{aquamarine1}{rgb}{0.50,1.00,0.83}
\definecolor{aquamarine2}{rgb}{0.46,0.93,0.78}
\definecolor{aquamarine3}{rgb}{0.40,0.80,0.67}
\definecolor{aquamarine4}{rgb}{0.27,0.55,0.45}
\definecolor{aquamarine}{rgb}{0.50,1.00,0.83}
\definecolor{azure1}{rgb}{0.94,1.00,1.00}
\definecolor{azure2}{rgb}{0.88,0.93,0.93}
\definecolor{azure3}{rgb}{0.76,0.80,0.80}
\definecolor{azure4}{rgb}{0.51,0.55,0.55}
\definecolor{azure}{rgb}{0.94,1.00,1.00}
\definecolor{beige}{rgb}{0.96,0.96,0.86}
\definecolor{bisque1}{rgb}{1.00,0.89,0.77}
\definecolor{bisque2}{rgb}{0.93,0.84,0.72}
\definecolor{bisque3}{rgb}{0.80,0.72,0.62}
\definecolor{bisque4}{rgb}{0.55,0.49,0.42}
\definecolor{bisque}{rgb}{1.00,0.89,0.77}
\definecolor{black}{rgb}{0.00,0.00,0.00}
\definecolor{blanchedalmond}{rgb}{1.00,0.92,0.80}
\definecolor{blue1}{rgb}{0.00,0.00,1.00}
\definecolor{blue2}{rgb}{0.00,0.00,0.93}
\definecolor{blue3}{rgb}{0.00,0.00,0.80}
\definecolor{blue4}{rgb}{0.00,0.00,0.55}
\definecolor{blueviolet}{rgb}{0.54,0.17,0.89}
\definecolor{blue}{rgb}{0.00,0.00,1.00}
\definecolor{brown1}{rgb}{1.00,0.25,0.25}
\definecolor{brown2}{rgb}{0.93,0.23,0.23}
\definecolor{brown3}{rgb}{0.80,0.20,0.20}
\definecolor{brown4}{rgb}{0.55,0.14,0.14}
\definecolor{brown}{rgb}{0.65,0.16,0.16}
\definecolor{burlywood1}{rgb}{1.00,0.83,0.61}
\definecolor{burlywood2}{rgb}{0.93,0.77,0.57}
\definecolor{burlywood3}{rgb}{0.80,0.67,0.49}
\definecolor{burlywood4}{rgb}{0.55,0.45,0.33}
\definecolor{burlywood}{rgb}{0.87,0.72,0.53}
\definecolor{cadetblue}{rgb}{0.37,0.62,0.63}
\definecolor{chartreuse1}{rgb}{0.50,1.00,0.00}
\definecolor{chartreuse2}{rgb}{0.46,0.93,0.00}
\definecolor{chartreuse3}{rgb}{0.40,0.80,0.00}
\definecolor{chartreuse4}{rgb}{0.27,0.55,0.00}
\definecolor{chartreuse}{rgb}{0.50,1.00,0.00}
\definecolor{chocolate1}{rgb}{1.00,0.50,0.14}
\definecolor{chocolate2}{rgb}{0.93,0.46,0.13}
\definecolor{chocolate3}{rgb}{0.80,0.40,0.11}
\definecolor{chocolate4}{rgb}{0.55,0.27,0.07}
\definecolor{chocolate}{rgb}{0.82,0.41,0.12}
\definecolor{coral1}{rgb}{1.00,0.45,0.34}
\definecolor{coral2}{rgb}{0.93,0.42,0.31}
\definecolor{coral3}{rgb}{0.80,0.36,0.27}
\definecolor{coral4}{rgb}{0.55,0.24,0.18}
\definecolor{coral}{rgb}{1.00,0.50,0.31}
\definecolor{cornflowerblue}{rgb}{0.39,0.58,0.93}
\definecolor{cornsilk1}{rgb}{1.00,0.97,0.86}
\definecolor{cornsilk2}{rgb}{0.93,0.91,0.80}
\definecolor{cornsilk3}{rgb}{0.80,0.78,0.69}
\definecolor{cornsilk4}{rgb}{0.55,0.53,0.47}
\definecolor{cornsilk}{rgb}{1.00,0.97,0.86}
\definecolor{cyan1}{rgb}{0.00,1.00,1.00}
\definecolor{cyan2}{rgb}{0.00,0.93,0.93}
\definecolor{cyan3}{rgb}{0.00,0.80,0.80}
\definecolor{cyan4}{rgb}{0.00,0.55,0.55}
\definecolor{cyan}{rgb}{0.00,1.00,1.00}
\definecolor{darkblue}{rgb}{0.00,0.00,0.55}
\definecolor{darkcyan}{rgb}{0.00,0.55,0.55}
\definecolor{darkgoldenrod}{rgb}{0.72,0.53,0.04}
\definecolor{darkgray}{rgb}{0.66,0.66,0.66}
\definecolor{darkgreen}{rgb}{0.00,0.39,0.00}
\definecolor{darkgrey}{rgb}{0.66,0.66,0.66}
\definecolor{darkkhaki}{rgb}{0.74,0.72,0.42}
\definecolor{darkmagenta}{rgb}{0.55,0.00,0.55}
\definecolor{darkolive}{rgb}{0.33,0.42,0.18}
\definecolor{darkorange}{rgb}{1.00,0.55,0.00}
\definecolor{darkorchid}{rgb}{0.60,0.20,0.80}
\definecolor{darkred}{rgb}{0.55,0.00,0.00}
\definecolor{darksalmon}{rgb}{0.91,0.59,0.48}
\definecolor{darksea}{rgb}{0.56,0.74,0.56}
\definecolor{darkslate}{rgb}{0.18,0.31,0.31}
\definecolor{darkslate}{rgb}{0.18,0.31,0.31}
\definecolor{darkslate}{rgb}{0.28,0.24,0.55}
\definecolor{darkturquoise}{rgb}{0.00,0.81,0.82}
\definecolor{darkviolet}{rgb}{0.58,0.00,0.83}
\definecolor{deeppink}{rgb}{1.00,0.08,0.58}
\definecolor{deepsky}{rgb}{0.00,0.75,1.00}
\definecolor{dimgray}{rgb}{0.41,0.41,0.41}
\definecolor{dimgrey}{rgb}{0.41,0.41,0.41}
\definecolor{dodgerblue}{rgb}{0.12,0.56,1.00}
\definecolor{firebrick1}{rgb}{1.00,0.19,0.19}
\definecolor{firebrick2}{rgb}{0.93,0.17,0.17}
\definecolor{firebrick3}{rgb}{0.80,0.15,0.15}
\definecolor{firebrick4}{rgb}{0.55,0.10,0.10}
\definecolor{firebrick}{rgb}{0.70,0.13,0.13}
\definecolor{floralwhite}{rgb}{1.00,0.98,0.94}
\definecolor{forestgreen}{rgb}{0.13,0.55,0.13}
\definecolor{gainsboro}{rgb}{0.86,0.86,0.86}
\definecolor{ghostwhite}{rgb}{0.97,0.97,1.00}
\definecolor{gold1}{rgb}{1.00,0.84,0.00}
\definecolor{gold2}{rgb}{0.93,0.79,0.00}
\definecolor{gold3}{rgb}{0.80,0.68,0.00}
\definecolor{gold4}{rgb}{0.55,0.46,0.00}
\definecolor{goldenrod1}{rgb}{1.00,0.76,0.15}
\definecolor{goldenrod2}{rgb}{0.93,0.71,0.13}
\definecolor{goldenrod3}{rgb}{0.80,0.61,0.11}
\definecolor{goldenrod4}{rgb}{0.55,0.41,0.08}
\definecolor{goldenrod}{rgb}{0.85,0.65,0.13}
\definecolor{gold}{rgb}{1.00,0.84,0.00}
\definecolor{gray0}{rgb}{0.00,0.00,0.00}
\definecolor{gray100}{rgb}{1.00,1.00,1.00}
\definecolor{gray10}{rgb}{0.10,0.10,0.10}
\definecolor{gray11}{rgb}{0.11,0.11,0.11}
\definecolor{gray12}{rgb}{0.12,0.12,0.12}
\definecolor{gray13}{rgb}{0.13,0.13,0.13}
\definecolor{gray14}{rgb}{0.14,0.14,0.14}
\definecolor{gray15}{rgb}{0.15,0.15,0.15}
\definecolor{gray16}{rgb}{0.16,0.16,0.16}
\definecolor{gray17}{rgb}{0.17,0.17,0.17}
\definecolor{gray18}{rgb}{0.18,0.18,0.18}
\definecolor{gray19}{rgb}{0.19,0.19,0.19}
\definecolor{gray1}{rgb}{0.01,0.01,0.01}
\definecolor{gray20}{rgb}{0.20,0.20,0.20}
\definecolor{gray21}{rgb}{0.21,0.21,0.21}
\definecolor{gray22}{rgb}{0.22,0.22,0.22}
\definecolor{gray23}{rgb}{0.23,0.23,0.23}
\definecolor{gray24}{rgb}{0.24,0.24,0.24}
\definecolor{gray25}{rgb}{0.25,0.25,0.25}
\definecolor{gray26}{rgb}{0.26,0.26,0.26}
\definecolor{gray27}{rgb}{0.27,0.27,0.27}
\definecolor{gray28}{rgb}{0.28,0.28,0.28}
\definecolor{gray29}{rgb}{0.29,0.29,0.29}
\definecolor{gray2}{rgb}{0.02,0.02,0.02}
\definecolor{gray30}{rgb}{0.30,0.30,0.30}
\definecolor{gray31}{rgb}{0.31,0.31,0.31}
\definecolor{gray32}{rgb}{0.32,0.32,0.32}
\definecolor{gray33}{rgb}{0.33,0.33,0.33}
\definecolor{gray34}{rgb}{0.34,0.34,0.34}
\definecolor{gray35}{rgb}{0.35,0.35,0.35}
\definecolor{gray36}{rgb}{0.36,0.36,0.36}
\definecolor{gray37}{rgb}{0.37,0.37,0.37}
\definecolor{gray38}{rgb}{0.38,0.38,0.38}
\definecolor{gray39}{rgb}{0.39,0.39,0.39}
\definecolor{gray3}{rgb}{0.03,0.03,0.03}
\definecolor{gray40}{rgb}{0.40,0.40,0.40}
\definecolor{gray41}{rgb}{0.41,0.41,0.41}
\definecolor{gray42}{rgb}{0.42,0.42,0.42}
\definecolor{gray43}{rgb}{0.43,0.43,0.43}
\definecolor{gray44}{rgb}{0.44,0.44,0.44}
\definecolor{gray45}{rgb}{0.45,0.45,0.45}
\definecolor{gray46}{rgb}{0.46,0.46,0.46}
\definecolor{gray47}{rgb}{0.47,0.47,0.47}
\definecolor{gray48}{rgb}{0.48,0.48,0.48}
\definecolor{gray49}{rgb}{0.49,0.49,0.49}
\definecolor{gray4}{rgb}{0.04,0.04,0.04}
\definecolor{gray50}{rgb}{0.50,0.50,0.50}
\definecolor{gray51}{rgb}{0.51,0.51,0.51}
\definecolor{gray52}{rgb}{0.52,0.52,0.52}
\definecolor{gray53}{rgb}{0.53,0.53,0.53}
\definecolor{gray54}{rgb}{0.54,0.54,0.54}
\definecolor{gray55}{rgb}{0.55,0.55,0.55}
\definecolor{gray56}{rgb}{0.56,0.56,0.56}
\definecolor{gray57}{rgb}{0.57,0.57,0.57}
\definecolor{gray58}{rgb}{0.58,0.58,0.58}
\definecolor{gray59}{rgb}{0.59,0.59,0.59}
\definecolor{gray5}{rgb}{0.05,0.05,0.05}
\definecolor{gray60}{rgb}{0.60,0.60,0.60}
\definecolor{gray61}{rgb}{0.61,0.61,0.61}
\definecolor{gray62}{rgb}{0.62,0.62,0.62}
\definecolor{gray63}{rgb}{0.63,0.63,0.63}
\definecolor{gray64}{rgb}{0.64,0.64,0.64}
\definecolor{gray65}{rgb}{0.65,0.65,0.65}
\definecolor{gray66}{rgb}{0.66,0.66,0.66}
\definecolor{gray67}{rgb}{0.67,0.67,0.67}
\definecolor{gray68}{rgb}{0.68,0.68,0.68}
\definecolor{gray69}{rgb}{0.69,0.69,0.69}
\definecolor{gray6}{rgb}{0.06,0.06,0.06}
\definecolor{gray70}{rgb}{0.70,0.70,0.70}
\definecolor{gray71}{rgb}{0.71,0.71,0.71}
\definecolor{gray72}{rgb}{0.72,0.72,0.72}
\definecolor{gray73}{rgb}{0.73,0.73,0.73}
\definecolor{gray74}{rgb}{0.74,0.74,0.74}
\definecolor{gray75}{rgb}{0.75,0.75,0.75}
\definecolor{gray76}{rgb}{0.76,0.76,0.76}
\definecolor{gray77}{rgb}{0.77,0.77,0.77}
\definecolor{gray78}{rgb}{0.78,0.78,0.78}
\definecolor{gray79}{rgb}{0.79,0.79,0.79}
\definecolor{gray7}{rgb}{0.07,0.07,0.07}
\definecolor{gray80}{rgb}{0.80,0.80,0.80}
\definecolor{gray81}{rgb}{0.81,0.81,0.81}
\definecolor{gray82}{rgb}{0.82,0.82,0.82}
\definecolor{gray83}{rgb}{0.83,0.83,0.83}
\definecolor{gray84}{rgb}{0.84,0.84,0.84}
\definecolor{gray85}{rgb}{0.85,0.85,0.85}
\definecolor{gray86}{rgb}{0.86,0.86,0.86}
\definecolor{gray87}{rgb}{0.87,0.87,0.87}
\definecolor{gray88}{rgb}{0.88,0.88,0.88}
\definecolor{gray89}{rgb}{0.89,0.89,0.89}
\definecolor{gray8}{rgb}{0.08,0.08,0.08}
\definecolor{gray90}{rgb}{0.90,0.90,0.90}
\definecolor{gray91}{rgb}{0.91,0.91,0.91}
\definecolor{gray92}{rgb}{0.92,0.92,0.92}
\definecolor{gray93}{rgb}{0.93,0.93,0.93}
\definecolor{gray94}{rgb}{0.94,0.94,0.94}
\definecolor{gray95}{rgb}{0.95,0.95,0.95}
\definecolor{gray96}{rgb}{0.96,0.96,0.96}
\definecolor{gray97}{rgb}{0.97,0.97,0.97}
\definecolor{gray98}{rgb}{0.98,0.98,0.98}
\definecolor{gray99}{rgb}{0.99,0.99,0.99}
\definecolor{gray9}{rgb}{0.09,0.09,0.09}
\definecolor{gray}{rgb}{0.75,0.75,0.75}
\definecolor{green1}{rgb}{0.00,1.00,0.00}
\definecolor{green2}{rgb}{0.00,0.93,0.00}
\definecolor{green3}{rgb}{0.00,0.80,0.00}
\definecolor{green4}{rgb}{0.00,0.55,0.00}
\definecolor{greenyellow}{rgb}{0.68,1.00,0.18}
\definecolor{green}{rgb}{0.00,1.00,0.00}
\definecolor{grey0}{rgb}{0.00,0.00,0.00}
\definecolor{grey100}{rgb}{1.00,1.00,1.00}
\definecolor{grey10}{rgb}{0.10,0.10,0.10}
\definecolor{grey11}{rgb}{0.11,0.11,0.11}
\definecolor{grey12}{rgb}{0.12,0.12,0.12}
\definecolor{grey13}{rgb}{0.13,0.13,0.13}
\definecolor{grey14}{rgb}{0.14,0.14,0.14}
\definecolor{grey15}{rgb}{0.15,0.15,0.15}
\definecolor{grey16}{rgb}{0.16,0.16,0.16}
\definecolor{grey17}{rgb}{0.17,0.17,0.17}
\definecolor{grey18}{rgb}{0.18,0.18,0.18}
\definecolor{grey19}{rgb}{0.19,0.19,0.19}
\definecolor{grey1}{rgb}{0.01,0.01,0.01}
\definecolor{grey20}{rgb}{0.20,0.20,0.20}
\definecolor{grey21}{rgb}{0.21,0.21,0.21}
\definecolor{grey22}{rgb}{0.22,0.22,0.22}
\definecolor{grey23}{rgb}{0.23,0.23,0.23}
\definecolor{grey24}{rgb}{0.24,0.24,0.24}
\definecolor{grey25}{rgb}{0.25,0.25,0.25}
\definecolor{grey26}{rgb}{0.26,0.26,0.26}
\definecolor{grey27}{rgb}{0.27,0.27,0.27}
\definecolor{grey28}{rgb}{0.28,0.28,0.28}
\definecolor{grey29}{rgb}{0.29,0.29,0.29}
\definecolor{grey2}{rgb}{0.02,0.02,0.02}
\definecolor{grey30}{rgb}{0.30,0.30,0.30}
\definecolor{grey31}{rgb}{0.31,0.31,0.31}
\definecolor{grey32}{rgb}{0.32,0.32,0.32}
\definecolor{grey33}{rgb}{0.33,0.33,0.33}
\definecolor{grey34}{rgb}{0.34,0.34,0.34}
\definecolor{grey35}{rgb}{0.35,0.35,0.35}
\definecolor{grey36}{rgb}{0.36,0.36,0.36}
\definecolor{grey37}{rgb}{0.37,0.37,0.37}
\definecolor{grey38}{rgb}{0.38,0.38,0.38}
\definecolor{grey39}{rgb}{0.39,0.39,0.39}
\definecolor{grey3}{rgb}{0.03,0.03,0.03}
\definecolor{grey40}{rgb}{0.40,0.40,0.40}
\definecolor{grey41}{rgb}{0.41,0.41,0.41}
\definecolor{grey42}{rgb}{0.42,0.42,0.42}
\definecolor{grey43}{rgb}{0.43,0.43,0.43}
\definecolor{grey44}{rgb}{0.44,0.44,0.44}
\definecolor{grey45}{rgb}{0.45,0.45,0.45}
\definecolor{grey46}{rgb}{0.46,0.46,0.46}
\definecolor{grey47}{rgb}{0.47,0.47,0.47}
\definecolor{grey48}{rgb}{0.48,0.48,0.48}
\definecolor{grey49}{rgb}{0.49,0.49,0.49}
\definecolor{grey4}{rgb}{0.04,0.04,0.04}
\definecolor{grey50}{rgb}{0.50,0.50,0.50}
\definecolor{grey51}{rgb}{0.51,0.51,0.51}
\definecolor{grey52}{rgb}{0.52,0.52,0.52}
\definecolor{grey53}{rgb}{0.53,0.53,0.53}
\definecolor{grey54}{rgb}{0.54,0.54,0.54}
\definecolor{grey55}{rgb}{0.55,0.55,0.55}
\definecolor{grey56}{rgb}{0.56,0.56,0.56}
\definecolor{grey57}{rgb}{0.57,0.57,0.57}
\definecolor{grey58}{rgb}{0.58,0.58,0.58}
\definecolor{grey59}{rgb}{0.59,0.59,0.59}
\definecolor{grey5}{rgb}{0.05,0.05,0.05}
\definecolor{grey60}{rgb}{0.60,0.60,0.60}
\definecolor{grey61}{rgb}{0.61,0.61,0.61}
\definecolor{grey62}{rgb}{0.62,0.62,0.62}
\definecolor{grey63}{rgb}{0.63,0.63,0.63}
\definecolor{grey64}{rgb}{0.64,0.64,0.64}
\definecolor{grey65}{rgb}{0.65,0.65,0.65}
\definecolor{grey66}{rgb}{0.66,0.66,0.66}
\definecolor{grey67}{rgb}{0.67,0.67,0.67}
\definecolor{grey68}{rgb}{0.68,0.68,0.68}
\definecolor{grey69}{rgb}{0.69,0.69,0.69}
\definecolor{grey6}{rgb}{0.06,0.06,0.06}
\definecolor{grey70}{rgb}{0.70,0.70,0.70}
\definecolor{grey71}{rgb}{0.71,0.71,0.71}
\definecolor{grey72}{rgb}{0.72,0.72,0.72}
\definecolor{grey73}{rgb}{0.73,0.73,0.73}
\definecolor{grey74}{rgb}{0.74,0.74,0.74}
\definecolor{grey75}{rgb}{0.75,0.75,0.75}
\definecolor{grey76}{rgb}{0.76,0.76,0.76}
\definecolor{grey77}{rgb}{0.77,0.77,0.77}
\definecolor{grey78}{rgb}{0.78,0.78,0.78}
\definecolor{grey79}{rgb}{0.79,0.79,0.79}
\definecolor{grey7}{rgb}{0.07,0.07,0.07}
\definecolor{grey80}{rgb}{0.80,0.80,0.80}
\definecolor{grey81}{rgb}{0.81,0.81,0.81}
\definecolor{grey82}{rgb}{0.82,0.82,0.82}
\definecolor{grey83}{rgb}{0.83,0.83,0.83}
\definecolor{grey84}{rgb}{0.84,0.84,0.84}
\definecolor{grey85}{rgb}{0.85,0.85,0.85}
\definecolor{grey86}{rgb}{0.86,0.86,0.86}
\definecolor{grey87}{rgb}{0.87,0.87,0.87}
\definecolor{grey88}{rgb}{0.88,0.88,0.88}
\definecolor{grey89}{rgb}{0.89,0.89,0.89}
\definecolor{grey8}{rgb}{0.08,0.08,0.08}
\definecolor{grey90}{rgb}{0.90,0.90,0.90}
\definecolor{grey91}{rgb}{0.91,0.91,0.91}
\definecolor{grey92}{rgb}{0.92,0.92,0.92}
\definecolor{grey93}{rgb}{0.93,0.93,0.93}
\definecolor{grey94}{rgb}{0.94,0.94,0.94}
\definecolor{grey95}{rgb}{0.95,0.95,0.95}
\definecolor{grey96}{rgb}{0.96,0.96,0.96}
\definecolor{grey97}{rgb}{0.97,0.97,0.97}
\definecolor{grey98}{rgb}{0.98,0.98,0.98}
\definecolor{grey99}{rgb}{0.99,0.99,0.99}
\definecolor{grey9}{rgb}{0.09,0.09,0.09}
\definecolor{grey}{rgb}{0.75,0.75,0.75}
\definecolor{honeydew1}{rgb}{0.94,1.00,0.94}
\definecolor{honeydew2}{rgb}{0.88,0.93,0.88}
\definecolor{honeydew3}{rgb}{0.76,0.80,0.76}
\definecolor{honeydew4}{rgb}{0.51,0.55,0.51}
\definecolor{honeydew}{rgb}{0.94,1.00,0.94}
\definecolor{hotpink}{rgb}{1.00,0.41,0.71}
\definecolor{indianred}{rgb}{0.80,0.36,0.36}
\definecolor{ivory1}{rgb}{1.00,1.00,0.94}
\definecolor{ivory2}{rgb}{0.93,0.93,0.88}
\definecolor{ivory3}{rgb}{0.80,0.80,0.76}
\definecolor{ivory4}{rgb}{0.55,0.55,0.51}
\definecolor{ivory}{rgb}{1.00,1.00,0.94}
\definecolor{khaki1}{rgb}{1.00,0.96,0.56}
\definecolor{khaki2}{rgb}{0.93,0.90,0.52}
\definecolor{khaki3}{rgb}{0.80,0.78,0.45}
\definecolor{khaki4}{rgb}{0.55,0.53,0.31}
\definecolor{khaki}{rgb}{0.94,0.90,0.55}
\definecolor{lavenderblush}{rgb}{1.00,0.94,0.96}
\definecolor{lavender}{rgb}{0.90,0.90,0.98}
\definecolor{lawngreen}{rgb}{0.49,0.99,0.00}
\definecolor{lemonchiffon}{rgb}{1.00,0.98,0.80}
\definecolor{lightblue}{rgb}{0.68,0.85,0.90}
\definecolor{lightcoral}{rgb}{0.94,0.50,0.50}
\definecolor{lightcyan}{rgb}{0.88,1.00,1.00}
\definecolor{lightgoldenrod}{rgb}{0.93,0.87,0.51}
\definecolor{lightgoldenrod}{rgb}{0.98,0.98,0.82}
\definecolor{lightgray}{rgb}{0.83,0.83,0.83}
\definecolor{lightgreen}{rgb}{0.56,0.93,0.56}
\definecolor{lightgrey}{rgb}{0.83,0.83,0.83}
\definecolor{lightpink}{rgb}{1.00,0.71,0.76}
\definecolor{lightsalmon}{rgb}{1.00,0.63,0.48}
\definecolor{lightsea}{rgb}{0.13,0.70,0.67}
\definecolor{lightsky}{rgb}{0.53,0.81,0.98}
\definecolor{lightslate}{rgb}{0.47,0.53,0.60}
\definecolor{lightslate}{rgb}{0.47,0.53,0.60}
\definecolor{lightslate}{rgb}{0.52,0.44,1.00}
\definecolor{lightsteel}{rgb}{0.69,0.77,0.87}
\definecolor{lightyellow}{rgb}{1.00,1.00,0.88}
\definecolor{limegreen}{rgb}{0.20,0.80,0.20}
\definecolor{linen}{rgb}{0.98,0.94,0.90}
\definecolor{magenta1}{rgb}{1.00,0.00,1.00}
\definecolor{magenta2}{rgb}{0.93,0.00,0.93}
\definecolor{magenta3}{rgb}{0.80,0.00,0.80}
\definecolor{magenta4}{rgb}{0.55,0.00,0.55}
\definecolor{magenta}{rgb}{1.00,0.00,1.00}
\definecolor{maroon1}{rgb}{1.00,0.20,0.70}
\definecolor{maroon2}{rgb}{0.93,0.19,0.65}
\definecolor{maroon3}{rgb}{0.80,0.16,0.56}
\definecolor{maroon4}{rgb}{0.55,0.11,0.38}
\definecolor{maroon}{rgb}{0.69,0.19,0.38}
\definecolor{mediumaquamarine}{rgb}{0.40,0.80,0.67}
\definecolor{mediumblue}{rgb}{0.00,0.00,0.80}
\definecolor{mediumorchid}{rgb}{0.73,0.33,0.83}
\definecolor{mediumpurple}{rgb}{0.58,0.44,0.86}
\definecolor{mediumsea}{rgb}{0.24,0.70,0.44}
\definecolor{mediumslate}{rgb}{0.48,0.41,0.93}
\definecolor{mediumspring}{rgb}{0.00,0.98,0.60}
\definecolor{mediumturquoise}{rgb}{0.28,0.82,0.80}
\definecolor{mediumviolet}{rgb}{0.78,0.08,0.52}
\definecolor{midnightblue}{rgb}{0.10,0.10,0.44}
\definecolor{mintcream}{rgb}{0.96,1.00,0.98}
\definecolor{mistyrose}{rgb}{1.00,0.89,0.88}
\definecolor{moccasin}{rgb}{1.00,0.89,0.71}
\definecolor{navajowhite}{rgb}{1.00,0.87,0.68}
\definecolor{navyblue}{rgb}{0.00,0.00,0.50}
\definecolor{navy}{rgb}{0.00,0.00,0.50}
\definecolor{oldlace}{rgb}{0.99,0.96,0.90}
\definecolor{olivedrab}{rgb}{0.42,0.56,0.14}
\definecolor{orange1}{rgb}{1.00,0.65,0.00}
\definecolor{orange2}{rgb}{0.93,0.60,0.00}
\definecolor{orange3}{rgb}{0.80,0.52,0.00}
\definecolor{orange4}{rgb}{0.55,0.35,0.00}
\definecolor{orangered}{rgb}{1.00,0.27,0.00}
\definecolor{orange}{rgb}{1.00,0.65,0.00}
\definecolor{orchid1}{rgb}{1.00,0.51,0.98}
\definecolor{orchid2}{rgb}{0.93,0.48,0.91}
\definecolor{orchid3}{rgb}{0.80,0.41,0.79}
\definecolor{orchid4}{rgb}{0.55,0.28,0.54}
\definecolor{orchid}{rgb}{0.85,0.44,0.84}
\definecolor{palegoldenrod}{rgb}{0.93,0.91,0.67}
\definecolor{palegreen}{rgb}{0.60,0.98,0.60}
\definecolor{paleturquoise}{rgb}{0.69,0.93,0.93}
\definecolor{paleviolet}{rgb}{0.86,0.44,0.58}
\definecolor{papayawhip}{rgb}{1.00,0.94,0.84}
\definecolor{peachpuff}{rgb}{1.00,0.85,0.73}
\definecolor{peru}{rgb}{0.80,0.52,0.25}
\definecolor{pink1}{rgb}{1.00,0.71,0.77}
\definecolor{pink2}{rgb}{0.93,0.66,0.72}
\definecolor{pink3}{rgb}{0.80,0.57,0.62}
\definecolor{pink4}{rgb}{0.55,0.39,0.42}
\definecolor{pink}{rgb}{1.00,0.75,0.80}
\definecolor{plum1}{rgb}{1.00,0.73,1.00}
\definecolor{plum2}{rgb}{0.93,0.68,0.93}
\definecolor{plum3}{rgb}{0.80,0.59,0.80}
\definecolor{plum4}{rgb}{0.55,0.40,0.55}
\definecolor{plum}{rgb}{0.87,0.63,0.87}
\definecolor{powderblue}{rgb}{0.69,0.88,0.90}
\definecolor{purple1}{rgb}{0.61,0.19,1.00}
\definecolor{purple2}{rgb}{0.57,0.17,0.93}
\definecolor{purple3}{rgb}{0.49,0.15,0.80}
\definecolor{purple4}{rgb}{0.33,0.10,0.55}
\definecolor{purple}{rgb}{0.63,0.13,0.94}
\definecolor{red1}{rgb}{1.00,0.00,0.00}
\definecolor{red2}{rgb}{0.93,0.00,0.00}
\definecolor{red3}{rgb}{0.80,0.00,0.00}
\definecolor{red4}{rgb}{0.55,0.00,0.00}
\definecolor{red}{rgb}{1.00,0.00,0.00}
\definecolor{rosybrown}{rgb}{0.74,0.56,0.56}
\definecolor{royalblue}{rgb}{0.25,0.41,0.88}
\definecolor{saddlebrown}{rgb}{0.55,0.27,0.07}
\definecolor{salmon1}{rgb}{1.00,0.55,0.41}
\definecolor{salmon2}{rgb}{0.93,0.51,0.38}
\definecolor{salmon3}{rgb}{0.80,0.44,0.33}
\definecolor{salmon4}{rgb}{0.55,0.30,0.22}
\definecolor{salmon}{rgb}{0.98,0.50,0.45}
\definecolor{sandybrown}{rgb}{0.96,0.64,0.38}
\definecolor{seagreen}{rgb}{0.18,0.55,0.34}
\definecolor{seashell1}{rgb}{1.00,0.96,0.93}
\definecolor{seashell2}{rgb}{0.93,0.90,0.87}
\definecolor{seashell3}{rgb}{0.80,0.77,0.75}
\definecolor{seashell4}{rgb}{0.55,0.53,0.51}
\definecolor{seashell}{rgb}{1.00,0.96,0.93}
\definecolor{sienna1}{rgb}{1.00,0.51,0.28}
\definecolor{sienna2}{rgb}{0.93,0.47,0.26}
\definecolor{sienna3}{rgb}{0.80,0.41,0.22}
\definecolor{sienna4}{rgb}{0.55,0.28,0.15}
\definecolor{sienna}{rgb}{0.63,0.32,0.18}
\definecolor{skyblue}{rgb}{0.53,0.81,0.92}
\definecolor{slateblue}{rgb}{0.42,0.35,0.80}
\definecolor{slategray}{rgb}{0.44,0.50,0.56}
\definecolor{slategrey}{rgb}{0.44,0.50,0.56}
\definecolor{snow1}{rgb}{1.00,0.98,0.98}
\definecolor{snow2}{rgb}{0.93,0.91,0.91}
\definecolor{snow3}{rgb}{0.80,0.79,0.79}
\definecolor{snow4}{rgb}{0.55,0.54,0.54}
\definecolor{snow}{rgb}{1.00,0.98,0.98}
\definecolor{springgreen}{rgb}{0.00,1.00,0.50}
\definecolor{steelblue}{rgb}{0.27,0.51,0.71}
\definecolor{tan1}{rgb}{1.00,0.65,0.31}
\definecolor{tan2}{rgb}{0.93,0.60,0.29}
\definecolor{tan3}{rgb}{0.80,0.52,0.25}
\definecolor{tan4}{rgb}{0.55,0.35,0.17}
\definecolor{tan}{rgb}{0.82,0.71,0.55}
\definecolor{thistle1}{rgb}{1.00,0.88,1.00}
\definecolor{thistle2}{rgb}{0.93,0.82,0.93}
\definecolor{thistle3}{rgb}{0.80,0.71,0.80}
\definecolor{thistle4}{rgb}{0.55,0.48,0.55}
\definecolor{thistle}{rgb}{0.85,0.75,0.85}
\definecolor{tomato1}{rgb}{1.00,0.39,0.28}
\definecolor{tomato2}{rgb}{0.93,0.36,0.26}
\definecolor{tomato3}{rgb}{0.80,0.31,0.22}
\definecolor{tomato4}{rgb}{0.55,0.21,0.15}
\definecolor{tomato}{rgb}{1.00,0.39,0.28}
\definecolor{turquoise1}{rgb}{0.00,0.96,1.00}
\definecolor{turquoise2}{rgb}{0.00,0.90,0.93}
\definecolor{turquoise3}{rgb}{0.00,0.77,0.80}
\definecolor{turquoise4}{rgb}{0.00,0.53,0.55}
\definecolor{turquoise}{rgb}{0.25,0.88,0.82}
\definecolor{violetred}{rgb}{0.82,0.13,0.56}
\definecolor{violet}{rgb}{0.93,0.51,0.93}
\definecolor{wheat1}{rgb}{1.00,0.91,0.73}
\definecolor{wheat2}{rgb}{0.93,0.85,0.68}
\definecolor{wheat3}{rgb}{0.80,0.73,0.59}
\definecolor{wheat4}{rgb}{0.55,0.49,0.40}
\definecolor{wheat}{rgb}{0.96,0.87,0.70}
\definecolor{whitesmoke}{rgb}{0.96,0.96,0.96}
\definecolor{white}{rgb}{1.00,1.00,1.00}
\definecolor{yellow1}{rgb}{1.00,1.00,0.00}
\definecolor{yellow2}{rgb}{0.93,0.93,0.00}
\definecolor{yellow3}{rgb}{0.80,0.80,0.00}
\definecolor{yellow4}{rgb}{0.55,0.55,0.00}
\definecolor{yellowgreen}{rgb}{0.60,0.80,0.20}
\definecolor{yellow}{rgb}{1.00,1.00,0.00}
 
\usepackage[usenames,dvipsnames]{xcolor}
\usepackage[psamsfonts]{amsfonts}
\usepackage{graphicx}
\usepackage{amssymb}
\usepackage{shortvrb} 
\usepackage{mathtools}  

\usepackage{hyperref} 
\hypersetup{pdfpagemode=FullScreen,   
backref,
colorlinks=true, 
citecolor=Bittersweet,
linkcolor=Bittersweet, 
urlcolor=Bittersweet
}
 
\setattribute{journal}{name}{}
\startlocaldefs

\theoremstyle{plain}
\newtheorem{Lem}{Lemma}[section]

\newtheorem{Theor}{Theorem}[section] 

\newtheorem{Defin}{Definition}[section]

\numberwithin{equation}{section}

\newcommand{\cqfd}{\hfill $\square$}
\newcommand{\R}{\mathbb R}

\newcommand{\ub}{{u}}
\newcommand{\zb}{{z}}
\newcommand{\yb}{{y}}
\newcommand{\xb}{{x}}
\newcommand{\eb}{{e}}
\newcommand{\vb}{{v}}
\newcommand{\qb}{{q}}
\newcommand{\tb}{{t}}
\newcommand{\bb}{{b}}
\newcommand{\Xb}{{X}}
\newcommand{\Yb}{{Y}}
\newcommand{\Zb}{{Z}}
\newcommand{\Ab}{{A}}
\newcommand{\Ob}{{O}}

\newcommand{\thetab}{{\theta}}
\newcommand{\mub}{{\mu}}
\newcommand{\Sigmab}{{\Sigma}}

\endlocaldefs

\begin{document}

\begin{frontmatter}
\title{From Halfspace M-depth to Multiple-output Expectile Regression}
\runtitle{From \mbox{M-depth} to Multiple-output Expectile Regression}

\begin{aug}
\author{\fnms{Abdelaati} \snm{Daouia}\thanksref{t1}\ead[label=e1]{abdelaati.daouia@tse-fr.eu}
\ead[label=u1,url]{https://www.tse-fr.eu/people/abdelaati-daouia}}
\,\! \and\,\hspace{1pt} 
\author{\fnms{Davy} \snm{Paindaveine}\thanksref{t2}
\ead[label=e2]{dpaindav@ulb.ac.be}
\ead[label=u2,url]{http://homepages.ulb.ac.be/dpaindav}}

\thankstext{t1}{Abdelaati Daouia's research is supported by the Toulouse School of Economics Individual Research Fund.}
\thankstext{t2}{Corresponding author. Davy Paindaveine's research is supported by a research fellowship from the Francqui Foundation and by the Program of Concerted Research Actions (ARC) of the Universit\'{e} libre de Bruxelles.}

\runauthor{A. Daouia and D. Paindaveine}

\affiliation{Universit\'{e} Toulouse Capitole and Universit\'{e} libre de Bruxelles}

\address{Universit\'{e} Toulouse Capitole\\
Manufacture des Tabacs\\
21, All\'{e}e de Brienne\\ 
31015 Toulouse Cedex 6\\ 
France\\
\printead{e1}\\  
%\printead{u1}\\
}

\address{Universit\'{e} libre de Bruxelles\\
ECARES and D\'{e}partement de Math\'{e}matique\\
Boulevard du Triomphe, CP210\\
B-1050, Brussels\\
Belgium\\ 
\printead{e2}\\ 
\printead{u2}\\
}

\end{aug}
\vspace{3mm}

\begin{abstract}
Despite the renewed interest in the \cite{NewPow1987} concept of \emph{expectiles} in fields such as econometrics, risk management, and extreme value theory, expectile regression---or, more generally, M-quantile regression---unfortunately remains limited to single-output problems. To improve on this, we introduce hyperplane-valued multivariate \mbox{M-quantiles} that show strong advantages, for instance in terms of equivariance, over the various point-valued multivariate \mbox{M-quantiles} available in the literature. Like their competitors, our multivariate M-quantiles are directional in nature and provide centrality regions when all directions are considered. These regions define a new statistical depth, the \emph{halfspace M-depth}, whose deepest point, in the expectile case, is the mean vector. Remarkably, the halfspace \mbox{M-depth} can alternatively be obtained by substituting, in the celebrated \cite{Tuk1975} halfspace depth, M-quantile outlyingness for standard quantile outlyingness, which supports a posteriori the claim that our multivariate M-quantile concept is the natural one. We investigate thoroughly the properties of the proposed multivariate \mbox{M-quantiles}, of halfspace \mbox{M-depth}, and of the corresponding regions. Since our original motivation was to define multiple-output expectile regression methods, we further focus on the expectile case. We show in particular that expectile depth is smoother than the Tukey depth and enjoys interesting monotonicity properties that are extremely promising for computational purposes. Unlike their quantile analogs, the proposed multivariate expectiles also satisfy the coherency axioms of multivariate risk measures. Finally, we show that our multivariate expectiles indeed allow performing multiple-output expectile regression, which is illustrated on simulated and real data. 
\end{abstract}



\end{frontmatter}

%%%%%%%%%%%%%%%%%%%%

\section{Introduction} 
\label{introsec}

Whenever one wants to assess the impact of a vector of covariates~$X$ on a scalar response~$Y$, mean regression, in its various forms (linear, nonlinear, or nonparametric), remains by far the most popular method. Mean regression, however, only captures the conditional mean
$$
\mu(x)
:=
{\rm E}[ Y | X=x]
=\arg\min_{\theta\in\mathbb{R}} {\rm E}\big[ (Y-\theta)^2
%-Y^2 
| X=x \big]
$$
of the response, hence fails to describe thoroughly the conditional distribution of~$Y$ given~$X$. Such a thorough description is given by the \cite{KoeBas1978} \emph{quantile regression}, that considers the conditional quantiles
\begin{equation}
\label{quantintro}
q_\alpha(x):=\arg\min_{\theta\in\mathbb{R}} {\rm E}\big[ \rho_{\alpha,L_1} (Y-\theta)
%-\rho_{\alpha,L_1}(Y) 
| X=x \big],
\quad
\alpha\in(0,1)
,
\end{equation}
where
$
%\begin{equation}
%	\label{defrho1andh}
\rho_{\alpha,L_1}(t)
:=
((1-\alpha) \mathbb{I}[t<0]
+\alpha \mathbb{I}[t>0]
) 
|t|
$
is the \emph{check function} (throughout, $\mathbb{I}[A]$ stands for the indicator function of~$A$). 
An alternative to quantile regression is the \cite{NewPow1987}
%; see also \cite{Aigetal1976}. 
\emph{expectile regression}, that focuses on the conditional expectiles
\begin{equation}
\label{expectintro}
e_\alpha(x):=\arg\min_{\theta\in\mathbb{R}} {\rm E}\big[ \rho_{\alpha,L_2}(Y-\theta)
| X=x \big],
\quad
\alpha\in(0,1)
,
\end{equation}
where 
$
%\begin{equation}
%	\label{defrho1andh}
\rho_{\alpha,L_2}(t)
:=
((1-\alpha) \mathbb{I}[t<0]
+\alpha \mathbb{I}[t>0]
) 
t^2
$
is an asymmetric quadratic loss function, in the same way the check function is an asymmetric absolute loss function. Conditional expectiles, like conditional quantiles, fully characterize the conditional distribution of the response and nicely include the conditional mean~$\mu(x)$ as a particular case. Sample conditional expectiles, unlike their quantile counterparts, are sensitive to extreme observations, but this may actually be an asset in some applications; in financial risk management, for instance, quantiles are often criticized for being too liberal (due to their insensitivity to extreme losses) and expectiles are therefore favoured in any prudent and reactive risk analysis (\citealp{Daoetal2018}).

Expectile regression shows other advantages over quantile regression, of which we mention only a few here. 
First, inference on quantiles requires estimating nonparametrically the conditional density of the response at the considered quantiles, which is notoriously difficult. In constrast, inference on expectiles can be performed without resorting to any smoothing technique, which makes it easy, e.g., to test for homoscedasticity or for conditional symmetry in linear regression models (\citealp{NewPow1987}). 
Second, since expectile regression includes classical mean regression as a particular case, it is closer to the least squares notion of explained variance and, in parametric cases, expectile regression coefficients can be interpreted with respect to variance heteroscedasticity. This is of particular relevance in complex regression specifications including nonlinear, random or spatial effects (\citealp{SobKne2012}).
Third, expectile smoothing techniques, based on kernel smoothing (\citealp{YaoTon1996}) or penalized splines (\citealp{SchEil2009}), show better smoothness and stability than their quantile counterparts and also make expectile crossings far more rare than quantile crossings;
see \cite{SchEil2009}, \cite{Eil2013} and \cite{Schetal2015}. 
%
%Finally, as pointed out in \cite{Efr1991}, 
%%\cite{YaoTon1996}, 
%\cite{SchEil2013} and \cite{Schetal2015}, sample quantiles and their intuitive appeal can actually be easily recovered from a set of sample expectiles.
% 
These points explain why expectiles recently regained much interest in econometrics; see, e.g., \cite{Kuaetal2009}, \cite{DeRHar2009}, and \cite{EmbHof2014}. 

Despite these nice properties, expectile regression still suffers from an important drawback, namely its limitation to single-output problems. %To the best of the authors' knowledge, indeed, no concept of multiple-output expectile regression has been proposed in the literature. 
In contrast, many works developed multiple-output quantile regression methods. We refer, among others, to 
\cite{Cha2003}, 
\cite{CheDeG2007}, 
\cite{Wei2008}, 
\cite{Haletal2010}, 
%\cite{PaiSim2011}, 
%\cite{McKetal2011},  
%\cite{KonMiz2012}, 
%\cite{ChaLai2013},
\cite{CouDiB2013},
\cite{WalKne2015}, 
\cite{Haletal2015},  
\cite{Caretal2016,Caretal2017}, 
%\cite{}, 
and 
\cite{Cha2017}. This is in line with the fact that defining a satisfactory concept of multivariate quantile is a classical problem that has attracted much attention in the literature (we refer to \cite{Ser2002} and to the references therein), whereas the literature on multivariate expectiles is much sparser. Some early efforts to define multivariate expectiles can be found in  \cite{Kol1997}, \cite{Breetal2001} and \cite{Koketal2002}, that all define more generally multivariate versions of the \emph{M-quantiles} from \cite{BreCha1988} (a first concept of multivariate M-quantile was actually already discussed in \cite{BreCha1988} itself). Recently, there has been a renewed interest in defining %specifically 
multivariate expectiles; we refer to \cite{CouDiB2014}, \cite{Mauetal2018,Mauetal2017}, and to \cite{Heretal2018}. Multivariate risk handling in finance and actuarial sciences is mostly behind this growing interest, as will be discussed in Section~\ref{secrisks} below.

This paper introduces multivariate expectiles---and, more generally, multivariate M-quantiles---that enjoy many desirable properties, particularly in terms of affine equivariance. While this equivariance property is a standard requirement in the companion problem of defining multivariate quantiles, the available concepts of multivariate expectiles or M-quantiles are at best orthogonal-equivariant. Like their competitors, our multivariate \mbox{M-quantiles} are directional quantities, but they are hyperplane-valued rather than point-valued. Despite this different nature, they still generate centrality regions when all directions are considered. While this has not been discussed in the multivariate M-quantile literature (nor in the multivariate expectile one), this defines an M-concept of statistical depth. The resulting \emph{halfspace \mbox{M-depth}} generalizes the \cite{Tuk1975} halfspace depth and satisfies the desirable properties of depth identified in \cite{ZuoSer2000A}. Remarkably, this \mbox{M-depth} can alternatively be obtained by replacing, in the halfspace Tukey depth, standard quantile outlyingness with M-quantile outlyingness, which a posteriori supports the claim that our multivariate M-quantile concept is the natural one. This is a key result that allows us to study the structural properties of \mbox{M-depth}. Compared to Tukey depth, the particular case of expectile depth shows interesting properties in terms of, e.g., smoothness and monotonicity, and should be appealing to practitioners due to its link with the most classical location functional, namely the mean vector. Our multivariate expectiles, unlike their quantile counterparts, also satisfy all axioms of coherent risk measures. Finally, in line with our original objective, they allow us to define multiple-output expectile regression methods that actually show strong advantages over their quantile counterparts. 
 
The outline of the paper is as follows. In Section~\ref{secMquantile}, we carefully define univariate M-quantiles through a theorem that extends a result from \cite{Jon1994} and is of independent interest.  
In Section~\ref{secproposedmultiquantile}, we introduce our concept of multivariate M-quantiles and compare the resulting M-quantile regions with those associated with alternative M-quantile concepts.  
In Section~\ref{secMdepth}, we define halfspace \mbox{M-depth} and investigate its properties, whereas, in Section~\ref{secExpectiledepth}, we focus on the particular case of expectile depth.
In Section~\ref{secrisks}, we discuss the relation between multivariate M-quantiles and risk measures, and we show that our expectiles satisfy the coherency axioms of multivariate risk measures.   
In Section~\ref{secregres}, we explain how these expectiles allow performing multiple-output expectile regression, which is illustrated on simulated and real data.
Final comments and perspectives for future research are provided in Section~\ref{secfinalcom}. 
Appendix~\ref{secAppA} describes some of the main competing multivariate \mbox{M-quantile} concepts, whereas Appendix~\ref{secAppB} collects all proofs.

%%%%%%%%%%%%%%%%%%%%%%%%%%%%%%

\section{On univariate M-quantiles} 
\label{secMquantile}

As mentioned in \cite{Jon1994}, the \cite{BreCha1988} M-quantiles are related to M-estimates (or \mbox{M-functionals}) of location in the same way standard quantiles are related to the median. In line with~(\ref{quantintro})-(\ref{expectintro}), the \mbox{order-$\alpha$} M-quantile of a probability measure~$P$ over~$\R$ may be thought of as
\begin{equation}
\label{Mquantfirstdef}	
\theta_\alpha^\rho(P):=
\arg\min_{\theta\in\mathbb{R}} {\rm E}\big[ \rho_{\alpha}(Z- \theta)
%-\rho_{\alpha}(Z) 
\big]
%,
%\quad
%\alpha\in(0,1)
,
\end{equation}
where~$\rho_\alpha(t)
%:= h_\alpha(t) \rho(t)
:=
(
(1-\alpha) \mathbb{I}[t<0]
%+ {\textstyle{\frac{1}{2}}} \mathbb{I}[t=0]
+\alpha \mathbb{I}[t>0]
)
\rho(t)$ is based on a suitable symmetric loss function~$\rho$ and where  the random variable~$Z$ has distribution~$P$.
% satisfying~$\rho(0)=0$.
% (the term involving~$\mathbb{I}[t=0]$ in the function~$h_\alpha$ could be safely dropped here, but some later claims in this section will require that $h_\alpha$ is defined as above).  
%
Standard quantiles are obtained for the absolute loss function~$\rho(t)=|t|$, whereas expectiles are associated with the quadratic loss function~$\rho(t)=t^2$. One may also consider the Huber loss functions
\begin{equation}
\label{rhoc}
%t\mapsto 
\rho_c(t)
:=
\frac{t^2}{2c} \, \mathbb{I}[|t|< c] +
\Big(|t|-\frac{c}{2}\Big)\mathbb{I}[|t|\geq c] 
,
\qquad 
c>0
, 
\end{equation}
that allow recovering, up to an irrelevant positive scalar factor, the absolute value and quadratic loss functions above. The resulting M-quantiles~$\theta_\alpha^{\rho_c}(P)$ thus offer a continuum between quantiles and expectiles. 

The M-quantiles in~(\ref{Mquantfirstdef}) may be non-unique: for instance, if~$\rho(t)=|t|$ and~$P=P_n$ is the empirical probability measure associated with a sample of size~$n$, then~$\theta_\alpha^\rho(P)$, for any~$\alpha=1/n,\ldots,(n-1)/n$, is an interval with non-empty interior. Another issue is that it is unclear what is the collection of probability measures~$P$ for which~$\theta_\alpha^\rho(P)$ is well-defined. We will therefore adopt an alternative definition of M-quantiles, that results from Theorem~\ref{Theorlemjones} below. The result, that is of independent interest, significantly extends the theorem in page~151 of \cite{Jon1994} (in particular, Jones' result excludes the absolute loss function and all Huber loss functions). To state the result, we define the class~$\mathcal{C}$ of loss functions~$\rho:\R\to\R^+$ that are convex, symmetric and such that~$\rho(t)=0$ for~$t=0$ only. For any~$\rho\in \mathcal{C}$, we write~$\psi_-$ for the left-derivative of~$\rho$ (existence follows from convexity of~$\rho$) and we denote as~$\mathcal{P}^{\rho}$ the collection of probability measures~$P$ over~$\R$ such that (i)~$P[\{\theta\}]<1$ for any~$\theta\in\R$ and (ii)~$\int_{-\infty}^\infty |\psi_-(z-\theta)| \,dP(z)<\infty$ for any~$\theta\in\R$. 

\begin{Theor}
\label{Theorlemjones}
Fix~$\alpha\in(0,1)$, $\rho\in\mathcal{C}$ and~$P\in\mathcal{P}^{\rho}$. Let~$Z$ be a random variable with distribution~$P$. Then, (i) $\theta
\mapsto O^\rho(\theta) 
:=
{\rm E}[ 
\rho_\alpha(Z- \theta) -\rho_{\alpha}(Z)
]
$
is well-defined for any~$\theta$, and it is left- and right-differentiable over~$\R$, hence also continuous over~$\R$. (ii) The corresponding left- and right-derivatives 
%~$O_-^{\rho\prime}(\theta)$ and right-derivative~$O_+^{\rho\prime}(\theta)$ 
satisfy~$O_-^{\rho\prime}(\theta)\leq O_+^{\rho\prime}(\theta)$ at any~$\theta$. (iii) The sign of~$O_+^{\rho\prime}(\theta)$ is the same as that of~$G^\rho(\theta)-\alpha$, where we let 
%\begin{equation}
%\label{defg}
$$
G^\rho(\theta)
:=
\frac{{\rm E}[|\psi_-(Z-\theta)| \mathbb{I}[Z\leq \theta]]}{{\rm E}[|\psi_-(Z-\theta)|]}
;
$$
%\end{equation}
%
(iv)
$\theta \mapsto G^\rho(\theta)$ is a cumulative distribution function over~$\R$.
(v) The \emph{order-$\alpha$ \mbox{M-quantile} of~$P$}, which we define as
\begin{equation}
\label{Mquantdef}
%\theta^\rho_{\alpha}
%:=
\theta^\rho_{\alpha}(P)
:=
\inf\big\{\theta\in\R: G^\rho(\theta)\geq\alpha\big\}
,
\end{equation}
minimizes~$\theta\mapsto O^\rho(\theta)$ over~$\R$, hence provides a unique representative of the argmin in~(\ref{Mquantfirstdef}). (vi) If~$\psi_-$ is continuous over~$\R$ (or if~$P$ is non-atomic), then~$G^{\rho}$ is continuous at~$\theta^\rho_{\alpha}(P)$, so that~$G^{\rho}(\theta^\rho_{\alpha}(P))=\alpha$.   
\end{Theor}
\vspace{2mm}
 
In this result, the objective function in~(\ref{Mquantfirstdef}) was replaced by the modified one~$O^\rho(\theta)$, which does not have any impact on the corresponding argmin. This is a classical trick from the quantile and expectile literatures that ensures that quantiles (resp., expectiles) are well-defined without any moment condition (resp., under a finite first-order moment condition), whereas the corresponding original objective function in~(\ref{Mquantfirstdef}) in principle imposes a finite first-order moment (resp., a finite second-order moment).  
 
When case~(vi) applies (as it does for the quadratic loss function and any Huber loss function), the equation~$G^{\rho}(\theta^\rho_{\alpha}(P))=\alpha$ plays the role of the first-order condition associated with~(\ref{Mquantfirstdef}). When case~(vi) does not apply (which is the case for empirical probability measures when using the absolute loss function), this first-order condition is to be replaced by the more general one in~(\ref{Mquantdef}). For the absolute and quadratic loss functions, one has
$$
G^\rho(\theta)
=
P[Z\leq \theta]
\quad
\textrm{ and }
\quad
G^\rho(\theta)
=
\frac{{\rm E}[|Z-\theta| \mathbb{I}[Z\leq \theta]]}{{\rm E}[|Z-\theta|]}
,
$$
respectively. For the absolute loss function,~$\theta^\rho_{\alpha}(P)$ in~(\ref{Mquantdef}) therefore coincides with the usual order-$\alpha$ quantile, whereas for the quadratic loss function, it similarly provides a uniquely defined order-$\alpha$ expectile. For our later purposes, it is important to note that the larger~$G^\rho(\theta)(\leq 1/2)$ (resp.,~$1-G^\rho(\theta-0)(\leq 1/2)$, where~$H(\theta-0)$ denotes the limit of~$H(t)$ as~$t\nearrow \theta$), the less~$\theta$ is outlying below (resp., above) the central location~$\theta_{1/2}^\rho(P)$. 
\vspace{-.6mm}
 Therefore, $M\!D^\rho(\theta,P):=\min(G^\rho(\theta),1-G^\rho(\theta-0))$ measures the centrality---as opposed to outlyingness---of the location~$\theta$ with respect to~$P$. In other words, $M\!D^\rho(\theta,P)$ defines a measure of \emph{statistical depth} over the real line; see \cite{ZuoSer2000A}. In the sequel, we will extend this ``\mbox{M-depth}" to the multivariate case. Note that, for~$d=1$ and~$\rho(t)=|t|$, the depth $M\!D^\rho(\theta,P)$ reduces to the \cite{Tuk1975} halfspace depth.

%%%%%%%%%%%%%%%%%%%%%%%%%%%%%%

\section{Our multivariate M-quantiles} 
\label{secproposedmultiquantile}

Since the original multivariate \mbox{M-quantiles} from \cite{BreCha1988}, that actually include the celebrated geometric quantiles from \cite{Cha1996} and the recent geometric expectiles from \cite{Heretal2018}, several concepts of multivariate M-quantiles have been proposed. For the sake of completeness, we will describe these multivariate M-quantiles, as well as those from \cite{Breetal2001} and \cite{Koketal2002}, in Appendix~\ref{secAppA}. For now, it is only important to mention that, possibly after an unimportant reparametrization, all aforementioned multivariate M-quantiles can be written as functionals~$P\mapsto \theta^\rho_{\alpha,u}(P)$ that take values in~$\R^d$ and are indexed by a scalar order~$\alpha\in(0,1)$ and a direction~$u\in\mathcal{S}^{d-1}:=\{\zb\in\R^d:\|\zb\|^2:=\zb'\zb=1\}$; here, $P$ is a probability measure over~$\R^d$. Typically, $\theta^\rho_{\alpha,u}(P)$ does not depend on~$u$ for~$\alpha=1/2$, and the resulting common location is seen as the center (the ``median") of~$P$. Our multivariate M-quantiles will also be of a directional nature but they will be \emph{hyperplane-valued} rather than \emph{point-valued}. For~$d=1$, it is often important to know whether some test statistic takes a value below or above a given quantile, that is used as a critical value; for~$d>1$, hyperplane-valued quantiles, unlike point-valued ones, could similarly be used as critical values with vector-valued test statistics.

Before describing our M-quantiles, we introduce the class of probability measures for which they will be well-defined. To this end, for~$\ub\in\mathcal{S}^{d-1}$ and a probability measure~$P$ over~$\R^d$,  denote as~$P_\ub$ the probability measure over~$\R$ that is defined through~$P_\ub[A]=P[\{\zb\in\R^d:\ub'\zb\in A\}]$; in other words, if~$\Zb$ is a random $d$-vector with distribution~$P$, then~$P_\ub$ is the distribution of~$\ub'\Zb$. Consider then the collection~$\mathcal{P}_d^{\rho}$ of probability measures~$P$ over~$\R^d$ such that (i) no hyperplane of~$\R^d$ has~$P$-probability mass one and such that (ii)~$\int_{-\infty}^\infty |\psi_-(z-\theta)| \,dP_\ub(z)<\infty$ for any~$\theta\in\R$ and~$\ub\in\mathcal{S}^{d-1}$. Note that~$\mathcal{P}_1^{\rho}$ coincides with the collection of probability measures~$\mathcal{P}^{\rho}$ introduced in Section~\ref{secMquantile}. Our concept of multivariate M-quantile is then the following.   

\begin{Defin}
\label{directMquant}
Fix~$\rho\in\mathcal{C}$ and~$P\in\mathcal{P}^\rho_d$. Let~$Z$ be a random $d$-vector with distribution~$P$. Then, for any~$\alpha \in(0,1)$ and $\ub\in\mathcal{S}^{d-1}$, the \emph{\mbox{order-$\alpha$} M-quantile of~$P$ in direction~$\ub$} is the hyperplane
$$
%\pi^\rho_{\alpha,\ub}
%= 
\pi^\rho_{\alpha,\ub}(P)
= 
\Big\{ \zb\in\R^d : \ub'\zb = \theta^\rho_{\alpha}(P_\ub)\Big\},
$$
where  
$
\theta^\rho_{\alpha}(P_\ub)
$
is the \mbox{order-$\alpha$} M-quantile of~$P_\ub$; see~(\ref{Mquantdef}). The corresponding upper-halfspace~$%H^{\rho}_{\alpha,\ub}=
H^{\rho}_{\alpha,\ub}(P)=\big\{ \zb\in\R^d : \ub'\zb \geq \theta^\rho_{\alpha}(P_\ub)\big\}$ will be called \emph{\mbox{order-$\alpha$} M-quantile halfspace of~$P$ in direction~$\ub$}.
\end{Defin}

For~$\rho(t)=|t|$, these quantile hyperplanes
%~$\pi^\rho_{\alpha,\ub}(P)$ 
%and halfspaces
%~$H^{\rho}_{\alpha,\ub}(P)$ 
reduce to those from~\cite{PaiSim2011} (see also
% \citealp{Haletal2010} and 
\citealp{KonMiz2012}), whereas~$\rho(t)=t^2$ provides the proposed multivariate expectiles. For any loss function~$\rho$, the hyperplanes~$\pi^\rho_{\alpha,\ub}$ are linked in a straightforward way to the direction~$\ub$: they are simply orthogonal to~$\ub$. In contrast, the point-valued competitors typically depend on~$u$ in an intricate way, and there is in particular no guarantee that~$\theta^\rho_{\alpha,u}(P)$ belongs to the halfline with direction~$u$ originating from the corresponding median (see above). Note that the ``intercepts" of our M-quantile hyperplanes are the univariate \mbox{M-quantiles} of the projection~$\ub'\Zb$ of~$\Zb$ onto~$\ub$ (where~$\Zb$ has distribution~$P$), hence also allow for a direct interpretation.  

Irrespective of the loss function~$\rho$, competing multivariate M-quantiles fail to be equivariant under affine transformations. As announced in the introduction, our \mbox{M-quantiles} improve on this. We have the following result. 

\begin{Theor}
\label{Theoraffequiv}
Fix~$\rho\in\mathcal{C}_{\rm aff}$ and~$P\in\mathcal{P}^\rho_d$, where~$\mathcal{C}_{\rm aff}(\subset \mathcal{C})$ is the collection of power loss functions~$\rho(t)=|t|^r$ with~$r\geq 1$. Let~$\Ab$ be an invertible~$d\times d$ matrix and $\bb$ be a $d$-vector. Then, for any~$\alpha \in(0,1)$ and~$\ub\in\mathcal{S}^{d-1}$,
$$
\pi^\rho_{\alpha,\ub_\Ab}(P_{\Ab,\bb}) = \Ab \pi^\rho_{\alpha,\ub}(P)+\bb
	\quad 
	\textrm{and}
	\quad
	H^{\rho}_{\alpha,\ub_\Ab}(P_{\Ab,\bb}) = \Ab H^{\rho}_{\alpha,\ub}(P)+\bb
,
$$
where~$\ub_\Ab:=(\Ab^{-1})'\ub/\|(\Ab^{-1})'\ub\|$ and where~$P_{\Ab,\bb}$ is the distribution of~$\Ab\Zb+\bb$ when~$\Zb$ is a random $d$-vector with distribution~$P$. 
\end{Theor}

In the univariate case, the M-quantiles associated with~$\rho\in\mathcal{C}_{\rm aff}$ are known as \emph{$L_r$-quantiles} and were used for testing symmetry in nonparametric regression (\citealp{Chen1996}); the estimation of extreme $L_r$-quantiles was also recently investigated in \cite{Daoetal2019}. While Theorem~\ref{Theoraffequiv} above shows in particular that quantile and expectile hyperplanes are affine-equivariant, the restriction to~$\mathcal{C}_{\rm aff}$ cannot be dropped. For instance, for fixed~$c>0$, the M-quantile hyperplanes~$\pi^{\rho_c}_{\alpha,\ub}(P)$, associated with the Huber loss function in~(\ref{rhoc}), fail to be affine-equivariant. Our multivariate extension is not to be blamed for this, however, since the corresponding univariate M-quantiles~$\theta^{\rho_c}_{\alpha}(P)$ themselves fail to be scale-equivariant. Actually, it can be checked that if~$\rho\in\mathcal{C}$ makes the univariate M-quantiles~$\theta^{\rho}_{\alpha}(P)$ scale-equivariant, then our multivariate M-quantiles associated with~$\rho$ are affine-equivariant in the sense described in Theorem~\ref{Theoraffequiv}. 
 
At first sight, a possible advantage of any point-valued \mbox{M-quantiles}~$\thetab^\rho_{\alpha,\ub}(P)$ is that they naturally generate contours and regions. More precisely, they allow  considering, for any~$\alpha\in(0,\frac{1}{2}]$, the \mbox{order-$\alpha$} M-quantile contour~$\{\thetab^\rho_{\alpha,\ub}(P): \ub\in\mathcal{S}^{d-1}\}$, the interior part of which is then the corresponding order-$\alpha$ \mbox{M-quantile} region. Our hyperplane-valued \mbox{M-quantiles}, however, also provide centrality regions, hence corresponding contours.

\begin{Defin}
\label{definMregions}
Fix~$\rho\in\mathcal{C}$ and~$P\in\mathcal{P}^\rho_d$. For any~$\alpha\in(0,1)$, the \mbox{order-$\alpha$} M-quantile region of~$P$ is
$
R^{\rho}_{\alpha}(P)
=
\bigcap_{\ub\in\mathcal{S}^{d-1}}  H^{\rho}_{\alpha,\ub}(P)
$
and the corresponding \mbox{order-$\alpha$} contour is the boundary~$\partial R^{\rho}_{\alpha}(P)$ of~$R^{\rho}_{\alpha}(P)$. 
\end{Defin}

Theorem~\ref{Theorlemjones} entails that the univariate M-quantiles~$\theta_\alpha^\rho(P)$ in~(\ref{Mquantdef}) are monotone non-decreasing functions of~$\alpha$. 
%(throughout, we write ``increasing" and ``decreasing" for ``non-decreasing" and ``non-increasing", respectively). 
A direct corollary is that the regions~$R^{\rho}_{\alpha}(P)$ are non-increasing with respect to inclusion.
%: if~$\alpha_1>\alpha_2$, then $R^{\rho}_{\alpha_1}(P)\subset R^{\rho}_{\alpha_2}$. 
The proposed regions enjoy many nice properties compared to their competitors resulting from point-valued M-quantiles, as we show on the basis of Theorem~\ref{TheorSupportAndcompactness} below. To state the result, we need to define the following concept: for a probability measure~$P$ over~$\R^d$, the \emph{$C$-support of~$P$} is
$
C_P 
:=
%\Big
\{ \zb\in\R^d : 
%\min( 
P[\ub'\Zb\leq \ub'\zb]  
%, P[\ub'\Zb\geq \ub'\zb] ) 
> 0 \textrm{ for any } \ub\in\mathcal{S}^{d-1} 
%\Big
\} 
,
$
where the random $d$-vector~$\Zb$ has distribution~$P$. Clearly, $C_P$ can be thought of as the convex hull of $P$'s support. We then have the following result. 

\begin{Theor}
\label{TheorSupportAndcompactness}
Fix~$\rho\in\mathcal{C}$ and~$P\in\mathcal{P}^\rho_d$. Then, for any~$\alpha\in(0,1)$, the region~$R_{\alpha}^{\rho}(P)$ is 
a convex and compact subset of~$C_P$. Moreover, if~$\rho\in\mathcal{C}_{\rm aff}$, then  
$
R^{\rho}_{\alpha}(P_{\Ab,\bb})
=
\Ab R^{\rho}_{\alpha}(P)+\bb
$
for any invertible~$d\times d$ matrix~$\Ab$ and $d$-vector~$\bb$.
\end{Theor}

No competing M-quantile regions combine these properties. For instance, the original M-quantile regions from \cite{BreCha1988}, hence also the geometric quantile regions from \cite{Cha1996} and their expectile counterparts from \cite{Heretal2018}, may extend far beyond the convex hull of the support; for geometric quantile regions, this was formally proved in \cite{GirStu2017}. This was actually the motivation for the alternative proposals in \cite{Breetal2001} and \cite{Koketal2002}. The regions introduced in these two papers, however, may fail to be convex, which is unnatural. More generally, none of the competing M-quantile or expectile regions are affine-equivariant. This may result in quite pathological behaviors: for instance, Theorem~2.2 from \cite{GirStu2017} implies that, if~$P$ is an elliptically symmetric probability measure admitting a density~$f$, then, for small values of~$\alpha$, the geometric quantile contours from \cite{Cha1996} are ``orthogonal" to the principal component structure of~$P$, in the sense that these contours are furthest (resp., closest) to the symmetry center of~$P$ in the last (resp., first) principal direction. In contrast, the affine-equivariance result in Theorem~\ref{TheorSupportAndcompactness} ensures that, in such a distributional setup, the shape of our M-quantile contours will match the principal component structure of~$P$.

\begin{figure}[htbp]
  \vspace*{-0mm} 
    \begin{center} 
    \includegraphics[width=\textwidth]{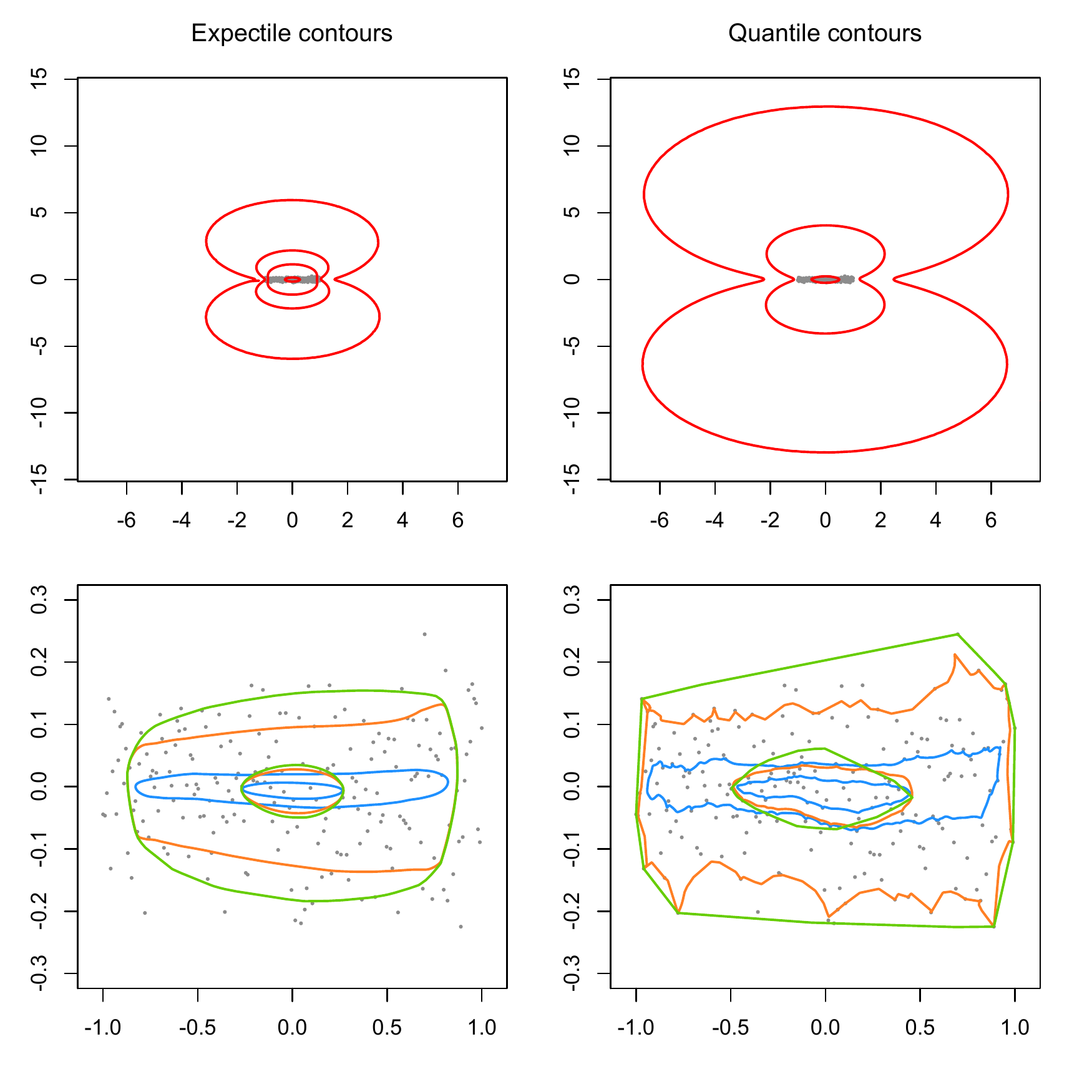}
     \end{center}
     \vspace*{-5mm}

%\caption{\it (Top left:) several types of quantile contours of orders~$\alpha=.005$ and~$.25$ associated with the cigar shaped dataset described in  Section~\ref{secproposedmultiquantile}; the quantile contours are those from  \cite{Breetal2001} (blue), the $(\delta=10)$-version of the \cite{Koketal2002} contours (orange), and the proposed contours (green; in this quantile case, these are the Tukey depth contours; see Section~\ref{secMdepth}).
%%
%(Top right:) the corresponding expectile contours. 
%%
%(Bottom:) the geometric quantile contours from \cite{Cha1996} (left) and the geometric expectile contours from \cite{Heretal2018} (right), for the same dataset and for orders~$\alpha=.00001$, $.0005$, $.005$, and $.25$ (the quantile contour with the smallest order is outside the plot).}
%
\caption{\it (Top:) the geometric expectile contours from \cite{Heretal2018} (left) and  geometric quantile contours from \cite{Cha1996} (right), for the cigar-shaped data described in Section~\ref{secproposedmultiquantile} and for orders~$\alpha=.00001$, $.0005$, $.005$, and $.25$ (for the smallest order~$\alpha$, the quantile contour is outside the plot).
(Bottom left:) the expectile contours from \cite{Breetal2001} (blue), the $(\delta=10)$-version of the \cite{Koketal2002} expectile contours (orange), and the proposed expectile contours (green),
%; these coincide with the Tukey depth contours; see Section~\ref{secMdepth}), 
for the same data and for orders~$\alpha=1/n=.005$ and~$.25$; we use the same value of~$\delta$ as in \cite{Koketal2002}.
(Bottom right:) the quantile versions of the contours in the bottom left panel. In each panel, the $n=200$ data points are plotted in grey.}
\label{Fig1}
\end{figure}

As an illustration, we consider the ``cigar-shaped" data example from \cite{Breetal2001} and \cite{Koketal2002}, for which~$P=P_n$ is the empirical probability measure associated with $n=200$ bivariate observations whose $x$-values form a uniform grid in~$[-1,1]$ and whose $y$-values are independently drawn from the normal distribution with mean 0 and variance~$.01$. Figure~\ref{Fig1} draws, for several orders~$\alpha$, the various quantile and expectile contours mentioned in the previous paragraph. Our contours~$\partial R^{\rho}_{\alpha}(P_n)$ were computed by replacing the intersection in Definition~\ref{definMregions} by an intersection over~$L=500$ equispaced directions~$\ub$ in~$\mathcal{S}^1$ (all competing contours require a similar discretization). The results show that, like the geometric quantiles from \cite{Cha1996}, their expectile counterparts from \cite{Heretal2018} may extend beyond the convex hull of the data points. The aforementioned pathological behavior of the extreme geometric quantiles \cite{Cha1996} relative to the principal component structure of~$P$ not only shows for these quantiles but also for the corresponding expectiles. Finally, the outer quantile/expectile regions from \cite{Breetal2001} and \cite{Koketal2002} are non-convex in most cases. In line with Theorem~\ref{TheorSupportAndcompactness}, our M-quantile regions and contours do not suffer from these deficiencies.

%%%%%%%%%%%%%%%%%%%%%%%%%%%%%%%%%%%%%%%%%%%%%%%%%%

%%%%%%%%%%%%%%%%%%%%%%%%%%%%%%%%%%%%%%%%%%%%%%%%%% 

\section{Halfspace \mbox{M-depth}} 
\label{secMdepth}
 
Our M-quantile regions~$R^{\rho}_{\alpha}(P)$ are \emph{centrality regions}, in the sense that they group locations~$\zb$ in the sample space~$\R^d$ according to their centrality with respect to the underlying distribution~$P$. This defines the following concept of \emph{depth} (throughout, we let~$\sup\emptyset:=0$). 

\begin{Defin}
\label{definMdepth}	
\ \hspace{-2.5mm}  
Fix~$\rho\in\mathcal{C}$ and~$P\in\mathcal{P}_d^{\rho}$. Then, the corresponding \emph{halfsp\-ace \mbox{M-depth}} of~$\zb$ with respect to~$P$ is 
$
%\label{deviceregdepth}	
M\!D^\rho(\zb,P)
=
\sup
\{
\alpha> 0 : \zb\in R^{\rho}_{\alpha}(P)
\}
$.  
\end{Defin}

While this was not done in the literature (neither for general \mbox{M-quantiles} nor for expectiles), other \mbox{M-depth} concepts could similarly be defined  from any collection of competing M-quantile regions.
% (it has been done only for geometric quantiles, where the corresponding depth is the \emph{spatial} depth; see, e.g., \citealp{VarZha2000} or \citealp{Ser2010}). 
However, these \mbox{M-depths} would, irrespective of~$\rho$, fail to meet one of the most classical requirements for depth, namely affine invariance; see \cite{ZuoSer2000A}. Our \mbox{M-depth} is better in this respect; see Theorem~\ref{Theordepthaxioms}(i) below. 

For any depth, the corresponding depth regions, that collect locations with depth larger than or equal to a given level~$\alpha$, are of particular interest. The following result shows that halfspace \mbox{M-depth} regions strictly coincide with the centrality regions introduced in the previous section.

\begin{Theor}
\label{TheorRegionsOrig}
Fix~$\rho\in\mathcal{C}$ and~$P\in\mathcal{P}_d^{\rho}$. Then, for any~$\alpha\in(0,1)$, the level-$\alpha$ depth region~$\{\zb\in\R^d:M\!D^\rho(\zb,P)\geq \alpha\}$ coincides with~$R^{\rho}_{\alpha}(P)$.  
%the centrality region~$R^{\rho}_{\alpha}(P)$.
\end{Theor}

This result has several interesting consequences. First, it implies that the depth~$M\!D^\rho$ reduces to the \cite{Tuk1975} halfspace depth for~$\rho(t)=|t|$ (since the corresponding centrality regions~$R^{\rho}_{\alpha}(P)$ are known to be the Tukey depth regions; see, e.g., Theorem~2 in \citealp{KonMiz2012}). 
Second, Theorems~\ref{TheorSupportAndcompactness}--\ref{TheorRegionsOrig}  show that halfspace \mbox{M-depth} regions are convex, so that our \mbox{M-depth} is quasi-concave: for any~$\zb_0,\zb_1\in\R^d$ and~$\lambda\in(0,1)$, one has 
$
M\!D^{\rho}((1-\lambda)\zb_0+\lambda \zb_1,P)
\geq 
\min(M\!D^{\rho}(\zb_0,P),M\!D^{\rho}(\zb_1,P))
$.
Third, since Theorems~\ref{TheorSupportAndcompactness}--\ref{TheorRegionsOrig} imply that the mapping~$\zb\mapsto M\!D^\rho(\zb,P)$ has closed upper level sets, this mapping is upper semicontinuous over~$\R^d$ (it is actually continuous  over~$\R^d$ if~$P$ assigns probability zero to all hyperplanes of~$\R^d$; see Lemma~\ref{LemGandMDepthContinuous}). 
Fourth, the compactness of halfspace \mbox{M-depth} regions with level~$\alpha>0$, which results again from Theorems~\ref{TheorSupportAndcompactness}--\ref{TheorRegionsOrig}, allows us to establish the existence of an M-deepest location.

\begin{Theor}
\label{Theordeepestpoint}
Fix~$\rho\in\mathcal{C}$ and~$P\in\mathcal{P}_d^{\rho}$. Then, $\sup_{\zb\in\R^d} M\!D^{\rho}(\zb,P)=M\!D^{\rho}(\zb_*,P)$ for some~$\zb_*\in\R^d$. 
\end{Theor}

The M-deepest location~$\zb_*$ may fail to be unique.
%, that is, the M-deepest region may fail to be a singleton. 
For the halfspace Tukey depth, whenever a unique representative of the deepest locations is needed, a classical solution consists in considering the \emph{Tukey median}, that is defined as the barycenter of the deepest region. The same solution can be adopted for our \mbox{M-depth} and the convexity of the M-deepest region will still ensure that this uniquely defined M-median has indeed maximal \mbox{M-depth}.

The following result shows that, for~$\rho\in\mathcal{C}_{\rm aff}$ (a restriction that is required only for Part~(i) of the result), the halfspace \mbox{M-depth}~$M\!D^{\rho}$ is a \emph{statistical depth function}, in the axiomatic sense of \cite{ZuoSer2000A}. 
%We have the following result. 

\begin{Theor}
\label{Theordepthaxioms}
Fix~$\rho\in\mathcal{C}_{\rm aff}$ and~$P\in\mathcal{P}_d^{\rho}$. Then, $M\!D^{\rho}(\zb,P)$ satisfies the following properties:
(i) (affine invariance:) for any invertible $d\times d$ matrix~$\Ab$ and $d$-vector~$\bb$, $M\!D^{\rho}(\Ab\zb+\bb,P_{\Ab,\bb})=M\!D^{\rho}(\zb,P)$, where~$P_{\Ab,\bb}$ was defined in Theorem~\ref{Theoraffequiv};
 %where~$P_{\Ab,\bb}[B]:=P[\Ab^{-1}(B-b)]$ for any $d$-Borel set~$B$;
%
(ii) (maximality at the center:) if~$P$ is centrally symmetric about~$\thetab_*$ (i.e., $P[\thetab_*+B]=P[\thetab_*-B]$ for any $d$-Borel set~$B$), then $M\!D^{\rho}(\thetab_*,P)\geq M\!D^{\rho}(\zb,P)$ for any~$d$-vector~$\zb$;
(iii) (monotonicity along rays:)  if~$\thetab$ has maximum \mbox{M-depth} with respect to~$P$, then, for any~$\ub\in\mathcal{S}^{d-1}$, $r \mapsto M\!D^{\rho}(\thetab+r \ub,P)$ is monotone non-increasing in~$r(\geq 0)$ and~$M\!D^{\rho}(\thetab+r \ub,P)=0$ for any~$r> r_\ub(P):=\sup \{r > 0 : \thetab + r \ub\in C_P \}
\linebreak
(\in (0,+\infty])$; 
(iv) (vanishing at infinity:) as~$\|\zb\|\to \infty$, $M\!D^{\rho}(\zb,P)\to 0$.
\end{Theor}

As mentioned above, $M\!D^\rho$ reduces to the Tukey depth for~$\rho(t)=|t|$. For any other~$\rho$ function, the depth~$M\!D^\rho$ is, to the authors' best knowledge, original. In particular, the (\emph{halfspace}) \emph{expectile depth} obtained for~$\rho(t)=t^2$ has not been considered so far. While, as already mentioned, competing concepts of multivariate expectiles would provide alternative concepts of expectile depth (through the corresponding expectile regions as in Definition~\ref{definMdepth}), the following result hints that our  
%expectile depth---or, more generally, our \mbox{M-depth}---
construction is the natural one.

\begin{Theor}
\label{reginter}
Fix~$\rho\in\mathcal{C}$ and~$P\in\mathcal{P}_d^{\rho}$. Then, for any~$z\in\R^d$,
\begin{equation}
\label{mequ}
M\!D^\rho(\zb,P) 
=
\inf_{\ub\in\mathcal{S}^{d-1}} 
%G^\rho_\ub(\ub'z),
%\textrm{ with }
%G^\rho_\ub(\ub'z)
%:=
\frac{{\rm E}[|\psi_-(\ub'(\Zb-\zb))| \mathbb{I}[\ub'(\Zb-\zb)\leq 0]]}{{\rm E}[|\psi_-(\ub'(\Zb-\zb))|]}
,\end{equation}
where~$\psi_-$ is the left-derivative of~$\rho$ and where~$\Zb$ has distribution~$P$. 
\end{Theor}

For~$\rho(t)=|t|$, we have~$\psi_-(t)=\mathbb{I}[t>0]-\mathbb{I}[t\leq 0]$, so that Theorem~\ref{reginter} confirms that our \mbox{M-depth} then coincides with the halfspace Tukey depth
%\begin{equation}
%\label{DefinHalfspacedepth}
$$
H\!D(\zb,P) 
=
\inf_{\ub\in\mathcal{S}^{d-1}}
P[\ub'\Zb\leq \ub'\zb] 
,
$$
%\end{equation}
that records the most extreme (lower-)outlyingness of~$\ub'\zb$ with respect to the distribution of~$\ub'\Zb$. The \mbox{M-depth} in~(\ref{mequ}) can be interpreted in the exact same way but replaces standard quantile outlyingness with M-quantile outlyingness; see the last paragraph of Section~\ref{secMquantile}. For~$\rho(t)=t^2$, Theorem~\ref{reginter} states that our expectile depth can be equivalently defined as 
\begin{equation}
\label{DefinExpectileDepth}
E\!D(\zb,P) 
=
\inf_{\ub\in\mathcal{S}^{d-1}}
\frac{{\rm E}[|\ub'(\Zb-\zb)| \mathbb{I}[\ub'(\Zb-\zb)\leq 0]]}{{\rm E}[|\ub'(\Zb-\zb)|]}
\cdot 
\end{equation}
Of course, we could similarly consider the continuum of halfspace \mbox{M-depth}s associated with the Huber loss functions~$\rho_c$ in~(\ref{rhoc}). However, since our work was mainly motivated by expectiles and multiple-output expectile regression, we will mainly focus on expectile depth in the next sections. 

Before doing so, we state three consistency results that can be proved on the basis of Theorem~\ref{reginter}. We start with the following uniform consistency result, that extends to an arbitrary halfspace M-depth the Tukey depth result from \cite{DonGas1992},  Section~6.

\begin{Theor}
\label{TheorMDepthConsistency}
Fix~$\rho\in\mathcal{C}$ and~$P\in\mathcal{P}_d^{\rho}$. Let~$P_n$ be the empirical probability measure associated with a random sample of size~$n$ from~$P$. Then, 
$
\sup_{\zb\in\R^d} |M\!D^\rho(\zb,P_n)-M\!D^\rho(\zb,P) |\to 0
$
almost surely as~$n\to\infty$. 
\end{Theor}

Jointly with a general result on the consistency of M-estimators (such as Theorem~2.12 in \citealp{Kos2008}), this uniform consistency property allows us to establish consistency of the sample M-deepest point. 

\begin{Theor}
\label{TheorMmedianConsistency}
Fix~$\rho\in\mathcal{C}$ and~$P\in\mathcal{P}_d^{\rho}$. Let~$P_n$ be the empirical probability measure associated with a random sample of size~$n$ from~$P$. Let~$\zb_*(P)$ be the halfspace M-median of~$P$, that is, the barycenter of~$\{\zb\in\R^d:M\!D^{\rho}(\zb,P)=\max_{\yb\in\R^d} M\!D^{\rho}(\yb,P)\}$ and let~$\zb_*(P_n)$ be the halfspace M-median of~$P_n$. Then,  $\zb_*(P_n)\to \zb_*(P)$ almost surely as~$n\to\infty$. 
\end{Theor}

Finally, we consider consistency of depth regions. This was first discussed in \cite{HeWan1997}, where the focus was mainly on elliptical distributions. While results under milder conditions were obtained in \cite{Kim2000} and \cite{ZuoSer2000B}, we will here exploit the general results from \cite{Dyc2016}. We need to introduce the following concept:~$P(\in\mathcal{P}_d^{\rho})$ is said to have a \emph{connected support} if and only if we have~$P[a<\ub'\Zb<b]>0$ whenever~$\min(P[\ub'\Zb\leq a],P[\ub'\Zb\geq b])>0$ for some~$\ub\in\mathcal{S}^{d-1}$ and~$a,b\in\R$ with~$a<b$  (as usual, $\Zb$ denotes a random $d$-vector with distribution~$P$). We then have the following result. 

\begin{Theor}
\label{TheorMRegionsConsistency}
Fix~$\rho\in\mathcal{C}$ and assume that~$P\in\mathcal{P}_d^{\rho}$ has a connected support and assigns probability zero to all hyperplanes in~$\R^d$. Let~$P_n$ be the empirical probability measure associated with a random sample of size~$n$ from~$P$. 
Then,  for any compact interval~$\mathcal{I}$ in~$(0,\alpha_*)$, with~$\alpha_*:=\max_{\zb\in\R^d} M\!D^{\rho}(\zb,P)$,
$$
\sup_{\alpha\in \mathcal{I}}
d_H(R^\rho_\alpha(P_n),R^\rho_\alpha(P))
\to
0
$$  
almost surely as~$n\to\infty$, where~$d_H$ denotes the Hausdorff distance. 
\end{Theor}

As announced, we now focus on the particular case of expectile depth.

%%%%%%%%%%%%%%%

\section{Halfspace expectile depth}  
\label{secExpectiledepth}

Below, we derive further properties of (halfspace) expectile depth that show why this particular \mbox{M-depth} should be appealing for practitioners. Throughout, we write~$\mathcal{P}_d$ for the collection~$\mathcal{P}_d^{\rho}$ of probability measures over~$\R^d$ for which expectile depth is well-defined, that is, the one associated with~$\rho(t)=t^2$. Clearly,~$\mathcal{P}_d$ collects the probability measures that (i) do not give $P$-probability one to any hyperplane of~$\R^d$ and that (ii) have finite first-order moments. Note that this moment assumption is required even for~$d=1$; as for~(i), it only rules out distributions that are actually over a lower-dimensional Euclidean space.

%%%%%%%%%%%%%%%%%%%%

\subsection{Further properties of expectile depth}
 \label{secexpectfurtherprop}

For the \mbox{depth}~$M\!D^\rho(\zb,P)$, a direction~$\ub_0$ is said to be \emph{minimal} if it achieves the infimum in~(\ref{mequ}). For the halfspace Tukey depth, 
%it is well-known that 
such a minimal direction does not always exist. A bivariate example is obtained for~$\zb=(1,0)'\in\R^2$ and~$P=\frac{1}{2}P_1+\frac{1}{2}P_2$, where~$P_1$ is the bivariate standard normal distribution and~$P_2$ is the Dirac distribution at~$(1,1)'$. In contrast, the continuity, for any~\mbox{$P\in\mathcal{P}_d$}, of the function whose infimum is considered over~$\mathcal{S}^{d-1}$ in~(\ref{DefinExpectileDepth}) (see Lemma~\ref{LemExpectileDepthSmoothness}(i)) and the compactness of~$\mathcal{S}^{d-1}$ imply that
\begin{equation}
\label{DefinExpectileDepth2}
E\!D(\zb,P) 
=
\min_{\ub\in\mathcal{S}^{d-1}} 
\frac{{\rm E}[ |\ub'(\Zb-\zb)| \mathbb{I}[\ub'(\Zb-\zb)\leq 0]]}{{\rm E}[|\ub'(\Zb-\zb)|]}
,
\end{equation}
so that a minimal direction always exists for expectile depth (in the mixture example above, $\ub_0=(-1,0)'$ is a minimal direction for~$E\!D(\zb,P)$).  

%For any halfspace \mbox{M-depth}, the deepest point is a location functional. F
Now, for halfspace Tukey depth, which is an $L_1$-concept, the deepest point is not always unique and its unique representative, namely the Tukey median, is a multivariate extension of the univariate median. Also, the depth of the Tukey median may depend on the underlying probability measure~$P$. Our expectile depth, that is rather of an $L_2$-nature, is much different.

\begin{Theor}
\label{TheorExpectDeepestPoint}
For any~$P\in\mathcal{P}_d$, the expectile depth $E\!D(\zb,P)$ is uniquely maximised at $\zb=\mub_P:={\rm E}[\Zb]$ $($where $\Zb$ is a random $d$-vector with distribution~$P)$ and the corresponding maximum depth is~$E\!D(\mub_P,P)=1/2$. 
\end{Theor}

Expectile depth regions therefore always provide nested regions around the mean vector~$\mub_P$, which should be appealing to practitioners. Since the maximal expectile depth is~$1/2$ for any~$P$, a natural affine-invariant test for~$\mathcal{H}_0:\mub_P=\mub_0$, where~$\mu_0\in\R^d$ is fixed, is the one rejecting~$\mathcal{H}_0$ for large values of
$
%T_n:=\frac{1}{2}-E\!D(\mub_0,P_n)
T_n:=(1/2)-E\!D(\mub_0,P_n)
$,
where~$P_n$ is the empirical probability measure associated with the sample~$\Zb_1,\ldots,\Zb_n$ at hand. Investigating the properties of this test is beyond the scope of the present work.

%%%%%%%%%%

We turn to another distinctive aspect of expectile depth. Theorem~\ref{Theordepthaxioms}(iii) shows that halfspace \mbox{M-depth} decreases monotonically when one moves away from a deepest point along any ray. This decrease, however, may fail to be strict (in the sample case, for instance, the halfspace Tukey depth is piecewise constant, hence will fail to be strictly decreasing). In contrast, expectile depth always offers a strict decrease (until, of course, the minimal depth value zero is reached, if it is). We have the following result.

\begin{Theor}
\label{TheorMonotExpectDepth}
Fix~$P\in\mathcal{P}_d$ and~$\ub\in\mathcal{S}^{d-1}$. Let~$r_\ub(P)=\sup \{r > 0 : \mub_P + r \ub\in C_P \}(\in (0,+\infty])$. Then, $r \mapsto E\!D(\mub_P+r \ub,P)$ is monotone strictly decreasing in~$[0,r_\ub(P)]$ and~$E\!D(\mub_P+r \ub,P)=0$ for~$r\geq r_\ub(P)$.  
\end{Theor}

Our \mbox{M-depth}s are upper semicontinuous functions of~$\zb$; see Section~\ref{secMdepth}. However, continuity does not hold in general (in particular, the piecewise constant nature of the halfspace Tukey depth  for empirical distributions rules out continuity in the sample case). Expectile depth is smoother.

\begin{Theor}
\label{TheorexpectileDepthContinuous}
Fix~$P\in\mathcal{P}_d$. Then, (i) $\zb\mapsto E\!D(\zb,P)$ is uniformly continuous over~$\R^d$; 
(ii) for~$d=1$,~$z\mapsto E\!D(z,P)$ is left- and right-differentiable over~$\R$;  
(iii) for~$d\geq 2$, if~$P$ is smooth in a neighbourhood~$\mathcal{N}\!$\!\, of~$\zb_0$ (meaning that for any~$\zb\in\mathcal{N}$, any hyperplane containing~$\zb$ has $P$-probability zero), then $\zb\mapsto E\!D(\zb,P)$ admits directional derivatives at~$\zb_0$ in all directions. 
%; \textcolor{red}{(iii) if, moreover,~$0<E\!D(z_0,P)<1/2$, then $z\mapsto E\!D(\zb,P)$ is continuously differentiable at~$z_0$}. 
\end{Theor}

%To further comment on the smoothness of expectile depth, we will need the following definitions. The probability measure~$P$ over~$\R^d$ is \emph{smooth at~$z_0(\in\R^d)$} if and only if any hyperplane containing~$z_0$ has~$P$-probability zero. It is said to be \emph{smooth in a neighbourhood of~$z_0(\in\R^d)$} if and only if there exists~$\delta>0$ such that~$P$ is smooth at any~$z$ such that~$\|z-z_0\|<\delta$. We then have the following result.  

\begin{figure}
\center
\includegraphics[width=\textwidth]{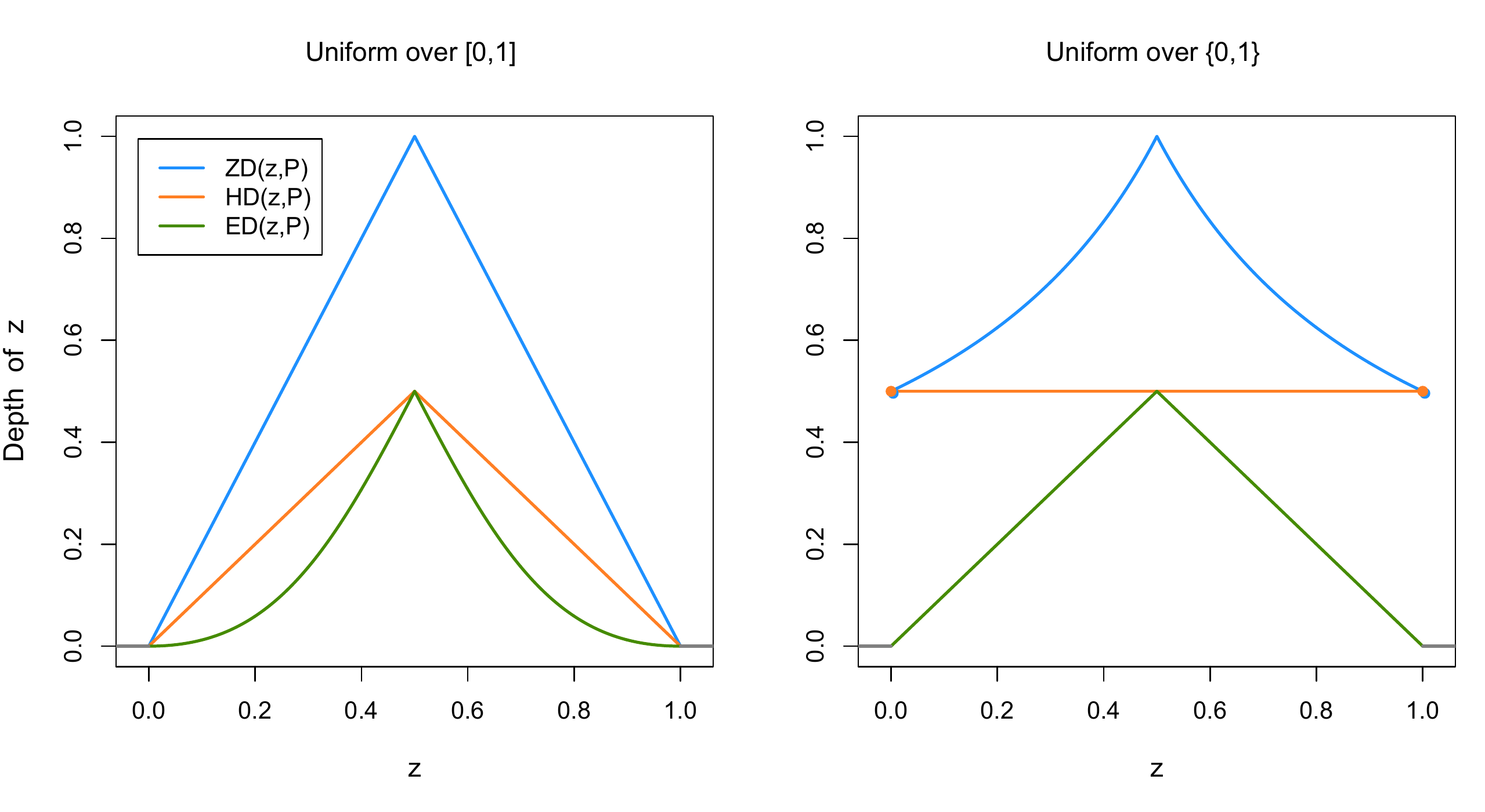}
\caption{Plots, as functions of~$z$, of the zonoid depth~$Z\!D(z,P)$ (blue), of the halfspace Tukey depth~$H\!D(z,P)$ (orange) and of the expectile depth~$E\!D(z,P)$ (green), when $P$ is the uniform over the interval~$[0,1]$ (left) and the uniform over the pair~$\{0,1\}$ (right). For both probability measures, all depth functions take value zero outside~$[0,1]$.}
\label{Fig2}
\end{figure}

We illustrate these results on the following univariate examples. It is easy to check that if~$P$ is the uniform measure over the interval~$\mathcal{I}=[0,1]$, then
\begin{equation}
\label{expecunivexam1}
E\!D(z,P) 
=
\frac{\min(z^2,(1-z)^2)}{z^2+(1-z)^2}\,\mathbb{I}[z \in\mathcal{I}]
%\quad
\ \textrm{ and }\
%\quad
H\!D(z,P)
=
\min(z,1-z)\mathbb{I}[z \in\mathcal{I}]
%=
%\frac{1}{2} Z\!D(z,P)  
,
\end{equation}
%$$
%Z\!D(z,P)
%=
%2H\!D(z,P)
%=
% \mathbb{I}[z \in\mathcal{I}]
%$$ 
whereas if~$P$ is the uniform over the pair~$\{0,1\}$, then 
%Let~$Z$ be uniform on~$\{0,1\}$. Then, for $0\leq z\leq 1$, we have 
%$$
%\frac{{\rm E}[ |Z-z| \mathbb{I}[Z\leq z]]}{{\rm E}[|Z-z|]}
%=
%\frac{\frac{1}{2}|z|}{\frac{1}{2}|z|+\frac{1}{2}|z-1|}
%=
%\frac{z}{z+1-z}
%=
%z
%$$
%and
%$$
%\frac{{\rm E}[ |-(Z-z)| \mathbb{I}[-Z\leq -z]]}{{\rm E}[|-(Z-z)|]}
%=
%\frac{\frac{1}{2}|z-1|}{\frac{1}{2}|z|+\frac{1}{2}|z-1|}
%=
%\frac{1-z}{z+1-z}
%=
%1-z,
%$$
%so that, for any~$z\in\R$,
\begin{equation}
\label{expecunivexam2}
E\!D(z,P)
%=
%\min_{\ub\in\mathcal{S}^{d-1}} 
%\frac{{\rm E}[ |\ub'(\Zb-\zb)| \mathbb{I}[\ub'\Zb\leq \ub'\zb]]}{{\rm E}[|\ub'(\Zb-\zb)|]}
=
\min(z,1-z) \mathbb{I}[z \in\mathcal{I}]
%\quad
\ \textrm{ and }\
%\quad
H\!D(z,P)
=
\frac{1}{2}\, \mathbb{I}[z \in\mathcal{I}]
;
\end{equation}
%$$
%Z\!D(z,P) 
%=
%\frac{1}{2}\min\Big(\frac{1}{1-z},\frac{1}{z}\Big)\,\mathbb{I}[z \in\mathcal{I}]
%$$
see Figure~\ref{Fig2}. This illustrates uniform continuity of expectile depth, as well as left- and right-differentiability. Although both distributions are smooth in a neighborhood of~$z_0=1/2$, plain differentiability does not hold at~$z_0$, which results from the non-uniqueness of the corresponding minimal direction; see \cite{Dem2004}. For the sake of comparison, the figure also plots the Tukey depth~$H\!D(z,P)$ and the zonoid depth~$Z\!D(z,P)$ from \cite{KosMos1997}. Comparison with the latter depth is natural as it is also maximized at~$\mub_P$, hence has an $L_2$-flavor. Both the Tukey and zonoid depths are less smooth than the expectile one, and in particular, the zonoid depth is not continuous for the discrete uniform. Also, the zonoid depth region~$\{z\in\R:Z\!D(z,P)\geq \alpha\}$ in the left panel is the ($L_1$) interquantile interval~$[\frac{\alpha}{2},1-\frac{\alpha}{2}]$, which is not so natural for a depth of an $L_2$-nature. In contrast, for any~$P$, the level-$\alpha$ expectile depth region is the interexpectile interval~$[\theta^\rho_{\alpha}(P),\theta^\rho_{1-\alpha}(P)]$, with~$\rho(t)=t^2$, which reflects, for any~$\alpha$ (rather than for the deepest level only), the $L_2$-nature of expectile depth.

%%%%%%%%%%%%%%%%%%%%%%%%%%%%%%%%%%%%%%%%

%%%%%%%%%%%%%%%%%%%%%%%%%%%%%%%%%%%%%%%%

\subsection{Some multivariate examples}

Consider the case where~$P(\in\mathcal{P}_d)$ is the distribution of~$\Zb=\Ab \Yb+\mub$, where~$\Ab$ is an invertible $d\times d$ matrix, $\mub$ is a $d$-vector and~$\Yb=(Y_1,\ldots,Y_d)'$ is a spherically symmetric random vector, meaning that the distribution of~$\Ob\Yb$ does not depend on the $d\times d$ orthogonal matrix~$\Ob$. In other words, $P$ is elliptical with mean vector~$\mub$ and scatter matrix~$\Sigmab=\Ab\Ab'$. In the standard case where~$\Ab=I_d$ (the \mbox{$d$-dimensional} identity matrix) and~$\mub=0$, Theorem~\ref{Theorlemjones}(iv) provides
\begin{eqnarray*}
%\lefteqn{
\hspace{-9mm} 
E\!D(\zb,P) 
&\!\!\!\!=\!\!\!\!&
\min_{\ub\in\mathcal{S}^{d-1}} 
\!
\frac{{\rm E}[ |Y_1-\ub'\zb| \mathbb{I}[Y_1\leq \ub'\zb]]}{{\rm E}[|Y_1-\ub'\zb|]}
\\[3mm]
&\!\!\!\!=\!\!\!\!&
-
\frac{{\rm E}[ (Y_1+\|\zb\|) \mathbb{I}[Y_1\leq -\|\zb\|]]}{{\rm E}[|Y_1+\|\zb\||]}
%=
%\frac{-\int_{-\infty}^{-\|\zb\|} (t+\|\zb\|) dP^{Y_1}(t)}{\int_{-
%\infty}^{\infty} |t+\|\zb\|| dP^{Y_1}(t)}
=:
g(\|\zb\|)
,
\end{eqnarray*}
so that, for arbitrary~$\mub$ and~$\Sigmab$, affine invariance entails that $E\!D(\zb,P) =g(\|\zb\|_{\mu,\Sigma})$, with~$\|\zb\|_{\mu,\Sigma}^2:=(\zb-\mub)'\Sigmab^{-1}(\zb-\mub)$. Expectile depth regions are thus concentric ellipsoids that, under absolute continuity of~$P$, coincide with equidensity contours. The function~$g$ depends on the distribution of~$\Yb$: if~$\Yb$ is $d$-variate standard normal, then it is easy to check that 
$$
g(r)
=
\frac{1}{2}
-
\frac{1}{2\{(2/r)\phi(r)+\, 2 \Phi(r)-1\}}
,
$$
where~$\phi$ and~$\Phi$ denote the probability density function and cumulative distribution function of the univariate standard normal distribution, respectively. If~$\Yb$ is uniform over the unit ball~$B^d:=\{\zb\in\R^d:\|\zb\|\leq 1\}$ or on the unit sphere~$\mathcal{S}^{d-1}$, then one can show that
%~$Y_1^2\sim {\rm Beta}(\frac{1}{2},\frac{d-1}{2})$ for~$d>1$ (for~$d=1$, $Y_1$ is uniform on~$\{-1,1\}$), which allows to show
%$$
%g(r)
%=
%\frac{1}{2}-
%\frac{\sqrt{\pi} r  (1-r^2)^{(1-d)/2} \Gamma(\frac{d+1}{2})}{2 \Gamma(\frac{d}{2}) (1+(d-1)r^2 {_2F_1}(1,\frac{d}{2};\frac{3}{2};r^2))}
%$$
$$
g(r)
=
\omega_{d}(r)
:=
\bigg(
\frac{1}{2}
-
\frac{\sqrt{\pi} r  (1-r^2)^{-(d+1)/2} \Gamma(\frac{d+3}{2})}{2 \Gamma(\frac{d+2}{2}) (1+(d+1)r^2 {_2F_1}(1,\frac{d+2}{2};\frac{3}{2};r^2))}
\bigg)
\mathbb{I}[r\leq 1]
%$$
%and
%$$
%g(r)
%=
%\frac{1}{2}-
%\frac{\sqrt{\pi} r  \Gamma(\frac{d+3}{2})}{d \Gamma(\frac{d}{2}) ((1-r^2)^{(d+1)/2}+(d+1)r^2 {_2F_1}(\frac{1}{2},\frac{1-d}{2};\frac{3}{2};r^2))}
%$$
%$$
%g(r)
%=
%\frac{1}{2}-
%\frac{\sqrt{\pi} r  (d+1)\Gamma(\frac{d+1}{2})}{2d \Gamma(\frac{d}{2}) ((1-r^2)^{(d+1)/2}+(d+1)r^2 {_2F_1}(\frac{1}{2},\frac{1-d}{2};\frac{3}{2};r^2))}
%$$
%$$
%g(r)
%=
%\frac{1}{2}-
%\frac{\sqrt{\pi} r  (d+1)\Gamma(\frac{d+1}{2})}{2d \Gamma(\frac{d}{2}) ((1-r^2)^{(d+1)/2}+(d+1)r^2 (1-r^2)^{(d+1)/2} {_2F_1}(1,\frac{d}{2}+1;\frac{3}{2};r^2))}
%$$
%$$
%g(r)
%=
%\frac{1}{2}-
%\frac{\sqrt{\pi} r (1-r^2)^{-(d+1)/2} (d+1)\Gamma(\frac{d+1}{2})}{2d \Gamma(\frac{d}{2}) (1+(d+1)r^2 {_2F_1}(1,\frac{d}{2}+1;\frac{3}{2};r^2))}
%$$
%$$
%\ \textrm{and} \
%g(r)
%=
%\omega_{d-2}(r)
%=
%\frac{1}{2}-
%\frac{\sqrt{\pi} r (1-r^2)^{-(d+1)/2} \Gamma(\frac{d+3}{2})}{2 \Gamma(\frac{d+2}{2}) (1+(d+1)r^2 {_2F_1}(1,\frac{d+2}{2};\frac{3}{2};r^2))}
$$
and 
$g(r)
=
\omega_{d-2}(r)
$,
respectively, where~$\Gamma$ is the Euler Gamma function and~$_2F_1$ is the hypergeometric function. From affine invariance, these expressions agree with those obtained for~$d=1$ in~(\ref{expecunivexam1}) and~(\ref{expecunivexam2}), respectively. 

Our last example is a non-elliptical one. Consider the probability measure~$P_\alpha(\in\mathcal{P}_d)$ having independent standard (symmetric) $\alpha$-stable marginals, with~$1<\alpha\leq 2$. If~$\Zb=(Z_1,\ldots,Z_d)'$ has distribution~$P_\alpha$, then~$\ub'\Zb$ 
%has characteristic function~${\rm E}[\exp(it\ub'Z)]
%%={\rm E}[\exp(i(tu_1)Z_1]{\rm E}[\exp(i(tu_2)Z_2]
%=\prod_{j=1}^d \exp(-|tu_j|^\alpha)
%=\exp(-|t|^\alpha \|u\|_{\alpha}^\alpha)
%={\rm E}[\exp(it \|u\|_\alpha Z_1)]
%,
%$
%so that~$\ub'Z$ 
is equal in distribution to~$\|\ub\|_\alpha Z_1$, where we let~$\|\xb\|_\alpha^\alpha:=\sum_{j=1}^d |x_j|^\alpha$. Thus,~(\ref{DefinExpectileDepth2}) provides 
%This yields
%$
%{\rm E}[ |\ub'(\Zb-\zb)| \mathbb{I}[\ub'\Zb\leq \ub'\zb]]
%%=
%%{\rm E}[ | \|u\|_\alpha Z_1 - \ub'z | \mathbb{I}[ \|u\|_\alpha Z_1 \leq \ub'z]]
%=
%\|u\|_\alpha {\rm E}[ ( Z_1 - v'z ) \mathbb{I}[  Z_1 \leq v'z]]
%,
%$
%with~$v:=u/\|u\|_\alpha$. 
$$
E\!D(\zb,P_\alpha)
%=
%\min_{\ub\in\mathcal{S}^{d-1}}
%\frac{{\rm E}[ |\ub'(\Zb-\zb)| \mathbb{I}[\ub'\Zb\leq \ub'\zb]]}{{\rm E}[|\ub'(\Zb-\zb)|]}
=
\min_{v\in\mathcal{S}_\alpha^{d-1}}
\frac{{\rm E}[ |Z_1-\vb'\zb| \mathbb{I}[Z_1 \leq \vb'\zb]]}{{\rm E}[ |Z_1-\vb'\zb|]}
,
$$
where~$\mathcal{S}^{d-1}_\alpha:=\{\vb\in\R^d: \|\vb\|_\alpha=1\}$ is the unit $L_\alpha$-sphere. Theorem~\ref{Theorlemjones}(iv) implies that the minimum is achieved when~$\vb'\zb$ takes its minimal value~$- \|\zb\|_{\beta}$, where~$\beta=\alpha/(\alpha-1)$ is the conjugate exponent to~$\alpha$; see Lemma~A.1 in~\cite{CheTyl2004}. 
%Assume without loss of generality that~$\zb_1,\zb_2\leq 0$. Then
%$$
%v'z = v_1 z_1 +v_2 z_2 = v_1 z_1 + {\rm Sgn}(v_2) ( 1 - |v_1|^\alpha)^{1/\alpha} z_2   
%$$
%should be minimized with non-negative values of~$v_1,v_2$. Now, the derivative of
%$
%v'z = v_1 z_1 + ( 1 - v_1^\alpha)^{1/\alpha} z_2   
%$
%with respect to~$v_1$ is~$z_1 - ( 1 - v_1^\alpha)^{(1/\alpha)-1} v_1^{\alpha-1} z_2$ and the second derivative is~$-(\alpha-1) v_1^{\alpha-2} (1-v_1^\alpha)^{(1/\alpha)-2} z_2\geq 0$, so that the function is convex in~$v_1$ over~$[0,1]$. Thus the minimum is indeed reached at any critical point, that is the solution of
%$$
% ( 1 - v_1^\alpha)^{1-(1/\alpha)} z_1 = v_1^{\alpha-1} z_2 
%,
%$$
%that is, of 
%$
%v_1^{\alpha} = |z_1|^{\beta}/(|z_1|^{\beta}+|z_2|^{\beta}) 
%,
%$ 
%where~$\beta:=\alpha/(\alpha-1)$ is the conjugate exponent of~$\alpha$. This provides
%$$
%v_{\rm min}
%=
%\bigg(
%\bigg(
%\frac{|z_1|^{\beta}}{|z_1|^{\beta}+|z_2|^{\beta}} 
%\bigg)^{1/\alpha}
%,
%\bigg(
%\frac{|z_2|^{\beta}}{|z_1|^{\beta}+|z_2|^{\beta}} 
%\bigg)^{1/\alpha}
%\bigg)
%=
%\frac{\big(
%|z_1|^{\beta/\alpha}
%,
%|z_2|^{\beta/\alpha}
%\big)}
%{(|z_1|^{\beta}+|z_2|^{\beta})^{1/\alpha}}
%,
%$$ 
%which leads to
%$$
%v_{\rm min}'z 
%=
%-
%\frac{
%|z_1|^{1+(\beta/\alpha)}
%+
%|z_2|^{1+(\beta/\alpha)}
%}
%{(|z_1|^{\beta}+|z_2|^{\beta})^{1/\alpha}}
%=
%-
%\frac{
%|z_1|^{\beta}
%+
%|z_2|^{\beta}
%}
%{(|z_1|^{\beta}+|z_2|^{\beta})^{1/\alpha}}
%%=
%%-
%%(|z_1|^{\beta}+|z_2|^{\beta})^{1/\beta}
%=
%-\|\zb\|_{\beta}
%.
%$$
Denoting as~$f_\alpha$ the marginal density of~$P_\alpha$, this yields
$$
E\!D(\zb,P_\alpha)
=
- \frac{{\rm E}[ (Z_1+\|\zb\|_{\beta} ) \mathbb{I}[  Z_1 \leq -\|\zb\|_{\beta}]]}{{\rm E}[ |Z_1+\|\zb\|_{\beta}|]}
=
-
\frac{\int_{-\infty}^0 x f_\alpha(x-\|\zb\|_{\beta})\, dx}{\int_{-\infty}^\infty |x| f_\alpha(x-\|\zb\|_{\beta})\, dx}
,
$$
which shows that expectile depth regions are concentric $L_{\beta}$-balls. For~$\alpha=2$, these results agree with those obtained in the Gaussian case above.

\subsection{An interesting monotonicity property}

The sample \mbox{M-depth} regions $R^{\rho}_{\alpha}(P_n)$ can be computed by replacing the intersection in Definition~\ref{definMregions} with an intersection over finitely many directions~$\ub_\ell$, $\ell=1,\ldots,L$, with~$L$ large; see Section~\ref{secproposedmultiquantile}. Some applications, however, do not require computing depth regions but rather the depth of a given location~$\zb$ only. An important example is supervised classification through the max-depth approach; see \cite{GhoCha2005B} or \cite{Lietal2012}. While the halfspace \mbox{M-depth} of~$\zb$ can in principle be obtained from the depth regions (recall that~$
%\label{deviceregdepth}	
M\!D^\rho(\zb,P)
=
\sup
\{
\alpha> 0 : \zb\in R^{\rho}_{\alpha}(P) 
\}
$), it will be much more efficient in such applications to compute~$M\!D^\rho(\zb,P)$ through the alternative expression in~(\ref{mequ}). Recall that, for the halfspace Tukey and halfspace expectile depths, this alternative expression reduces to 
$$
H\!D(\zb,P) 
=
\inf_{\ub\in\mathcal{S}^{d-1}} h_{\zb}(\ub)
\quad
\textrm{ and } 
\quad 
E\!D(\zb,P) 
=
\min_{\ub\in\mathcal{S}^{d-1}} e_{\zb}(\ub)
,
$$
respectively, where we let
\begin{equation}
	\label{EDcomput}
h_{\zb}(\ub):=
P[\ub'\Zb\leq \ub'\zb] 
\
\textrm{ and } 
\ 
e_{\zb}(\ub):=\frac{{\rm E}[ |\ub'(\Zb-\zb)| \mathbb{I}[\ub'(\Zb-\zb)\leq 0]]}{{\rm E}[|\ub'(\Zb-\zb)|]}
\cdot
\end{equation}
There is a vast literature dedicated to the evaluation of halfspace Tukey depth and it is definitely beyond the scope of the paper to thoroughly discuss the computational aspects of our expectile depth. Yet we state a monotonicity property that is extremely promising for expectile depth evaluation.

\begin{Theor}
\label{Theordepthuniqueu}
Fix~$P\in\mathcal{P}_d$ and~$\zb\in\R^d$ such that~$E\!D(\zb,P)>0$. Assume that~$P[\Pi\setminus \{\zb\}]=0$ for any hyperplane~$\Pi$ containing~$\zb$. Fix a great circle~$\mathcal{G}$ of~$\mathcal{S}^{d-1}$ and let~$\ub_0$ be an arbitrary minimizer of~$e_{\zb}(\cdot)$ on~$\mathcal{G}$. Let~$\ub_t$, 
\linebreak
$t\in[0,\pi]$, be a path on~$\mathcal{G}$ from~$\ub_0$ to~$-\ub_0$. Then, there exist~$t_a, t_b$ with~$0\leq t_a\leq t_b\leq \pi$ such that~$t\mapsto e_{\zb}(\ub_t)$ is constant over~$[0,t_a]$, admits a strictly positive 
%right-
derivative at any~$t\in (t_a,t_b)$ $($hence is strictly increasing over~$[t_a,t_b])$, and is constant over~$[t_b,\pi]$. Moreover, letting~$Z$ be a random $d$-vector with distribution~$P$, the minimal direction~$\ub_0$ is such that~$\zb$ belongs to the line segment with endpoints~${\rm E}[ \Zb \mathbb{I}[\ub_0'\Zb\leq \ub_0'\zb]]$ and~${\rm E}[ \Zb \mathbb{I}[\ub_0'\Zb\geq \ub_0'\zb]]$.
\end{Theor}

\begin{figure}
\center
\includegraphics[width=\textwidth]{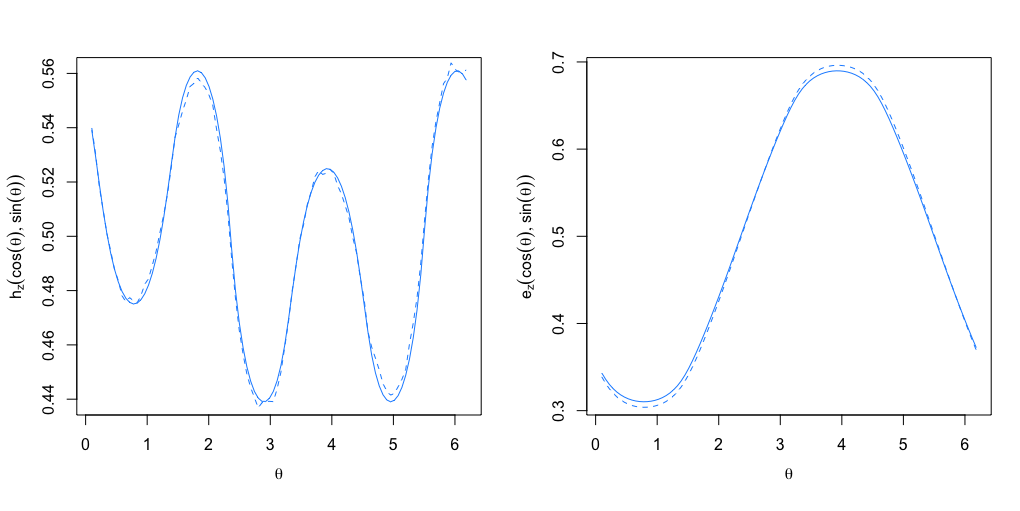}
\caption{(Left:) plots of the Tukey depth outlyingness~$\theta\mapsto h_\zb((\cos \theta,\sin \theta)')$ in~(\ref{EDcomput}),   
% over~$[0,2\pi]$, 
with~$\zb=(.8,.8)'$, for the probability measure~$P$ over~$\R^2$ whose marginals are independent exponentials with mean one (solid curve) and for the empirical probability measure associated with a random sample of size~$n=10,\!000$ from~$P$ (dashed curve). (Right:) the corresponding plots of the expectile depth outlyingness~$\theta\mapsto e_\zb((\cos \theta,\sin \theta)')$; see also~(\ref{EDcomput}).}
\label{Fig3}
\end{figure}

As an example, we consider the probability measure~$P$ over~$\R^2$ whose marginals are independent exponentials with mean one. Figure~\ref{Fig3} draws, for~$\zb=(.8,.8)'$, the plot of~$\theta\mapsto h_\zb((\cos \theta,\sin \theta)')$ and~$\theta\mapsto e_\zb((\cos \theta,\sin \theta)')$ over~$[0,2\pi]$. Clearly, this illustrates the monotonicity property in Theorem~\ref{Theordepthuniqueu} and shows that monotonicity may fail for Tukey depth. The figure also plots these functions evaluated on the empirical probability measure~$P_n$ associated with a random sample of size~$n=10,\!000$ from~$P$, which shows that monotonicity extends to the sample case. Obviously, this monotonicity opens the door to computation of expectile depth through standard optimization algorithms, while the lack of monotonicty for Tukey depth could lead such algorithms to return local minimizers only.

%%%%%%%%%%%%%%%%%%%%%%%%%%%%%%%%%%%%%%%%%%%%%%%%%%

%%%%%%%%%%%%%%%%%%%%%%%%%%%%%%

\section{Multivariate expectile risks} 
\label{secrisks}

 The risk of a collection of financial assets is typically assessed by aggregating these assets, using their monetary values, into a combined random univariate portfolio~$Z$. It is then sufficient to consider univariate risk measures~$\varrho(Z)$; see \cite{Artetal1999} and \cite{Del2002}. More and more often, however, the focus is on the more realistic situation where the risky portfolio is a random $d$-vector whose components relate to different security markets. In such a context, liquidity problems and/or transaction costs between the various security markets typically prevent investors from aggregating their portfolio into a univariate portfolio (\citealp{Jouetal2004}). This calls for multivariate risk measures~$\varrho(\Zb)$, where~$\Zb$ is a random $d$-vector.  
 
Extensions of the axiomatic foundation for \emph{coherent} univariate risk measures to the $d$-variate framework have been studied in \cite{Jouetal2004} and \cite{CasMol2007}. Such extensions usually involve set-valued risk measures, as in the following definition (we restrict here to bounded random vectors as in \citealp{Jouetal2004}, but the extension to the general case could be achieved as in \citealp{Del2002}).

\begin{Defin}
\label{defincoher}
Let~$L^{\infty}_d$ be the set of (essentially) bounded random \mbox{$d$-vectors} and~$\mathcal{B}_d$ be the Borel sigma-algebra on~$\R^d$. Then a \emph{coherent $d$-variate risk measure} is a function $R:L^{\infty}_d\to \mathcal{B}_d$ satisfying the following properties:
%\begin{itemize}
%\item[(i)] 
(i) 
(translation invariance:)  $R(\Zb+\zb) = R(\Zb)+\zb$ for any $\Zb\in L^{\infty}_d$ and~$\zb \in\R^d$;
%\item[(ii)] 
(ii) 
(positive homogeneity:) $R(\lambda \Zb) = \lambda R(\Zb)$ for any~$\Zb\in L^{\infty}_d$ and~$\lambda>0$;
%\item[(iii)] 
(iii) 
(monotonicity:) 
if $\Xb\leq \Yb$ almost surely in the componentwise sense, then
$R(\Yb) \subset  R(\Xb)  \oplus \R^{d}_+$ and 
$R(\Xb)  \subset  R(\Yb)  \oplus \R^{d}_-$,
where~$\oplus$ denotes the Minkowski sum and where we let~$\R^d_\pm:=\{\xb\in\R^d\!:\pm x_1\geq 0,\ldots,\pm x_d\geq 0\}$;
%\item[(iv)] 
(iv) 
(subadditivity:)
$R(\Xb+\Yb) \subset R(\Xb) \oplus  R(\Yb)$ for any~$\Xb,\Yb\in L^{\infty}_d$;
%\item[(v)] 
(v) 
(connectedness/closedness:)
$R(\Xb)$ is connected and closed for any~$\Xb\in L^{\infty}_d$.
%\end{itemize}
%	
\end{Defin} 

In the univariate case, such coherent set-valued risk measures can be obtained as~$R(Z)=[-\varrho(Z),\infty)$, where~$\varrho(Z)$ is a real-valued coherent risk measure in the sense of \cite{Artetal1999} and \cite{Del2002}; see Remark~2.2 in \cite{Jouetal2004}. For the most classical risk measure, namely the \!\textit{Value at Risk}, the resulting set is $R(Z)=[-\textrm{VaR}_{\alpha}(Z),\infty)$, where~$-\textrm{VaR}_{\alpha}(Z)=q_{\alpha}(Z)$ is the standard $\alpha$-quantile of~$Z$. The sign convention in~$\textrm{VaR}_{\alpha}(Z)$ corresponds to an implicit specification of the positive direction~$u=1$, which associates a positive risk measure with the typically negative profit---that is, loss---$q_\alpha(Z)$ obtained for small values of~$\alpha$. 

In this univariate setting, M-quantiles have recently received a lot of attention since the resulting risk measures share the important property of {\it elicitability}, which corresponds to the existence of a natural backtesting methodology (\citealp{Gne2011}). In this framework, expectiles play a special role as they are the only M-quantiles providing coherent risk measures (\citealp{Beletal2014}). Actually, expectiles define the only coherent risk measure that is also elicitable (\citealp{Zie2016}). In the $d$-variate case, a natural M-quantile set-valued risk measure is given by our M-quantile halfspace~$H^{\rho}_{\alpha,\ub}(\Zb)$ in Definition~\ref{directMquant} (in this section, $H^{\rho}_{\alpha,\ub}(\Zb)$, $R^{\rho}_{\alpha}(\Zb)$,\,\ldots respectively stand for~$H^{\rho}_{\alpha,\ub}(P)$, $R^{\rho}_{\alpha}(P)$, \ldots, where~$P$ is the distribution of~$Z$). For~$d=1$, $\rho(t)=|t|$ and the positive direction~$u=1$, this M-quantile set-valued risk measure reduces to the risk measure~$[-\textrm{VaR}_{\alpha}(Z),\infty)$ above, which, as already mentioned, also relies on the choice of a positive direction. For~$d>1$, it is similarly natural to restrict to ``positive" directions~$\ub$, that is, to~$\ub\in\mathcal{S}^{d-1}_+:=\mathcal{S}^{d-1}\cap\R^d_+$.

Now, already for~$d=1$, the VaR risk measure fails to be subadditive in general (\citealp{Ace2002}). It is also often criticized for its insensitivity to extreme losses, since it depends on the frequency of tail losses but not on their severity. Denoting as~$e_\alpha(Z)$ the order-$\alpha$ expectile of~$Z$, the expectile risk measure~$R(Z)=[e_{\alpha}^{\rho}(Z),\infty)$, with~$\alpha\in(0,\frac{1}{2}]$, improves over VaR on both accounts since it is coherent (\citealp{Beletal2014}) and depends on the severity of tail losses (\citealp{Kuaetal2009}). Our expectile $d$-variate risk measure, namely the halfspace~$H^{\rho}_{\alpha,\ub}(\Zb)$ based on~$\rho(t)=t^2$, extends this univariate expectile risk measure to the $d$-variate setup and, quite nicely, turns out to be coherent for any~$\alpha\in(0,\frac{1}{2}]$ and any direction $\ub\in\mathcal{S}^{d-1}_+$: %It is therefore the only multivariate M-quantile risk measure $H^{\rho}_{\alpha,\ub}$ that is coherent since, already in the one-dimensional case, $H^{\rho}_{\alpha,\ub}\equiv[\theta_{\alpha}^{\rho},\infty)$ is coherent only for the quadratic loss function $\rho(t)=t^2$.
since connectedness/closedness holds trivially ($H^{\rho}_{\alpha,\ub}(Z)$ is a closed halfspace) and since translation invariance and positive homogeneity directly follow from Theorem~\ref{Theoraffequiv}, we focus on monotonicity and subadditivity (see Definition~\ref{defincoher}) and further cover some other properties from \cite{DycMos2011}.

\begin{Theor}
\label{Theoraxioms}
Fix~$\rho(t)=t^2$ and
let~$\Xb,\Yb$ be random $d$-vectors with respective distributions~$P,Q$ in~$\mathcal{P}_d^\rho$. Then, we have the following properties: (i) (monotonicity) if $
\Xb\leq \Yb$ almost surely in a componentwise sense, then
$
H^{\rho}_{\alpha,\ub}(\Yb)\subset H^{\rho}_{\alpha,\ub}(\Xb) \oplus \R^d_+
$
and
$
H^{\rho}_{\alpha,\ub}(\Xb)\subset H^{\rho}_{\alpha,\ub}(\Yb) \oplus \R^d_-
$ for any~$\alpha\in(0,1)$ and~$\ub\in \mathcal{S}^{d-1}_+$;
(ii) (subadditivity) for any~$\alpha\in(0,\frac{1}{2}]$ and~$\ub\in\mathcal{S}^{d-1}$,
$
H^{\rho}_{\alpha,\ub}(\Xb+\Yb) 
\subset
H^{\rho}_{\alpha,\ub}(\Xb) \oplus H^{\rho}_{\alpha,\ub}(\Yb) 
$;
(iii) (superadditivity) for any~$\alpha\in[\frac{1}{2},1)$ and~$\ub\in\mathcal{S}^{d-1}$,
$
H^{\rho}_{\alpha,\ub}(\Xb) \oplus H^{\rho}_{\alpha,\ub}(\Yb) 
\subset
H^{\rho}_{\alpha,\ub}(\Xb+\Yb) 
$;
(iv) (nestedness:) for any~$\ub\in\mathcal{S}^{d-1}$, %the mapping 
$\alpha\mapsto H^{\rho}_{\alpha,\ub}(\Xb)$ is non-increasing with respect to inclusion.
\end{Theor}

\begin{figure}[htbp]   
  \vspace*{-0.4 cm}   
    \begin{center}            
        \includegraphics[width=\textwidth]{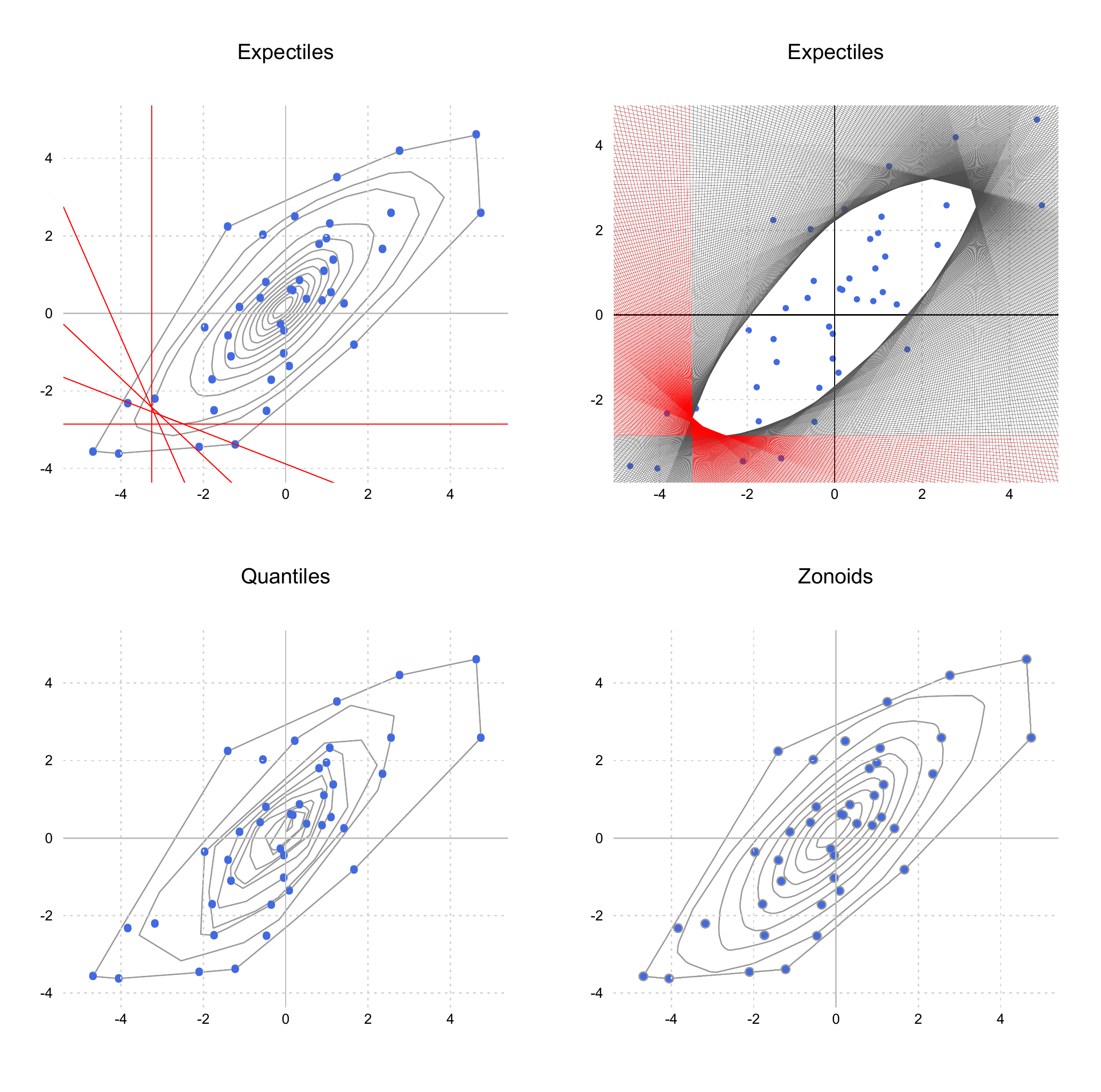}
     \end{center}
     \vspace*{-.8cm}
\caption{(Top left:) boundary hyperplanes of the expectile halfspaces~$H^{\rho}_{\alpha,\ub}(P_n)$, for~$\alpha=.02$ and~$\ub=(\cos \frac{\ell\pi}{8},\sin\frac{\ell\pi}{8})'$ with~$\ell=0,1,2,3,4$, along with the expectile depth regions~$R^{\rho}_\alpha(P_n)$ associated with~$\alpha=.0001, .01, .02, .06, \ldots, .38, .42$; as in section~\ref{secproposedmultiquantile}, each such region was computed from~$L=500$ equispaced directions~$u\in\mathcal{S}^{1}$. 
(Top right): the boundary hyperplanes of the  $500$ expectile halfspaces that led to the construction of~$R^{\rho}_\alpha(P_n)$ with~$\alpha=.02$; hyperplanes associated with (positive) directions~$\ub\in\mathcal{S}^{1}_+$ are drawn in red.
(Bottom left): halfspace Tukey depth regions of order~$\alpha=.02, .06, \ldots, .38, .42$.
(Bottom right:) zonoid depth regions of order~$\alpha=.01, .1, .2,\ldots, .8, .9$.
}
\label{Fig4}
\end{figure}

In order to illustrate these $d$-variate M-quantile risk measures, we briefly consider the daily returns on \mbox{Intel~Corp.} and Adobe \mbox{Systems~Inc.} shares in May--June 2008. We chose the institutions, frequency of data and time horizon exactly as in \cite{DycMos2011}. The data, kindly sent to us by \mbox{Prof.} Rainer Dyckerhoff, 
were taken from the historical stock market database at the Center for Research in Security Prices (CRSP), University of Chicago. Figure~\ref{Fig4} shows the resulting $n = 42$ bivariate observations along with some of the corresponding expectile risk measures~$H^{\rho}_{\alpha,u}(P_n)$ (more precisely, the figure only displays their boundary hyperplanes) and some expectile depth regions~$R^{\rho}_{\alpha}(P_n)$. We also provide there a few halfspace Tukey depth regions and zonoid depth regions. For~$d=1$, the latter are related to \emph{expected shorftall}, hence are also connected to risk measures. However, while zonoid regions formally are coherent risk measures (\citealp{CasMol2007,DycMos2011}), these regions, like their quantile and expectile counterparts, can hardly be interpreted in terms of riskiness as such centrality regions not only trim joint returns with extreme losses but also those with extreme profits (accordingly, a univariate zonoid depth region is not an interval of the form~$[-\varrho(Z),\infty)$ but rather a compact interval). In contrast, our M-quantile risk measures~$H^{\rho}_{\alpha,\ub}(P_n)$, $\ub\in\mathcal{S}^{d-1}_+$, protect against adverse joint returns only. They also offer an intuitive interpretation for the multivariate risk in the sense that the required capital reserve should cover any loss associated with joint returns inside~$H^{\rho}_{\alpha,\ub}(P_n)$, that is, above the hyperplane $\pi^{\rho}_{\alpha,\ub}(P_n)$. For these risks, the choice of a suitable security level~$\alpha$ and direction~$\ub\in\mathcal{S}^{d-1}_+$ is a decision that should be made by risk managers and regulators. Other $d$-variate set-valued risk measures that trim unfavorable returns only yet do not require the choice of a direction~$u$, are the upper envelopes~$\cap_{\ub\in\mathcal{S}^{d-1}_+} H^{\rho}_{\alpha,\ub}(P_n)$ of our directional M-quantile risk measures.

%%%%%%%%%%%%%%%%%%%%%%%%%%%%%%%%%%%%%%%%%%%%%%%%%%

%%%%%%%%%%%%%%%%%%%%%%%%%%%%%%

\section{Multiple-output expectile regression} 
\label{secregres}

We now consider the mul\-tiple-output regression framework involving a $d$-vector~$Y$ of responses and a $p$-vector~$X$ of (random) covariates. For any possible value~$x$ of~$X$, denote as~$P_x$ the conditional distribution of~$Y$ given~$X=x$. Our interest then lies in the conditional M-quantile halfspaces and regions
$$
H^\rho_{\alpha,u,x}=H^\rho_{\alpha,u}(P_x)
\ \
\textrm{ and }
\ \
R^\rho_{\alpha,x}=R^\rho_{\alpha}(P_x),
$$
with~$\alpha\in(0,1)$ and~$u\in\mathcal{S}^{d-1}$. If a random sample~$(X_1,Y_1),\ldots,(X_n,Y_n)$ is available, then one may consider the estimates
\begin{equation}
	\label{samplecondregions}
\hat{H}^{\rho(n)}_{\alpha,u,x}
:=
\big\{ \yb\in\R^d : \ub'\yb \geq \theta^{\rho(n)}_{\alpha,u,x} \big\}
\
\textrm{ and }
\
\hat{R}^{\rho(n)}_{\alpha,x}
:=
\cap_{u\in \mathcal{S}^{d-1}}
\hat{H}^{\rho(n)}_{\alpha,u,x}
,
\end{equation}
%hyperplanes~$\pi_{\alpha,u,x}=\pi^\rho_{\alpha,u}(P_x)$,
where~$\theta^{\rho(n)}_{\alpha,u,x}$ is the estimate of~$\theta^{\rho}_{\alpha,u}(P^{u'Y|[X=x]})$ obtained from a single-output, linear or nonparametric, regression using the responses~$u'Y_1,\ldots,u'Y_n$ and covariates~$X_1,\ldots,X_n$ (in the examples below, that focus on~$d=2$, the intersection in~(\ref{samplecondregions}) was replaced with an intersection over~$L=200$ equispaced directions in~$\mathcal{S}^1$). For expectiles, single-output linear and nonparametric regression can respectively be performed via the functions \texttt{expectreg.ls} and \texttt{expectreg.boost} from the R package {\it expectreg} (nonparametric regression here is thus based on the expectile boosting approach from \citealp{SobKne2012}). Multiple-output quantile regression can be achieved in the same way, by performing single-output linear quantile regression (via the function \texttt{rq} in the R package {\it quantreg}) or single-output nonparametric quantile regression (via, e.g., the function \texttt{cobs} in the R package {\it cobs}, which relies on the popular quantile smoothing spline approach). Whenever we use \texttt{expectreg.boost} and \texttt{cobs} below, it is with the corresponding default automatic selection of smoothing parameters.

%%%%%%%%%%%%%%%%%%%

\subsection{Simulated data illustration} 
\label{sec71}

\begin{figure}[htbp]
  \vspace*{-0.4 cm} 
  \begin{center}$
\begin{array}{cc}
\hspace{-5mm}
\includegraphics[width=6.901cm]{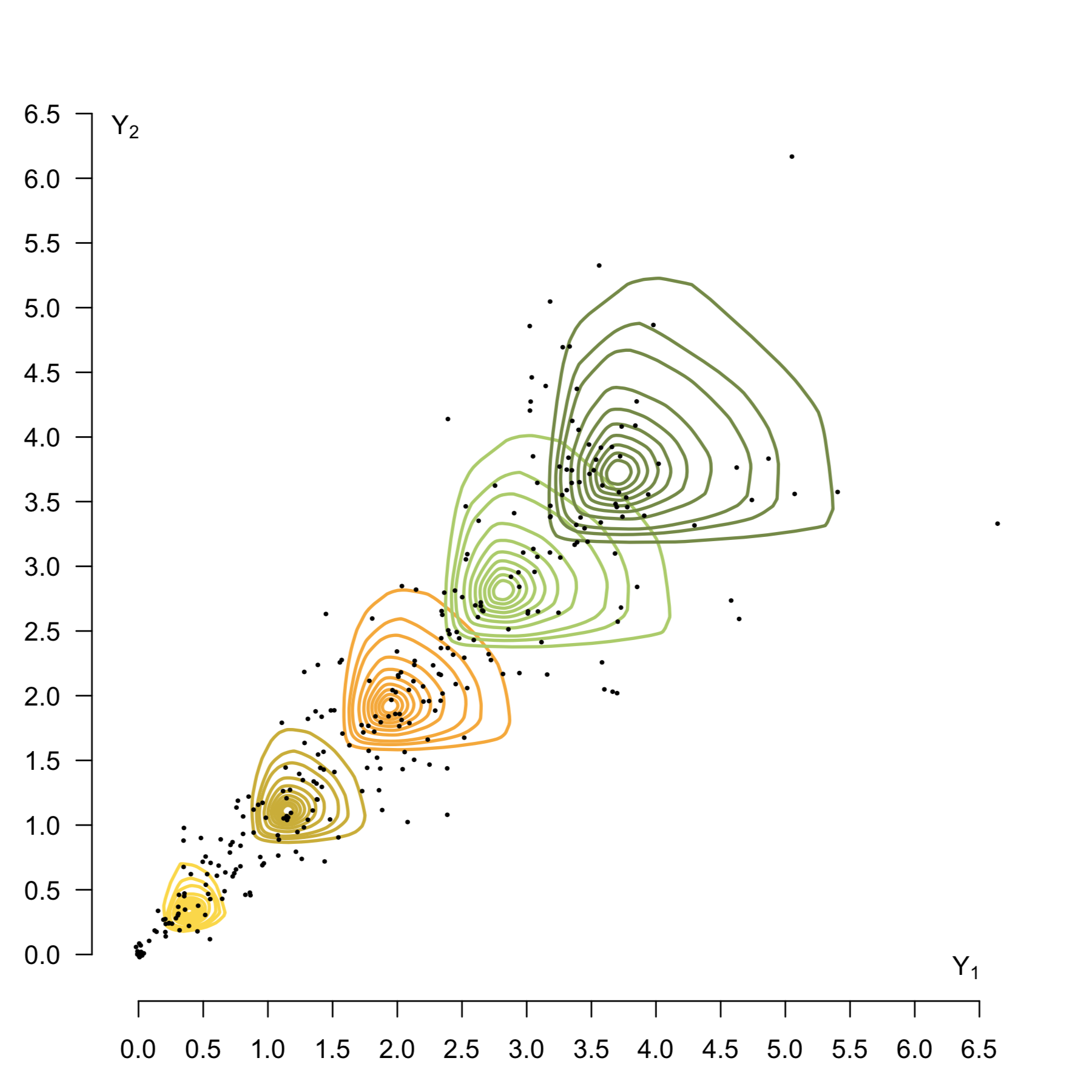}
&
\hspace{-7mm}
\includegraphics[width=6.901cm]{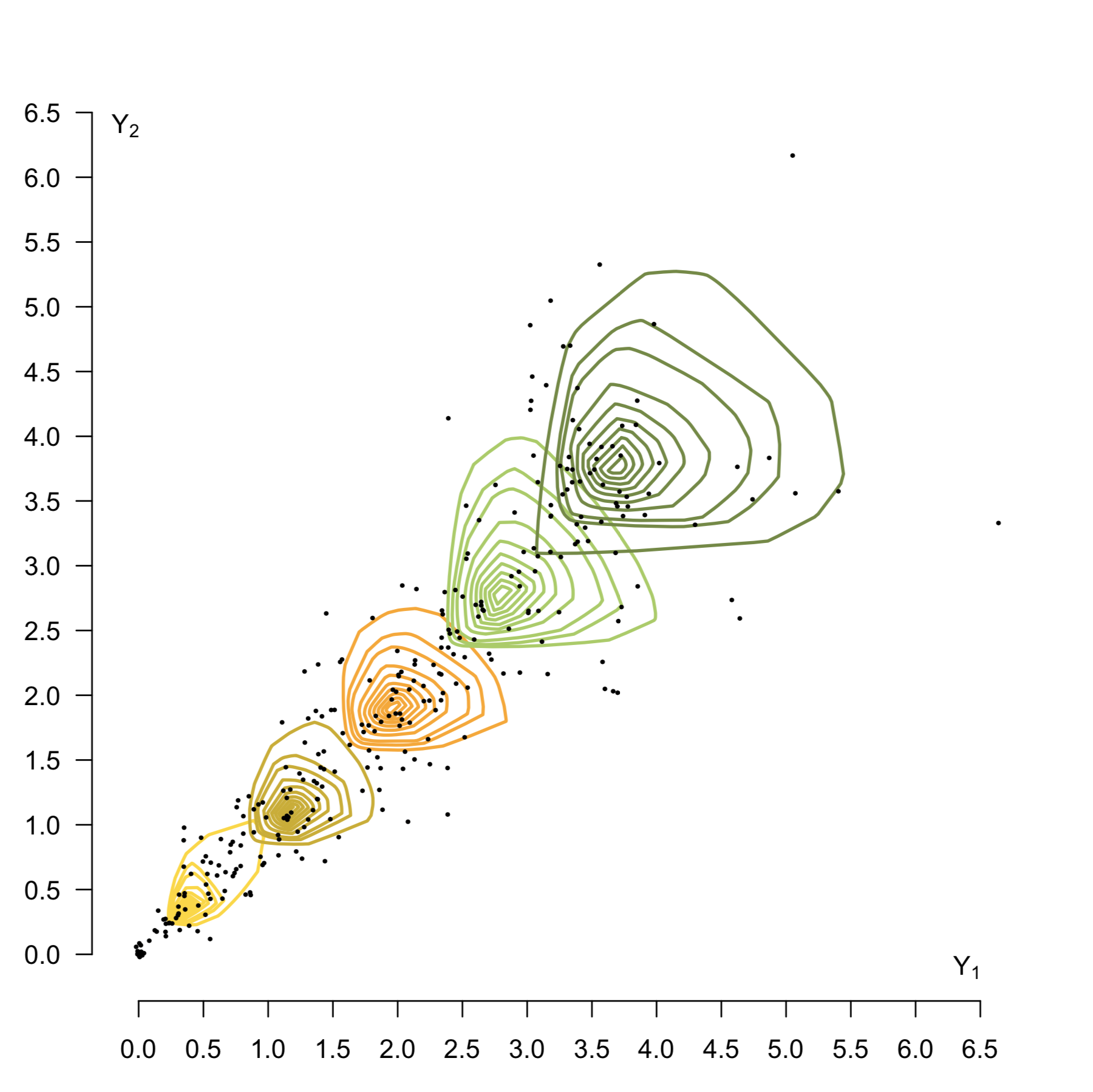}
\\
\hspace{-5mm}
\includegraphics[width=6.901cm]{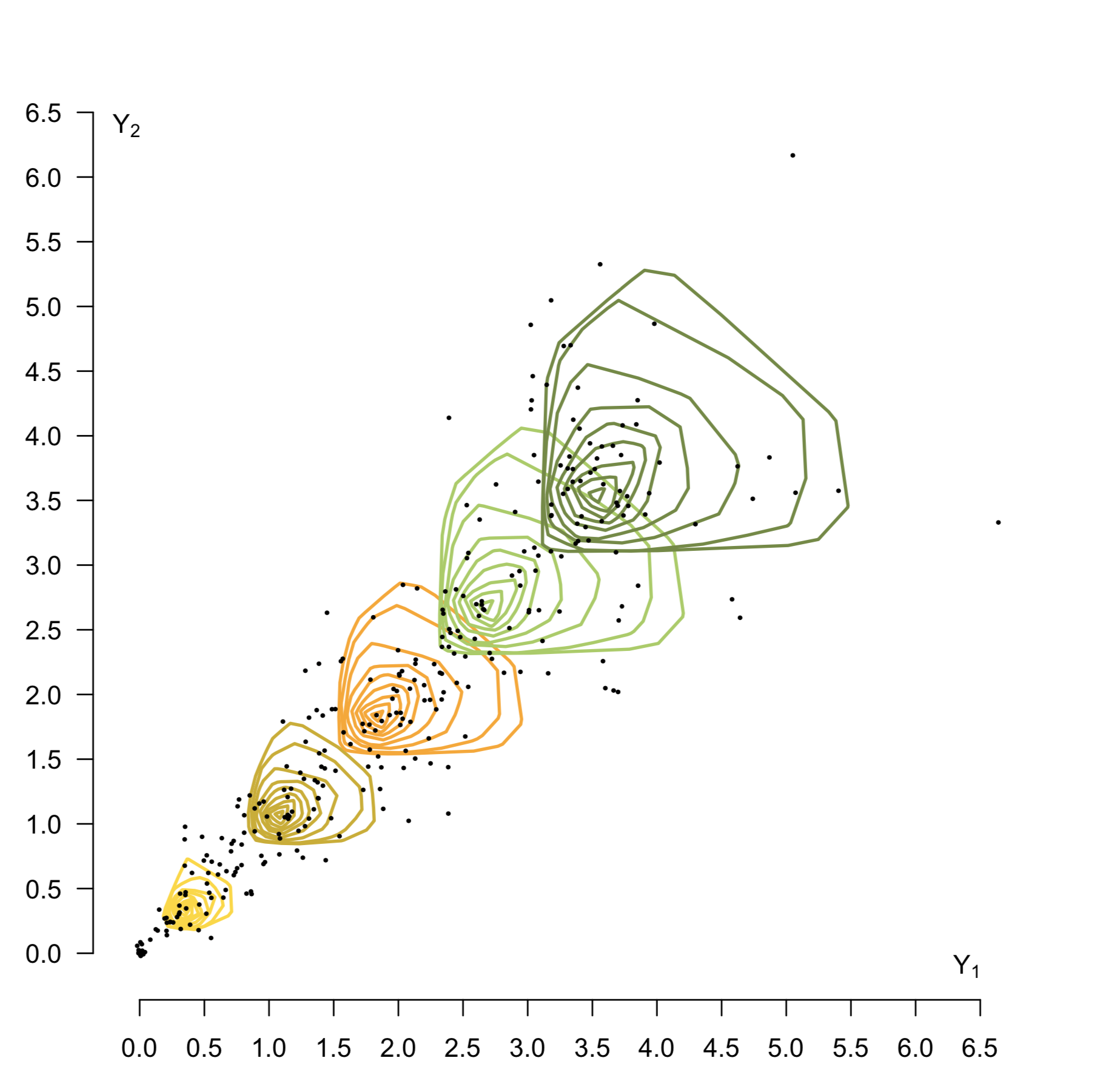}
&
\hspace{-7mm}
\includegraphics[width=6.901cm]{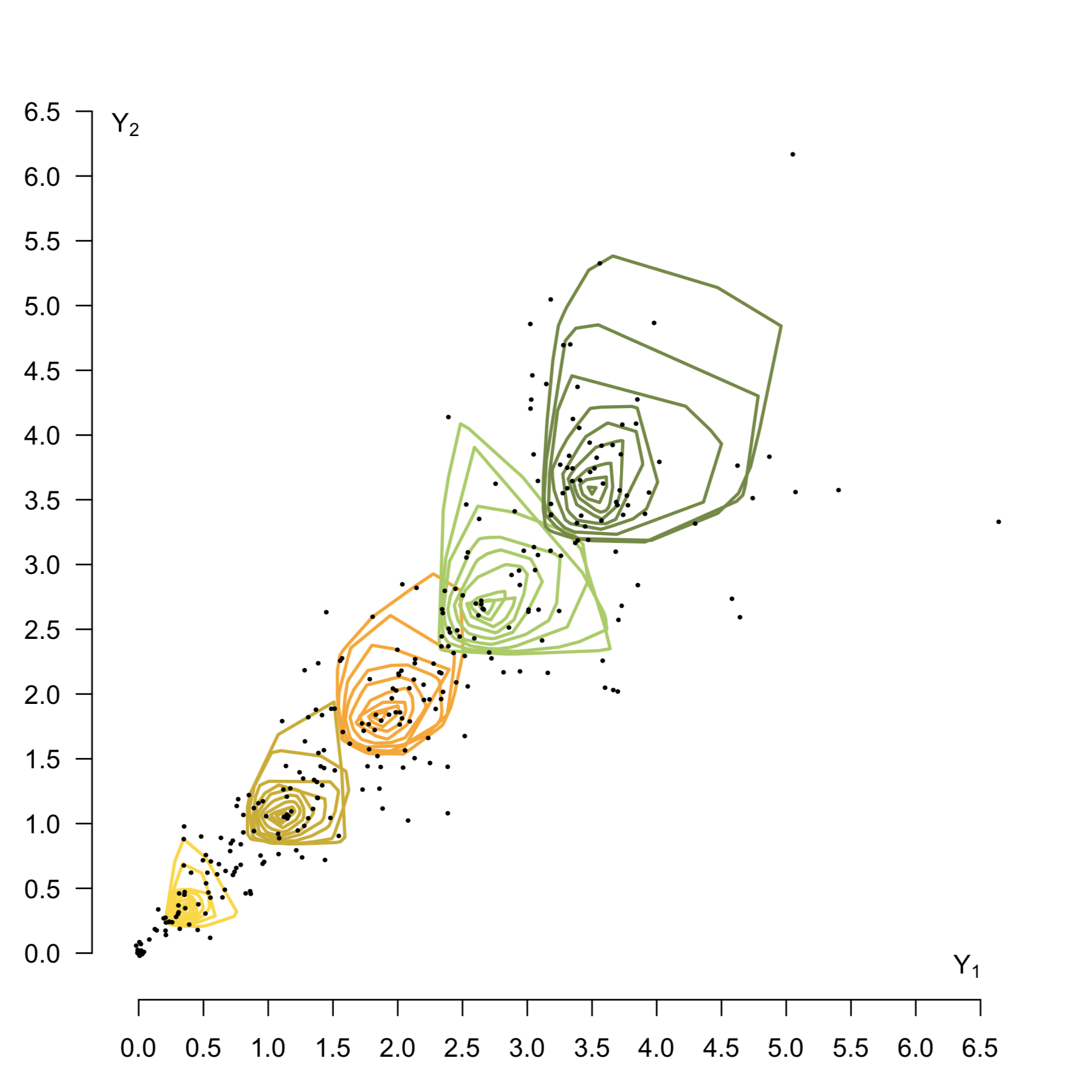}
\end{array}$
\end{center}
     \vspace*{1mm}
\caption{(Top:) conditional expectile contours~$\partial \hat{R}^{\rho(n)}_{\alpha,x}$, for $\alpha\in\{.01,.03,.05, .10,.15,\ldots, .40\}$ and for values of~$x$ that are the $10\%$ (yellow), $30\%$ (brown), $50\%$ (orange),
$70\%$ (light green) and $90\%$ (dark green) empirical quantiles of $X_1,\ldots,X_n$, obtained by applying a linear regression method (left) or a nonparametric regression method (right) to a random sample of size~$n=300$ from the linear regression model in~(\ref{simumodel}). (Bottom:) conditional quantile contours associated with the same values of~$\alpha$ (but~$.01$) and the same values of~$x$. Again, both linear regression (left) and nonparametric regression (right) are considered; see Section~\ref{sec71} for details. Bivariate responses~$(Y_{i1},Y_{i2})$, $i=1,\ldots,n$, are shown in black.
}
\label{Fig5}
\end{figure}

To illustrate these multiple-output regression methods on simulated data, we generated a random sample of size~$n=300$ from the heteroscedastic linear regression model
\begin{equation}
\label{simumodel}
\bigg(
\!\!
\begin{array}{c}
	Y_1 \\[.5mm] 
	Y_2
\end{array}
\!\!
\bigg)
=
4
\bigg(
\!\!
\begin{array}{c}
	X \\[.5mm] 
	X
\end{array}
\!\!
\bigg)
 + 
 \sqrt{\frac{X}{3}}
% 2 \sqrt{X}
\bigg(
\!\!
\begin{array}{c}
	\varepsilon_1 \\[.5mm] 
	\varepsilon_2
\end{array}
\!\!
\bigg)
 ,
\end{equation}
where the covariate $X$ is uniform over~$[0,1]$, $\varepsilon_1+1,\varepsilon_2+1$ are exponential with mean one,
%uniform over~$[-.5, .5]$, 
and~$X,\varepsilon_1,\varepsilon_2$ are mutually independent. For 
 several orders~$\alpha$ and several values of~$x$, we evaluated
\vspace{-.6mm}
 the conditional quantile and expectile regions~$\hat{R}^{\rho(n)}_{\alpha,x}$, in each case both from the corresponding linear and nonparametric regression methods above. The resulting contours are provided in Figure~\ref{Fig5}. Both expectile and quantile methods capture the trend and heteroscedasticity. Unlike quantiles, however, expectiles provide linear and nonparametric regression fits that are very similar (which is desirable, since the underlying model is a linear one). Expectile contours are also smoother than the quantile ones. Inner expectile contours, that do not have the same location as their quantile counterparts, are easier to interpret since they relate to conditional means of the marginal responses (inner quantile contours refer to the Tukey median, which is not directly related to marginal medians). Finally, it should be noted that, unlike expectile contours, several quantile contours associated with a common value of~$x$ do unpleasantly cross, which is incompatible with what occurs at the population level.

%%%%%%%%%%%%%%%%%%%%%%%%%%%%%%%%%%%%%%%%%%%%%%%%%%

%%%%%%%%%%%%%%%%%%%%%%%%%%%%%%%%%%%%%%%%%%%%%%%%%%

%%%%%%%%%%%%%%%%%%%

\subsection{Real data illustration}
\label{sec72}

We now conduct multiple-output expectile regression to investigate risk factors for early childhood malnutrition in India. Prior studies typically focused on children's {\it height} as an indicator of nutritional status (\citealp{Koe2011,Fenetal2011}). Given that the prevalence of underweighted children in India is among the highest worldwide, we consider here determinants of children's {\it weight} ($Y_1$; in kilograms) and {\it height} ($Y_2$; in centimeters) simultaneously. We use a selected sample of 37,623 observations, coming from the 2005/2006 Demographic and Health Survey (DHS) conducted in India for children under five years of age. Since a thorough case study is beyond the scope of this paper, we restrict to assessing the separate effects of the following covariates on the response~$Y=(Y_1,Y_2)'$: (a) the child's age (in months), (b) the duration of breastfeeding (in months), and (c) the mother's Body Mass Index (defined as~$\textrm{BMI}:=\textrm{weight}/\textrm{height}^2$, in $\textrm{kilograms}/\textrm{meters}^2$). \cite{Koe2011} investigated the additive effects of these covariates on low levels of the single response {\it height} through a quantile regression with small~$\alpha$; see \cite{Fenetal2011} for a similar quantile regression analysis of this dataset. 

\begin{figure}[htbp]
  \vspace*{-0.4 cm} 
  \begin{center}$
\begin{array}{cc} 
\hspace{-5mm} 
\includegraphics[width=7.001cm]{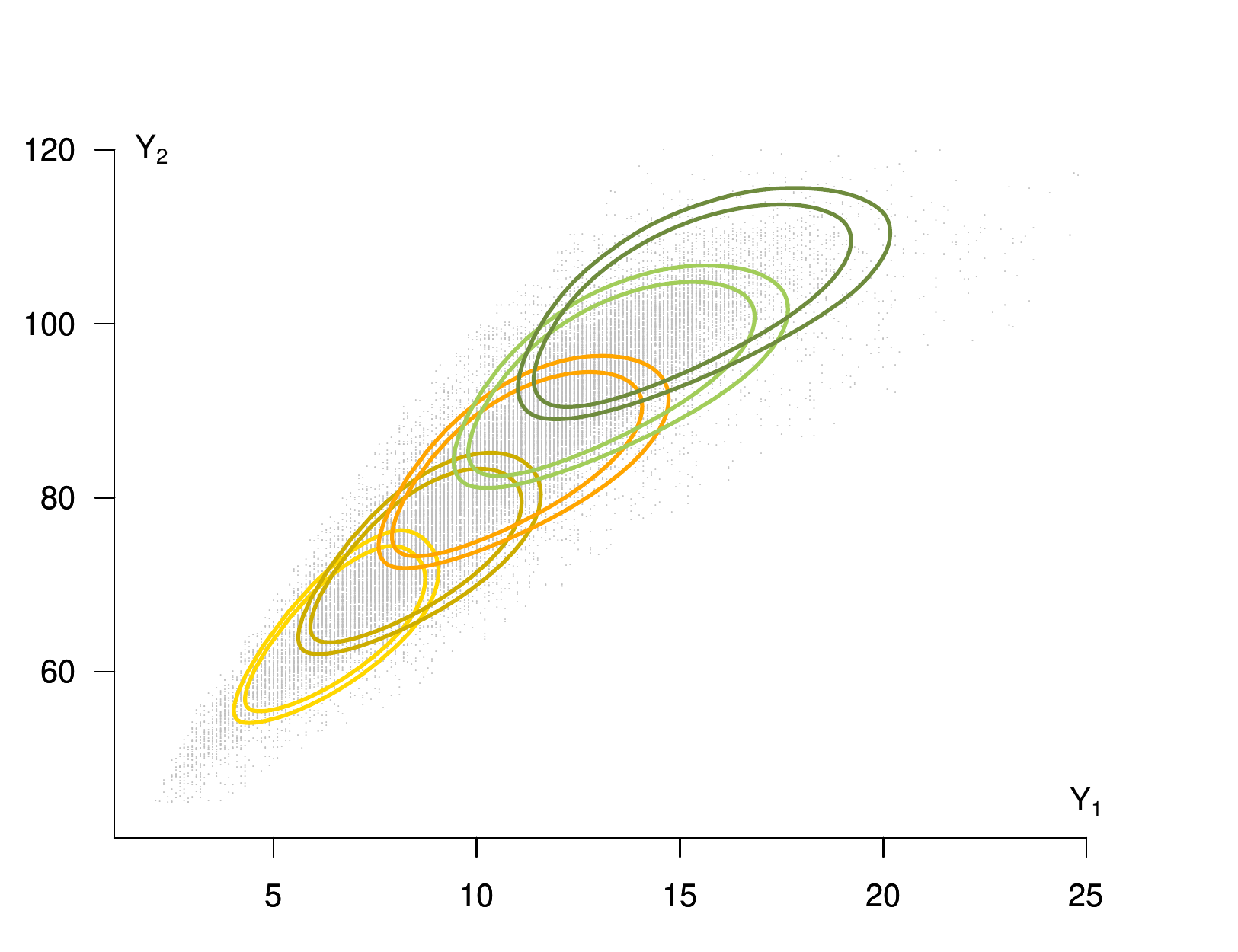}
&
\hspace{-7mm}
\includegraphics[width=7.001cm]{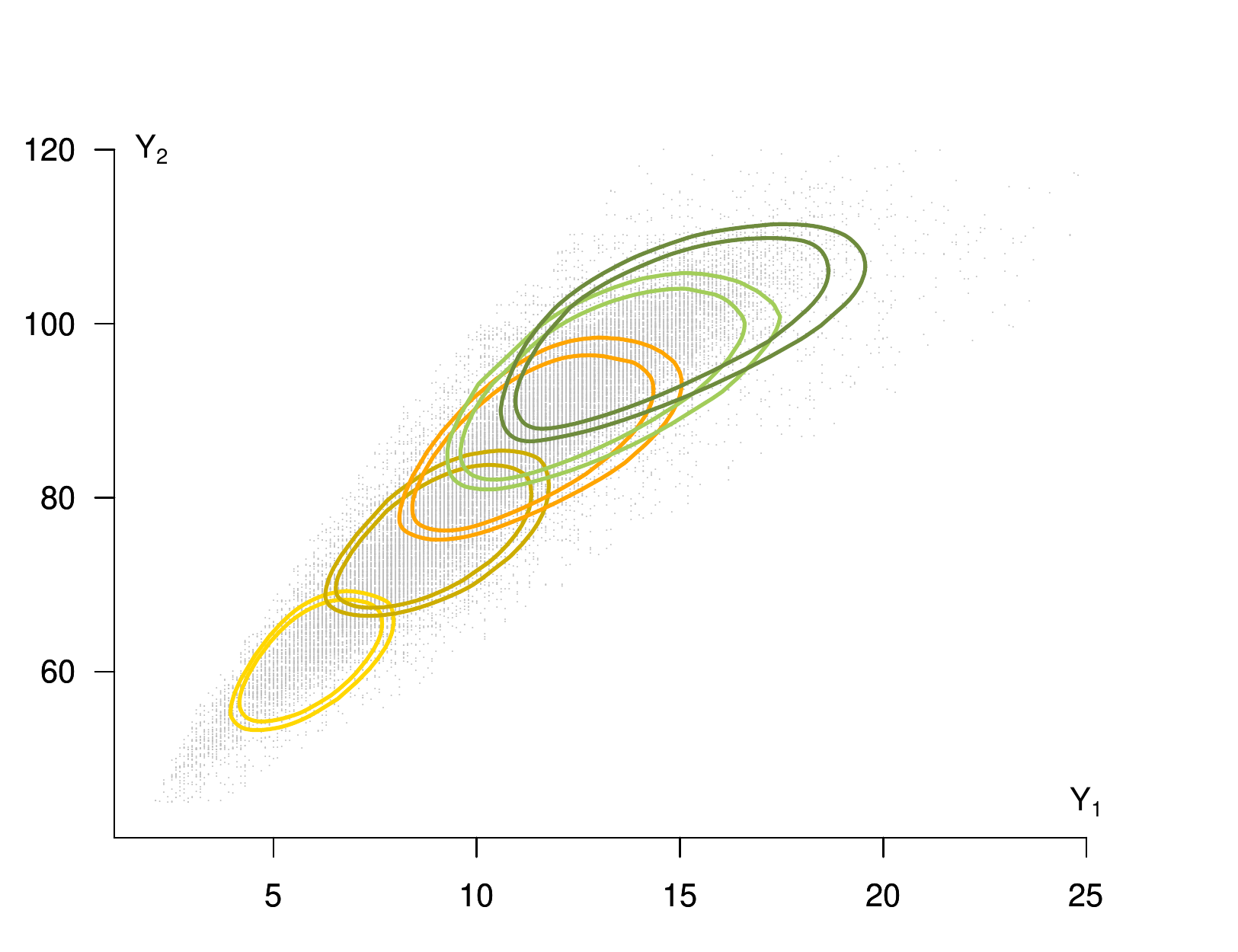}
\\
\hspace{-5mm}
\includegraphics[width=7.001cm]{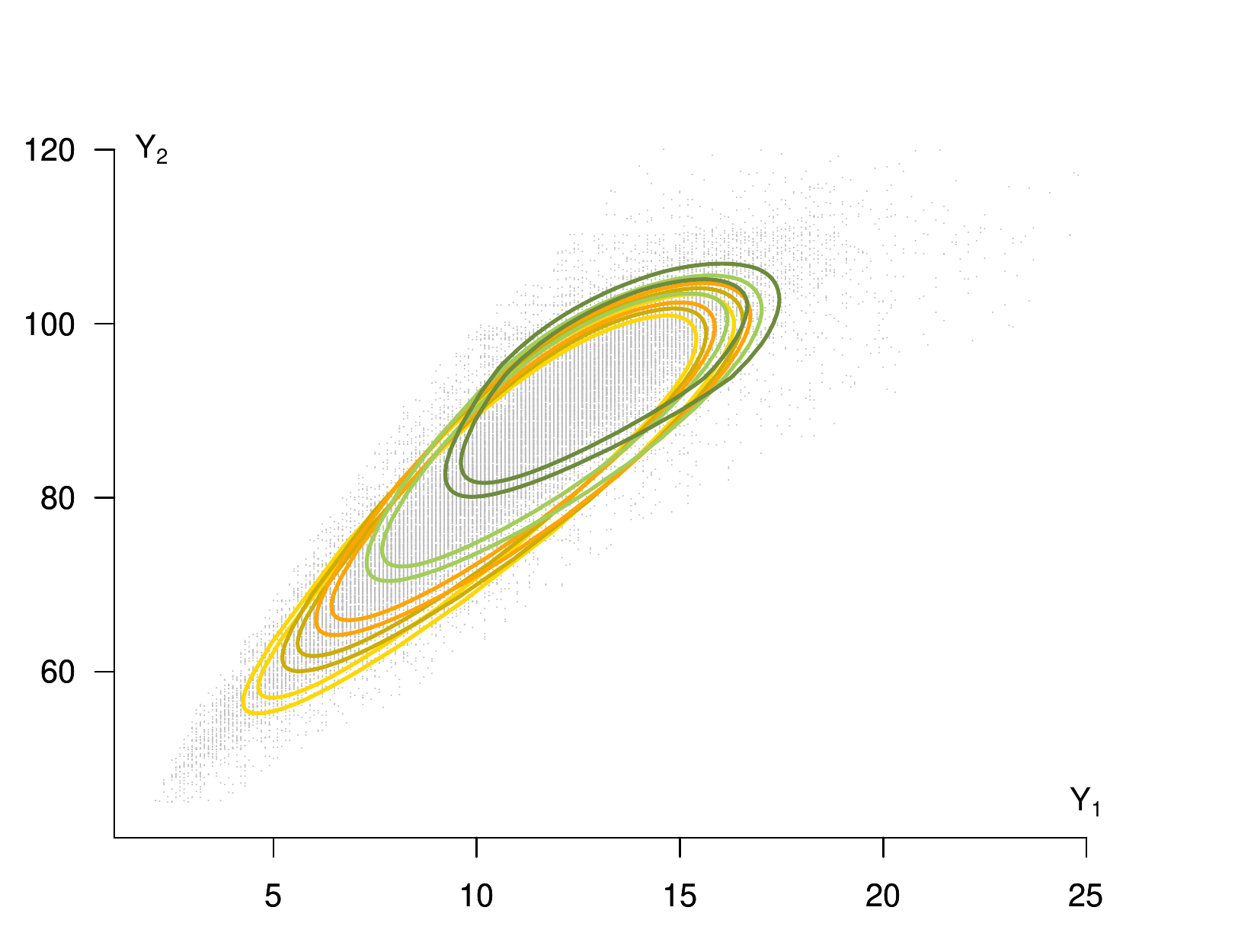}
&
\hspace{-7mm}
\includegraphics[width=7.001cm]{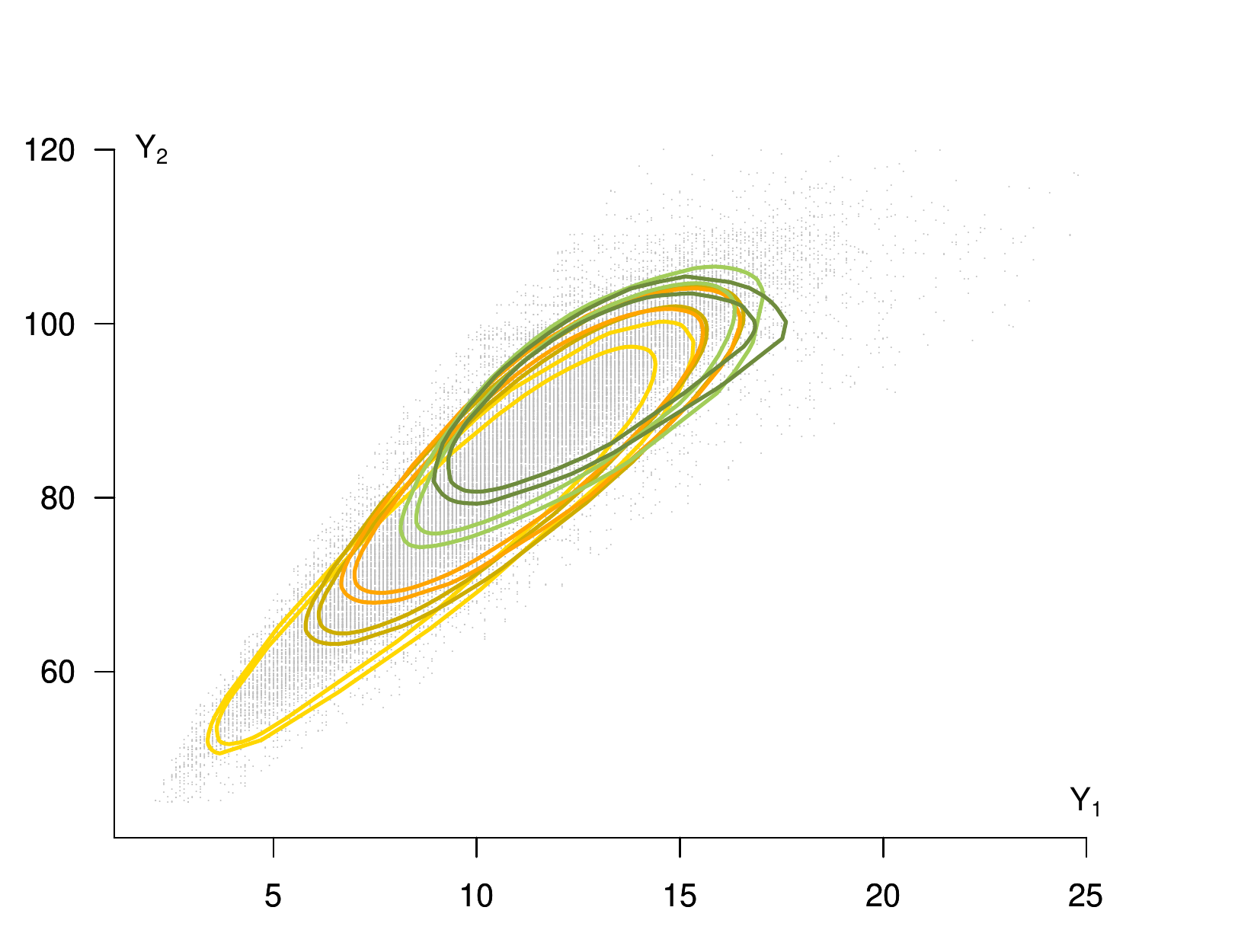}
\\
\hspace{-5mm}
\includegraphics[width=7.001cm]{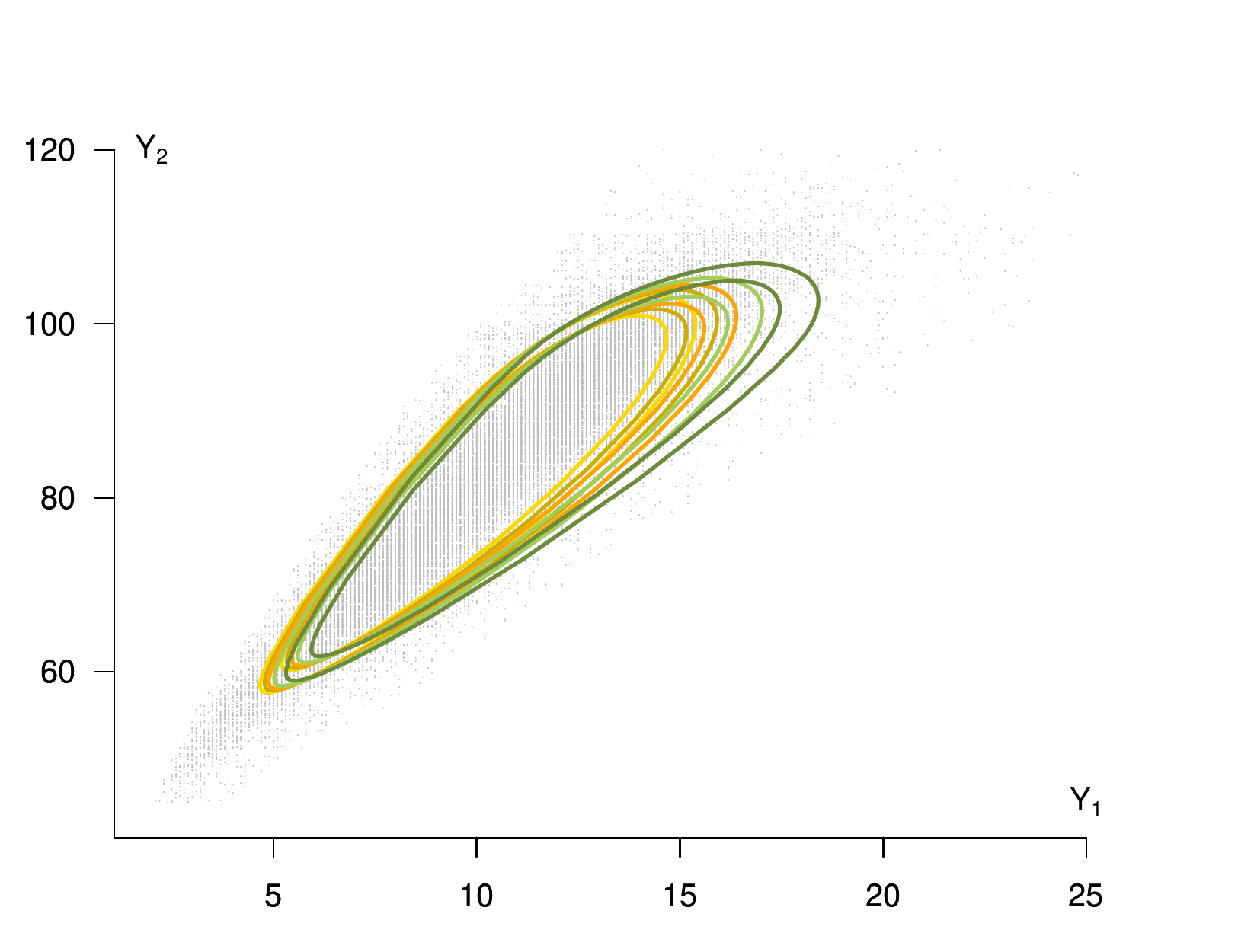}
&
\hspace{-7mm}
\includegraphics[width=7.001cm]{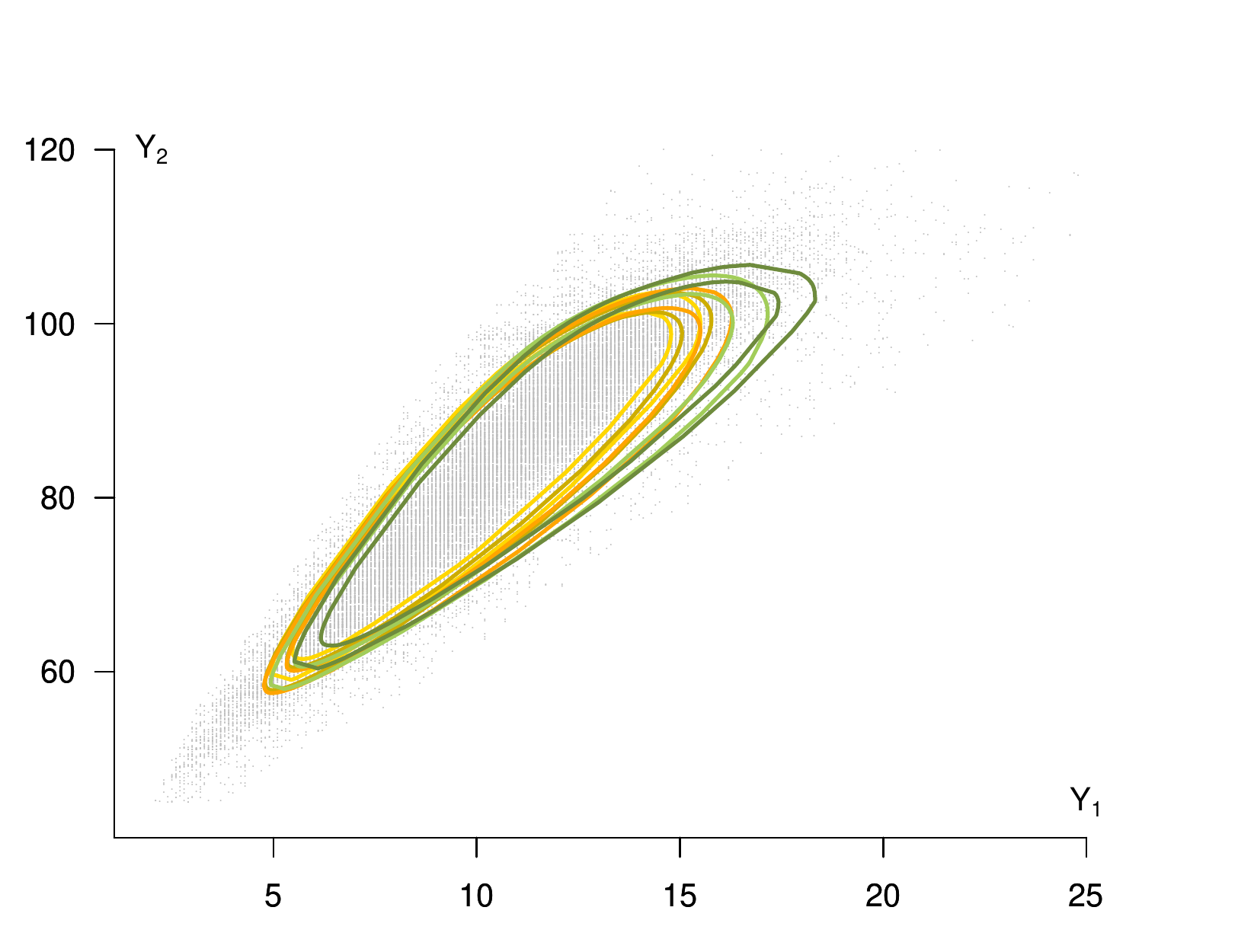}
\end{array}$
\end{center}
     \vspace*{1mm}
\caption{(Top:) conditional expectile contours~$\partial \hat{R}^{\rho(n)}_{\alpha,x}$ for the joint children's {\it weights} $(Y_1)$ and {\it heights} $(Y_2)$ conditional on their {\it age}, for $\alpha\in\{.005,.01\}$ and for values of~$x$ that are the $5\%$ (yellow), $25\%$ (brown), $50\%$ (orange),
$75\%$ (light green) and $95\%$ (dark green) empirical quantiles of the covariate, obtained from a linear regression method (left) or from a nonparametric regression method (right). (Middle:) results when the covariate is the duration of {\it breastfeeding}. (Bottom:) results when the covariate is the mother's {\it BMI}; see Section~\ref{sec72} for details. Bivariate responses~$(Y_{i1},Y_{i2})$, $i=1,\ldots,n$, are shown in gray.
}
\label{Fig6}
\end{figure}

For each of the three covariates, Figure~\ref{Fig6} plots both linear and nonparametric conditional expectile contours associated with the extreme levels~$\alpha\in\{.005, .01\}$ and several covariate values~$x$. Like for the simulated example, these contours are smooth and nested. For each covariate, we could comment on trend and heteroscedasticity; for instance, age provides, as expected, a monotone increase in both trend and variability. We could also compare the linear and nonparametric fits to identify nonlinear effects; in particular, this reveals that the effect of age is linear only after the first age quartile, that is, after 16 months. Regarding the specific impact of covariates on the joint distribution of the responses, it is seen, e.g., that the first principal direction of the bivariate response distribution becomes more and more horizontal as children get older. There does not seem to be a strong trend nor heteroscedasticity for the BMI covariate. But it is seen that mothers with a BMI above median may lead to overweigthed tall children, but not to overweighted short ones; at such BMI levels, both green expectile contours indeed show an asymmetry to obesity in the upper-right direction, but not in the lower-left one (this is arguably not related to malnutrition, but it is still pointing to some risk factor). Similar comments could be given for the breastfeeding covariate. Finally, while we do not show here the quantile results (they are available from the authors on simple request), we mention that the corresponding quantile contours are more elliptical than the expectile ones, hence, e.g., do not reveal subtle risks such as the one related to overweight of tall children from mothers with a large BMI.

%%%%%%%%%%%%%%%%%%%%

%%%%%%%%%%

\section{Final comments and perspectives for future research}
\label{secfinalcom}

As shown in Section~\ref{secregres}, we introduced in this paper a methodology that allows practitioners to easily conduct multiple-output M-quantile regression. The expectile case is of particular interest due to its relation with classical mean regression and to its advantages over quantile regression: multiple-output expectile regression provides smoother and more flexible contours, that are easier to interpret and do not show crossings. This expectile regression method has potential applications in, e.g., financial risk management, econometrics, or any field where extreme values and tail information are relevant. Our construction led to introduce new statistical depths, namely the halfspace M-depths, and new multivariate, affine-equivariant, \mbox{M-quantiles}. Recently,  \cite{Caretal2016}
%,Caretal2017} 
and \cite{Cheetal2017} defined multivariate quantiles that are equivariant under large groups of diffeomorphic transformations, collecting gradients of convex functions, which extends to the multivariate setup the equivariance of quantiles under monotone increasing transformations of the real line. It should be noted that other univariate \mbox{M-quantiles} are not equivariant under monotone increasing transformations, so that the fact that our multivariate \mbox{M-quantiles} are equivariant under affine transformations but not under such diffeomorphic transformations is a limitation for standard quantiles only, hence not for expectiles nor for other M-quantiles. 

Perspectives for future research are rich and diverse. On the inferential side, it would be natural to study how the properties of the M-location functional~$\zb_*(P)$ (see Theorem~\ref{TheorMmedianConsistency}) depend on the loss function~$\rho$. Also, the corresponding estimators~$\zb_*(P_n)$ include the sample average (for $\rho(t)=t^2$) and the Tukey median (for $\rho(t)=|t|$) as particular cases, which leads to investigating whether or not a nice trade-off between efficiency and robustness can be achieved by choosing a suitable Huber loss function~$\rho_c(t)$ in~(\ref{rhoc}) or, in case affine equivariance of~$\zb_*(P_n)$ is a requirement, a suitable power loss function~$\rho(t)=|t|^r$. As for hypothesis testing, we suggested in Section~\ref{secexpectfurtherprop} 
%a nonparametric 
an affine-invariant location test based on expectile depth. It would be desirable to investigate the asymptotic null and non-null properties of this test. Further questions of interest are: does the monotonicity property in Theorem~\ref{Theordepthuniqueu} hold in the sample case (with plain differentiability replaced by left- and right-differentiability)? Extensive numerical exercises lead us to conjecture that the answer is positive. We closed Section~\ref{secrisks} by suggesting alternative set-valued risk-measures that do not require choosing a specific direction~$u$ yet have the topology of natural risk measures (unlike centrality regions); are these alternative risk measures coherent in the sense of Definition~\ref{defincoher}? Finally, is it possible to define an expectile---or, more generally, an M-quantile---concept of scatter depth that would extend the concept of halfspace (Tukey) scatter depth considered in \cite{PaiBev18} and \cite{Chenetal16}? The same question holds for the shape depth concept recently introduced in \cite{PVB19}, that also relies on the halfspace Tukey depth.

%%%%%%%%%%

\section*{Acknowledgements}

%We would like to thank the Editor Edward \mbox{I.} George, the Associate Editor and two anonymous referees for their insightful comments and suggestions, that led to a substantial improvement of a previous version of this work. 
 We are very grateful to \mbox{Prof.} Rainer Dyckerhoff for providing the data used in Section~\ref{secrisks} and to \mbox{Prof.} Philip Kokic for sharing his Matlab code computing the M-quantiles from \cite{Breetal2001} and \cite{Koketal2002}. We also want to thank \mbox{Prof.} Germain Van Bever for his comments on an earlier version of this work. 
%\vspace{-2mm}

%%%%%%%%%%

%\newpage
\bibliographystyle{imsart-nameyear.bst} 
\bibliography{Paper.bib}

%%%%%%%%%%%%%%%%%%%%%%%%%%%%%%%%%%%%%%%%%%%%%%%%%%

%%%%%%%%%%%%%%%%%%%%  
%\clearpage
%%%%%%%%%%%%%%%%%%%%  

\appendix

\section{Competing M-quantile concepts}
\label{secAppA}

In this first appendix, we define some of the main concepts of multivariate M-quantiles and multivariate expectiles available in the literature, with a particular emphasis on the concepts we used above for comparison with the proposed multivariate M-quantiles and multivariate expectiles. 

Before proceeding, it is needed to introduce an alternative parametrization of the univariate M-quantiles~$\theta^\rho_\alpha=\theta^\rho_\alpha(P)$ from Section~\ref{secMquantile}; the dependence on~$P$ will play no role in this section, hence will be dropped in the notation. 
\vspace{-.6mm}
This alternative parametrization is~$\theta_{\tau,u}^{\rho}:=\theta_{(1-\tau u)/2}^{\rho}$ and 
%$\tau u=1-2\alpha$ 
 indexes univariate M-quantiles by an order~$\tau\in[0,1)$ and a direction~$u\in\{-1,1\}$, or equivalently by an order~$\tau u$ that belongs to the open unit ``ball"~$(-1,1)$ of~$\R$. In this directional parametrization, the most central M-quantile corresponds to~$\tau=0$ and the most extreme ones are obtained as~$\tau\to 1$. In the $d$-dimensional case ($d\geq 2$), where there are no left nor right, it is natural to similarly index M-quantiles by an order~$\tau\in[0,1)$ and a direction~$\ub\in\mathcal{S}^{d-1}$, or equivalently by a vectorial order~$\tau \ub$ belonging to the open unit ball~$\{\zb\in\R^d:\|\zb\|<1\}$ of~$\R^d$, with the same idea that~$\tau=0$ will yield the most central M-quantile and that~$\tau \ub$ with~$\tau\to 1$ will provide extreme M-quantiles in direction~$\ub$. It is then standard (see the references below) to consider contours generated by M-quantiles of a fixed order~$\tau$. As in the body of the paper, these contours are the boundaries of ``centrality regions" that provide a center-outward ordering of points in~$\R^d$. There are alternative directional parametrizations for M-quantiles; typically, these involve an order~$\alpha\in(0,1)$ and a direction~$\ub\in\mathcal{S}^{d-1}$, and are such that the M-quantile of order~$\alpha$ in direction~$\ub$ is equal to the \mbox{M-quantile} of order~$1-\alpha$ in direction~$-\ub$. For such a parametrization, central M-quantiles are associated with~$\alpha=1/2$, whereas extreme ones are  obtained as~$\alpha\to 0$ and $\alpha\to 1$ (this is the parametrization of multivariate M-quantiles that was used in the paper). In the rest of this section, we discriminate between these two types of parametrization by using the notation~$\tau$ and~$\alpha$ in a consistent way. 

These general considerations allow us to review some of the main concepts of multivariate M-quantiles and expectiles. The first concept of multivariate M-quantiles can be found in \cite{BreCha1988}.
%, where the univariate M-quantiles themselves were introduced. 
There, for any~$\alpha\in(0,\frac{1}{2}]$, the order-$\alpha$ M-quantile of~$P$ in direction~$\ub$ is defined as the ``geometric" quantity  
$$
\thetab_{\alpha,\ub}^{\rho,{\rm geom}}
:=
\arg\min_{\thetab\in\mathbb{R}^d} 
{\rm E}\bigg[ 
\bigg\{
1
-(1-2\alpha)
\frac{\ub'(\Zb- \thetab)}{\|\Zb- \thetab\|}
\bigg\}
\rho(\|\Zb- \thetab\|) 
\bigg]
$$
(throughout this appendix, $Z$ is a random $d$-vector with distribution~$P$), whereas the definition for~$\alpha\in(\frac{1}{2},1)$ results from the identity~$\thetab_{\alpha,\ub}^{\rho,{\rm geom}}=\thetab_{1-\alpha,-\ub}^{\rho,{\rm geom}}$. In the univariate case, $\thetab_{\alpha,1}^{\rho,{\rm geom}}$ reduces to the M-quantile~$\theta_{\alpha}^{\rho}$ in~(\ref{Mquantfirstdef}). 
%$$
%{\rm E}\bigg[ 
%\bigg(
%1
%-(1-2\alpha)
%\frac{\ub'(Z- \theta)}{\|\Zb- \thetab\|}
%\bigg)
%\rho(\|\Zb- \thetab\|) 
%\bigg]
%$$
%$$
%=
%{\rm E}\bigg[ 
%\bigg(
%1
%-(1-2\alpha)
%S(Z- \theta)
%\bigg)
%\rho(|Z- \theta|) 
%\bigg]
%$$
%$$
%=
%{\rm E}\bigg[ 
%\bigg(
%2(1-\alpha)
%\mathbb{I}[Z-\theta<0]
%+
%\mathbb{I}[Z-\theta=0]
%+
%2\alpha
%\mathbb{I}[Z-\theta>0]
%\bigg)
%\rho(|Z- \theta|) 
%\bigg]
%$$
%$$
%=
%2
%{\rm E}\bigg[ 
%h_\alpha(Z-\theta)
%\rho(|Z- \theta|) 
%\bigg]
%$$
%$$
%=
%2
%{\rm E}\bigg[ 
%h_\alpha(Z-\theta)
%\rho(Z- \theta) 
%\bigg]
%$$
%(recall $\rho$ is symmetric)
The term ``geometric" above is justified by the fact that, for~$\rho(t)=|t|$ and~$\rho(t)=t^2$, these M-quantiles reduce to the geometric quantiles
\begin{equation}
	\label{quantgeom}
\qb_{\alpha,\ub}^{\rm geom}
:=
\arg\min_{\thetab\in\mathbb{R}^d} 
{\rm E}\big[ 
\|\Zb- \thetab\|
-
(1-2\alpha) 
\ub'(\Zb- \thetab)
\big]
\end{equation}
and geometric expectiles
\begin{equation}
	\label{expectgeom}
\eb_{\alpha,\ub}^{\rm geom}
:=
\arg\min_{\thetab\in\mathbb{R}^d} 
{\rm E}\big[ 
\|\Zb- \thetab\|
\{\|\Zb- \thetab\|
-
(1-2\alpha) 
 \ub'(\Zb- \thetab) \}
\big]
\end{equation}
from \cite{Cha1996} and \cite{Heretal2018}, respectively; these papers, that actually rather rely on the~$(\tau,\ub)$-directional parame\-trization, would refer to~(\ref{quantgeom})/(\ref{expectgeom}) as quantiles/expectiles of order~$\tau=1-2\alpha$ in direction~$-\ub$.  %As recently showed theoretically in \cite{GirStu2015,GirStu2017}, geometric quantiles exhibit undesirable properties. In particular, (a) the extreme quantile contours obtained as~$\alpha\to 0$ may extend far outside the support of the distribution. Also, (b) such extreme quantile contours exhibit a structure that is incompatible with the principal component structure of the underlying distribution: more precisely, they will be furthest (resp., closest) to the center of the distribution in the last (resp., first) principal direction, which is orthogonal to what one would expect from quantile contours. As we show empirically in Figure~\ref{Fig1}, these pathological features unfortunately extend to geometric expectiles (the empirical version of~(\ref{expectgeom}) is simply obtained by replacing the expectation with a sample average over the observations~$\Zb_1,\ldots,\Zb_n$ at hand). 

As we mentioned in the paper, geometric M-quantiles may extend far outside the support of the distribution as~$\alpha\to 0$. To improve on this,  \cite{Breetal2001} and \cite{Koketal2002} introduced alternative concepts of multivariate M-quantiles, actually only for Huber loss functions. To define these quantiles, we need to introduce the following notation: let~$S(t):=\mathbb{I}[t>0]-\mathbb{I}[t<0]$ \label{pagedefsign} be the sign function and~$
%\tb\mapsto
\psi_c(\tb) 
:=
(\tb/c) \mathbb{I}[\|\tb\|<c] 
+  
(\tb/\|\tb\|) \mathbb{I}[\|\tb\|\geq c]
%= (\max(\|x\|,c))^{-1} x,
$
be a $d$-variate extension of Huber's \mbox{$\psi$-function}. Write also~$h_\alpha(t):=
(1-\alpha) \mathbb{I}[t<0]
+ (1/2) \mathbb{I}[t=0]
+\alpha \mathbb{I}[t>0]
$.
Then, for~$c,\delta>0$, 
%the \cite{Breetal2001} order-$\alpha$ M-quantile~$\thetab^{c}_{\alpha,\ub}$ of~$\Zb$ in direction~$\ub$ is defined as the value of~$\theta(\in\R^d)$ that makes the quantity 
%\begin{eqnarray*}
%\label{Kokicpre}	
%{\rm E}\Bigg[
%\bigg\{
%%\big({\textstyle{\frac{1}{2}}}-\alpha\big)
%(1-2\alpha)
%S(\ub'(\Zb-\thetab)) 
%\bigg( 1-  
%\frac{|\ub'(\Zb-\thetab)|}{\|\Zb-\thetab\|} \bigg)
%+
%2\tilde{h}_\alpha(\ub'(\Zb-\thetab))
%\bigg\} 
%\psi_c(\Zb-\thetab)
%\Bigg]
%\end{eqnarray*}
%equal to zero. 
%%. 
%With the same notation, 
the \cite{Koketal2002} order-$\alpha$ \mbox{M-quantile}~$\thetab^{\delta,c}_{\alpha,\ub}$ of~$P$ in direction~$\ub$ is the solution~$\theta(\in\R^d)$ of 
\begin{equation}
\label{Kokic}	
{\rm E}\Bigg[
\bigg\{
%\big({\textstyle{\frac{1}{2}}}-\alpha\big)
(1-2\alpha)
S(\ub'(\Zb-\thetab)) 
\bigg( 1-  
\frac{|\ub'(\Zb-\thetab)|}{\|\Zb-\thetab\|} \bigg)^\delta
+
2h_\alpha(\ub'(\Zb-\thetab))
\bigg\} 
\psi_c(\Zb-\thetab)
\Bigg]=0
.
\end{equation}
%\begin{equation*}
%\eta_{\delta}(p_i) := 
% \left\{
%\begin{array}{lll}
%\, \, \,  \, (1-\cos p_i)^{\delta} (1-2\alpha) + 2\alpha          & \quad \mbox{if} &  p_i\in (-\frac{\pi}{2},\frac{\pi}{2})  \\
%- (1+\cos p_i)^{\delta} (1-2\alpha) + 2(1-\alpha)          & \quad \mbox{if} &   p_i\in [-\pi,-\frac{\pi}{2}] \cup [\frac{\pi}{2},\pi],
%\end{array}
%\right 
%.
%\end{equation*}
In the univariate case, $\psi_c$ is the derivative of the Huber loss function~$\rho_c$ in~(\ref{rhoc}), which allows showing that, for~$u=1$ and any~$\delta>0$,~(\ref{Kokic}) reduces to the first-order condition~\mbox{$G^{\rho_c}(\theta)=\alpha$};
\vspace{-.4mm}
 see the last paragraph of Section~\ref{secMquantile}.
%
%In the univariate case, we obtain (for~$u=1$):
%$$
%{\rm E}\Big[
%h_\alpha(Z-\theta)
%\psi_c(Z-\theta)
%\Big]=0
%$$
%$$
%(1-\alpha)
%{\rm E}\Big[
%\mathbb{I}[Z-\theta<0]
%\psi_c(Z-\theta)
%\Big]
%+
%\frac{1}{2}
%{\rm E}\Big[
%\mathbb{I}[Z-\theta=0]
%\psi_c(Z-\theta)
%\Big]
%+
%\alpha
%{\rm E}\Big[
%\mathbb{I}[Z-\theta>0]
%\psi_c(Z-\theta)
%\Big]
%=0
%$$
%$$
%(1-\alpha)
%{\rm E}\Big[
%\mathbb{I}[Z-\theta<0]
%\psi_c(Z-\theta)
%\Big]
%+
%\alpha
%{\rm E}\Big[
%\mathbb{I}[Z-\theta>0]
%\psi_c(Z-\theta)
%\Big]
%=0
%$$
%$$
%{\rm E}\Big[
%\mathbb{I}[Z-\theta<0]
%\psi_c(Z-\theta)
%\Big]
%+
%\alpha
%{\rm E}\Big[
%(\mathbb{I}[Z-\theta>0]-\mathbb{I}[Z-\theta<0])
%\psi_c(Z-\theta)
%\Big]
%=0
%$$
%$$
%{\rm E}\Big[
%\mathbb{I}[Z-\theta<0]
%(-\psi_c(Z-\theta))
%\Big]
%-
%\alpha
%{\rm E}\Big[
%(\mathbb{I}[Z-\theta>0]-\mathbb{I}[Z-\theta<0])
%\psi_c(Z-\theta)
%\Big] 
%=0
%$$
%$$
%{\rm E}\Big[
%|\mathbb{I}[Z-\theta<0] \psi_c(Z-\theta)|
%\Big]
%-
%\alpha
%{\rm E}\Big[
%|\psi_c(Z-\theta)|
%\Big]
%=0
%$$
%$$
%\alpha
%=
%\frac{
%{\rm E}\Big[
%|\mathbb{I}[Z-\theta<0] \psi_c(Z-\theta)|
%\Big]
%}
%{
%{\rm E}\Big[
%|\psi_c(Z-\theta)|
%\Big]
%}
%$$
%
Consequently, $\theta^{\delta,c}_{\alpha,1}$, for any~$\delta>0$, reduces
\vspace{-.8mm}
 to the univariate \mbox{M-quantile}~$\theta_\alpha^{\rho_c}$ from~(\ref{Mquantdef}), so that~$\thetab^{\delta,c}_{\alpha,\ub}$ may indeed be considered as a multivariate \mbox{M-quantile}. The multivariate \mbox{M-quantiles} from \cite{Breetal2001} correspond to the particular case obtained for~$\delta=1$. For both the M-quantiles from \cite{Breetal2001} and \cite{Koketal2002}, the parameter~$c>0$, as it was the case for univariate M-quantiles, allows going continuously from quantiles to expectiles, that are obtained as~$c\to 0$ and~$c\to \infty$, respectively. 
  %In the sequel, these multivariate M-quantiles will be used as the main competitors of the multivariate M-quantiles we will propose. 

%Note that the limiting case obtained as~$\delta\to\infty$ solves
%$$
%{\rm E}\Big[ 
%h_\alpha(\ub'(Z-\theta))
%\psi_c(Z-\theta)
%\Big]=0
%,
%$$ 
%
%\cite{CouDiB2013,CouDiB2014}

As explained in the main body of the paper, an important drawback of the aforementioned multivariate M-quantiles (and of other multivariate M-quantiles, such as those from \citealp{Kol1997}) is their weak equivariance properties. More precisely, these M-quantiles are equivariant under orthogonal transformations, but they fail to be equivariant under general affine transformations. Actually, other recent proposals enjoy even weaker equivariance properties; for instance, the multivariate expectiles from \cite{Mauetal2018,Mauetal2017} are not equivariant under orthogonal transformations.    
%\cite{Sai2016,Mauetal2017,Mauetal2018}

%%%%%%%%%%%%%%%

\section{Proofs}
\label{secAppB}

This second appendix presents the proofs of all results stated in the paper. 
%
%\subsection{Proof of Theorem~\ref{Theorlemjones}} 
%
The proofs 
%of Theorem~\ref{Theorlemjones} 
will make use of the following comments related to any loss function~$\rho\in\mathcal{C}$. Since~$\rho$ is convex on~$\R$, it is continuous and it admits a left-derivative function~$\psi_-$ and right-derivative function~$\psi_+$ that are both monotone non-decreasing and satisfy~$\psi_-(z)\leq \psi_+(z)$ at any~$z$; see Theorem~1.3.3 in \cite{NicPer2006}. The functions~$\psi_-$ and~$\psi_+$ are left- and right-continuous, respectively, and both have an at most countable number of discontinuity points; see pages~6-7 in \cite{RobVar1973}. Since~$\rho$ has a global minimum at~$0$, we have~$\psi_\pm(z)\leq 0$ for~$z<0$, $\psi_-(0)\leq 0\leq \psi_+(0)$, and $\psi_\pm(z)\geq 0$ for~$z>0$. 

We will need the following mean-value theorem for one-sided derivatives; see, e.g., \cite{LeiCoh2007}. 

\begin{Lem}%[\cite{LeiCoh2007}]
\label{LemOneSidedMeanValueTheor}
Let~$f:[a,b]\to\R$ be a continuous function.  
(i) Assume that $f$ is left-differentiable on~$(a,b)$, with left-derivative~$f'_-$. Then,
$$
f'_-(c_1)
\leq
\frac{f(b)-f(a)}{b-a}
\leq
f'_-(c_2)
$$  
for some~$c_1,c_2\in(a,b)$.
(ii) Assume that $f$ is right-differentiable on~$(a,b)$, with right-derivative~$f'_+$. Then,
$$
f'_+(c_1)
\leq
\frac{f(b)-f(a)}{b-a}
\leq
f'_+(c_2)
$$  
for some~$c_1,c_2\in(a,b)$.
\end{Lem}

%We will make repeated use of this result. For now, 
This result implies in particular that~$\psi_-(z)$ may be zero at~$z=0$ only. By contradiction, assume  that there exists~$z_0\neq 0$ with~$\psi_-(z_0)=0$. Monotonicity of~$\psi_-$ then implies that~$\psi_-(z)=0$ for any~$z$ strictly between~$0$ and~$z_0$, which, from Lemma~\ref{LemOneSidedMeanValueTheor}, entails that~$\rho(z_0)=0$, a contradiction since~$\rho(z)=0$ at~$z=0$ only. Of course, one shows similarly that~$\psi_+(z)$ may be zero at~$z=0$ only.

\begin{Lem}
\label{LemDCTpsi}
\
\hspace{-5mm} 
\mbox{Fix~$\rho\in\mathcal{C}$ and~$P\in\mathcal{P}^{\rho}$. Then, for any~$\theta_0,a,b\in\R$ with~$a<b$,}
\begin{eqnarray*}	
\lefteqn{
\lim_{\theta\stackrel{>}{\to} \theta_0} 
\int_{-\infty}^{\infty}
\psi_-(z-\theta) 
\mathbb{I}[\theta+a< z < \theta+b]
\,dP(z)
}
\\[2mm]
& & 
\hspace{13mm} 
=
\int_{-\infty}^{\infty}
\psi_-(z-\theta_0) 
\mathbb{I}[\theta_0+a < z\leq \theta_0+b]
\,dP(z)
\end{eqnarray*}
and
\begin{eqnarray*}	
\lefteqn{
\lim_{\theta\stackrel{<}{\to} \theta_0} 
\int_{-\infty}^{\infty}
\psi_-(z-\theta) 
\mathbb{I}[\theta+a< z< \theta+b]
\,dP(z)
}
\\[2mm]
& & 
\hspace{13mm} 
=
\int_{-\infty}^{\infty}
\psi_-(z-\theta_0+0) 
\mathbb{I}[\theta_0+a\leq z< \theta_0+b]
\,dP(z)
,
\end{eqnarray*}
 where~$\psi_-(z-\theta_0+0)$ denotes the limit of~$\psi_-(t)$ as~$t$ converges to~$z-\theta_0$ from above.     
\end{Lem}

\noindent {\sc Proof of Lemma~\ref{LemDCTpsi}.}
For any~$h\in(0,1]$, let~$r(h,z):=\psi_-(z-\theta_0-h) 
\mathbb{I}[\theta_0+h+a< z < \theta_0+h+b]
-
\psi_-(z-\theta_0) 
\mathbb{I}[\theta_0+a < z\leq \theta_0+b]
$. Note that, for any~$z$,~$r(h,z)\to 0$ as~$h$ converges to zero from above. Since
$$
|r(h,z)|
\leq
2\max(|\psi_-(a)|,|\psi_-(b)|)
%\leq
%2|\psi_-(z-a)|+2|\psi_-(z-b)|
$$
for any~$h\in(0,1]$, the Lebesgue Dominated Convergence Theorem thus implies that
$
\int_{-\infty}^{\infty}
|r(h,z)|
%
%\linebreak
%
\,dP(z)
\to 0$ as~$h$ converges to zero from above, which establishes the first statement of the lemma. The second statement follows similarly by using the fact that~$\ell(h,z):=\psi_-(z-\theta_0-h) 
\mathbb{I}[\theta_0+h+a< z < \theta_0+h+b]
-
\psi_-(z-\theta_0+0) 
\mathbb{I}[\theta_0+a \leq z< \theta_0+b]
\to 0$ as~$h$ converges to zero from below.  
\cqfd
\vspace{3mm}

We can now prove Theorem~\ref{Theorlemjones}. 
\vspace{2mm}
 
{\sc Proof of Theorem~\ref{Theorlemjones}}.
(i) 
%Consider the function mapping~$\theta$ to 
%$$
%O(\theta)
%:=
%{\rm E}\big[ 
%\rho_\alpha(Z- \theta) -\rho_\alpha(Z) 
%\big]
%.
%$$
We first show that~$O(\theta)$ is well-defined for any~$\theta\in\R$. To do so, note that~$\theta \mapsto \rho_\alpha(z-\theta)$ is continuous and right-differentiable with right-derivative~$((1-\alpha)\mathbb{I}[z\leq \theta]+\alpha \mathbb{I}[z>\theta]) \psi_-(z-\theta)$, so that Lemma~\ref{LemOneSidedMeanValueTheor} shows that there exists~$\lambda\in(0,1)$ such that
\begin{eqnarray*}
\big|
\rho_\alpha(z-\theta) -\rho_\alpha(z) 
\big|
&\leq &
|\theta|
\max(\alpha,1-\alpha)
|\psi_-(z-\lambda \theta)|
\\[1mm]
&\leq &
|\theta|
\max(\alpha,1-\alpha)
|\psi_-(-|z|-|\theta|)|
,
\end{eqnarray*}
where we used the fact that, for any~$a>0$, the upper bound of~$|\psi_-(z)|$ on~$[-a,a]$ is~$|\psi_-(-a)|$. Since~$P\in\mathcal{P}^\rho$, this entails that~$z\mapsto \big|
\rho_\alpha(z-\theta) -\rho_\alpha(z)\big|$ is~$P$-integrable for any~$\theta$, so that the mapping~$\theta\mapsto O(\theta)$ is indeed well-defined for any~$\theta$. 

We turn to right-differentiability of~$\theta\mapsto O(\theta)$ and show that
\begin{eqnarray*}
	\lefteqn{
O_+^{\rho\prime}(\theta)
=
(1-\alpha)
\int_{-\infty}^{\infty}
|\psi_-(z-\theta)| \mathbb{I}[z\leq \theta]
\,dP(z)
}
\\[2mm]
& & 
\hspace{40mm} 
-
\alpha
\int_{-\infty}^{\infty}
|\psi_-(z-\theta)|\mathbb{I}[z>\theta]
\,dP(z)
.
\end{eqnarray*}
To do so, fix~$\theta_0\in\R$, $h\in(0,1]$, and write
\begin{eqnarray}
\lefteqn{
\frac{O(\theta_0+h)-O(\theta_0)}{h}
-
O_+^{\rho\prime}(\theta_0)
}	
\nonumber
\\[2mm]
& & 
\hspace{9mm} 
=
\int_{-\infty}^{\infty}
\bigg\{
\frac{\rho_\alpha(z-\theta_0-h)-\rho_\alpha(z-\theta_0)}{h} - O_+^{\rho\prime}(\theta_0)
\bigg\}
\,dP(z)
\nonumber
\\[2mm]
& & 
\hspace{9mm} 
=
(1-\alpha )
\int_{-\infty}^{\infty}
L_{\theta_0}(h,z)
\,dP(z)
+
\alpha 
\int_{-\infty}^{\infty}
R_{\theta_0}(h,z)
\,dP(z)
,
\label{rightdiffO}
\end{eqnarray}
where we let
\begin{eqnarray*}
	\lefteqn{
\hspace{-33mm} 
L_{\theta_0}(h,z)
:=
\frac{\rho(z- \theta_0-h) \mathbb{I}[z<\theta_0+h]-\rho(z- \theta_0) \mathbb{I}[z<\theta_0]}{h}
}
\\[2mm]
& & 
\hspace{-5mm} 
-
|\psi_-(z-\theta_0)|\mathbb{I}[z\leq \theta_0]
\end{eqnarray*}
and
\begin{eqnarray*}
	\lefteqn{
\hspace{-33mm} 
R_{\theta_0}(h,z)
:=
\frac{\rho(z- \theta_0-h) \mathbb{I}[z\geq \theta_0+h]-\rho(z- \theta_0) \mathbb{I}[z\geq \theta_0]}{h}
}
\\[2mm]
& & 
\hspace{-5mm} 
+
|\psi_-(z-\theta_0)|\mathbb{I}[z> \theta_0]
.
\end{eqnarray*}
For any~$z$, $\theta\mapsto \rho(z- \theta) \mathbb{I}[z<\theta]$ is right-differentiable at~$\theta_0$, with right-derivative~$\psi_+(\theta_0-z)\mathbb{I}[z\leq \theta_0]=-\psi_-(z-\theta_0)\mathbb{I}[z\leq \theta_0]=|\psi_-(z-\theta_0)|\mathbb{I}[z\leq \theta_0]$, so that, for any~$z$, we have that~$L_{\theta_0}(h,z)\to 0$ as~$h$ converges to zero from above. Lemma~\ref{LemOneSidedMeanValueTheor} then shows that, for any~$h\in(0,1]$, there exists~$\lambda\in(0,1)$ such that
\begin{eqnarray*}
	|L_{\theta_0}(h,z)|
&\leq &
|\psi_-(z-\theta_0-\lambda h)| \mathbb{I}[z\leq \theta_0+\lambda h]
+
|\psi_-(z-\theta_0)|\mathbb{I}[z\leq \theta_0]
\\[2mm]
%&\leq &
%|\psi_-(z-\theta_0-1)| \mathbb{I}[z\leq \theta_0+\lambda h]
%-
%|\psi_-(z-\theta_0)|\mathbb{I}[z\leq \theta_0]
%\\[2mm]
&\leq &
|\psi_-(z-\theta_0-1)|
+
|\psi_-(z-\theta_0)|
,
\end{eqnarray*}
where this upper-bound, which does not depend on~$h$, is a $P$-integrable function of~$z$. Therefore, Lebesgue's DCT
%Dominated convergence Theorem 
entails that~$\int_{-\infty}^{\infty}L_{\theta_0}(h,z)\,dP(z)\to 0$ as~$h$ converges to zero from above. Using the fact that, for any~$z$, $\theta\mapsto \rho(z- \theta) \mathbb{I}[z\geq \theta]$ is right-differentiable at~$\theta_0$, with right-derivative~$\psi_+(\theta_0-z)\mathbb{I}[z>\theta_0]=-\psi_-(z-\theta_0)\mathbb{I}[z>\theta_0]=-|\psi_-(z-\theta_0)|\mathbb{I}[z> \theta_0]$, one can similarly show that~$\int_{-\infty}^{\infty} R_{\theta_0}(h,z)\,dP(z)\to 0$ as~$h$ converges to zero from above. From~(\ref{rightdiffO}), this establishes that~$\theta\mapsto O(\theta)$ is right-differentiable, with right-derivative
\begin{eqnarray*}
	\lefteqn{
O_+^{\rho\prime}(\theta)
=
(1-\alpha)
\int_{-\infty}^{\infty}
|\psi_-(z-\theta)| \mathbb{I}[z\leq \theta]
\,dP(z)
}
\\[2mm]
& & 
\hspace{28mm} 
-
\alpha
\int_{-\infty}^{\infty}
|\psi_-(z-\theta)|\mathbb{I}[z>\theta]
\,dP(z)
.
\end{eqnarray*}
Proceeding in the same way, 
%With the same argument, we obtain
%$$
%R'_+(\theta)
%=
%-
%\int_{-\infty}^{\infty}
%|\psi_-(z-\theta)| \mathbb{I}[z>\theta]
%\,dP(z)
%$$
%and
%$$
%R'_-(\theta)
%=
%-
%\int_{-\infty}^{\infty}
%|\psi_+(z-\theta)| \mathbb{I}[z\geq\theta]
%\,dP(z)
%%=
%%- P[Z\geq \theta]
%.
%$$
one can show that~$\theta\mapsto O(\theta)$ is left-differentiable, with left-derivative
\begin{eqnarray*}
	\lefteqn{
O_-^{\rho\prime}(\theta) 
=
(1-\alpha)
\int_{-\infty}^{\infty}
|\psi_+(z-\theta)| \mathbb{I}[z< \theta]
\,dP(z)
}
\\[2mm]
& & 
\hspace{28mm} 
-
\alpha
\int_{-\infty}^{\infty}
|\psi_+(z-\theta)| \mathbb{I}[z\geq\theta]
\,dP(z)
.
\end{eqnarray*}

(ii) 
The expressions of~$O_+^{\rho\prime}(\theta)$ and~$O_-^{\rho\prime}(\theta)$ above provide
\begin{eqnarray*}
	\lefteqn{
O_+^{\rho\prime}(\theta) - O_-^{\rho\prime}(\theta)
=
\{ (1-\alpha) |\psi_-(0)| + \alpha |\psi_+(0)| \} P[\{\theta\}]
}
\\[2mm]
& & 
\hspace{25mm} 
+
(1-\alpha)
\int_{-\infty}^{\infty}
\{ |\psi_-(z-\theta)|  - |\psi_+(z-\theta)|  \}
\mathbb{I}[z< \theta]
\,dP(z)
\\[2mm]
& & 
\hspace{25mm} 
+
\alpha
\int_{-\infty}^{\infty}
\{ |\psi_+(z-\theta)|  - |\psi_-(z-\theta)|  \}
\mathbb{I}[z> \theta]
\,dP(z)
.
\end{eqnarray*}
Since~$\psi_-(t)\leq \psi_+(t)\leq 0$ for any~$t<0$ and~$\psi_+(t)\geq \psi_-(t)\geq 0$ for any~$t>0$, we conclude that~$O_+^{\rho\prime}(\theta) - O_-^{\rho\prime}(\theta) \geq 0$.  

(iii) As we showed below Lemma~\ref{LemOneSidedMeanValueTheor}, $\psi_-(z)$ may take value zero at~$z=0$ only. This and the assumption that~$P[\{\theta\}]<1$ for any~$\theta\in\R$ ensure that
$$
\int_{-\infty}^{\infty}
|\psi_-(z-\theta)| 
\,dP(z)
\geq
\int_{-\infty}^{\infty}
|\psi_-(z-\theta)| \mathbb{I}[z\neq \theta]
\,dP(z)
>0
.
$$
Therefore, we may rewrite 
\begin{eqnarray*}
O_+^{\rho\prime}(\theta)
%&\!\!=\!\!&
%(1-\alpha)
%\int_{-\infty}^{\infty}
%|\psi_-(z-\theta)| \mathbb{I}[z\leq \theta]
%\,dP(z)
%-
%\alpha
%\int_{-\infty}^{\infty}
%|\psi_-(z-\theta)|\mathbb{I}[z>\theta]
%\,dP(z)
%\\[2mm]
&\!\!=\!\!&
\int_{-\infty}^{\infty}
|\psi_-(z-\theta)|\mathbb{I}[z\leq \theta]
\,dP(z)
-
\alpha
\int_{-\infty}^{\infty}
|\psi_-(z-\theta)| 
\,dP(z)
\\[2mm]
&\!\!=\!\!&
\left(
G^\rho(\theta)
-
\alpha
\right)
%\left(
\int_{-\infty}^{\infty}
|\psi_-(z-\theta)|
\,dP(z)
%\right)
,
\end{eqnarray*}
where~$G^\rho(\theta)$ was defined in the statement of the theorem. It trivially follows that~$O_+^{\rho\prime}(\theta)$ and~$G^\rho(\theta)-\alpha$ have the same sign. 

(iv) 
Write
\begin{eqnarray}
	G^\rho(\theta)
%=
%\frac{\int_{-\infty}^{\infty}
%|\psi_-(z-\theta)| \mathbb{I}[z\leq \theta]
%\,dP(z)}{\int_{-\infty}^{\infty}
%|\psi_-(z-\theta)| 
%\,dP(z)
%}
&\!\!=\!\!&
\frac{
-
\int_{-\infty}^\infty
\psi_-(z-\theta)\mathbb{I}[z\leq \theta] 
\,dP(z)
}{
-
\int_{-\infty}^\infty
\psi_-(z-\theta) \mathbb{I}[z\leq \theta]
\,dP(z)
+
\int_{-\infty}^{\infty}
\psi_-(z-\theta) 
\mathbb{I}[z> \theta]
\,dP(z)
}
\nonumber
%$$
%=
%\frac{
%\int_{-\infty}^{\infty}
%\{
%\psi_-(0)-\psi_-(z-\theta) 
%\}
%\mathbb{I}[z\leq \theta]
%\,dP(z)
%-
%\psi_-(0)
%P[Z\leq \theta]
%}{
%\int_{-\infty}^{\infty}
%\{
%\psi_-(0)-\psi_-(z-\theta) 
%\}
%\mathbb{I}[z\leq \theta]
%\,dP(z)
%+
%\int_{-\infty}^{\infty}
%\{
%\psi_-(z-\theta) 
%+
%\psi_-(0)
%\}
%\mathbb{I}[z> \theta]
%\,dP(z)
%-
%\psi_-(0)
%}
%,
%$$
\\[2mm]
&\!\!=\!\!&
\frac{
H^{\rho}_2(\theta)
-
\psi_-(0)
P[Z\leq \theta]
}{
H^{\rho}_2(\theta)
+
H^{\rho}_1(\theta)
-
\psi_-(0)
}
=
\frac{
H^{\rho}_2(\theta)
+
\delta
F(\theta)
}{
H^{\rho}_1(\theta)
+
H^{\rho}_2(\theta)
+
\delta
}
,
\label{expressionG}
\end{eqnarray}
where we let $\delta:=-\psi_-(0)$, $F(z):=P[Z\leq z]$,
$$
H^{\rho}_1(\theta)
:=
\int_{-\infty}^{\infty}
\{
\psi_-(z-\theta) + \psi_-(0) 
\}
\mathbb{I}[z> \theta]
\,dP(z)
%(\geq 0)
,
$$
and
$$
H^{\rho}_2(\theta)
:=
\int_{-\infty}^{\infty}
\{
\psi_-(0) - \psi_-(z-\theta) 
\}
\mathbb{I}[z\leq \theta]
\,dP(z)
%(\geq 0)
.
$$
Note that~$\delta\geq 0$ 
%(since~$-\delta$ is the limit of~$\psi_-(z)$ as $z$ goes to~$0$ from below) 
and that~$H^{\rho}_1(\theta)$ and~$H^{\rho}_2(\theta)$ are non-negative for any~$\theta$. 

Fix then~$\theta_b>\theta_a$. Since~$\psi_-$ is monotone non-decreasing, we have~$\{\psi_-(z-\theta_b) + \psi_-(0) \}\mathbb{I}[z> \theta_b]-\{\psi_-(z-\theta_a) + \psi_-(0) \}\mathbb{I}[z> \theta_a]\leq 0$, so that~$H^{\rho}_1(\theta_b)\leq H^{\rho}_1(\theta_a)$. Similarly, the monotonicity of~$\psi_-$ implies that~$\{\psi_-(0) - \psi_-(z-\theta_b) \}\mathbb{I}[z\leq \theta_b]-\{\psi_-(0) - \psi_-(z-\theta_a) \}\mathbb{I}[z\leq \theta_a]\geq 0$, so that~$H^{\rho}_2(\theta_b)\geq H^{\rho}_2(\theta_a)$. Therefore, $H^{\rho}_1$ and~$H^{\rho}_2$ are monotone non-increasing and non-decreasing, respectively. Since direct computations allow to check that
\begin{eqnarray}
\lefteqn{
\hspace{-10mm} 
\{\delta+H^{\rho}_1(\theta_a)+H^{\rho}_2(\theta_a)\}
\{\delta+H^{\rho}_1(\theta_b)+H^{\rho}_2(\theta_b)\}
\{G^\rho(\theta_b)-G^\rho(\theta_a)\}
}
\nonumber
\\[2mm]
& & 
\hspace{3mm} 
=
\{\delta(1-F(\theta_a))+H^{\rho}_1(\theta_a)\}
\{H^{\rho}_2(\theta_b)-H^{\rho}_2(\theta_a)\}
\nonumber
\\[2mm]
& & 
\hspace{17mm} 
-
\{\delta F(\theta_a)+H^{\rho}_2(\theta_a)\}
\{H^{\rho}_1(\theta_b)-H^{\rho}_1(\theta_a)\}
\label{MonotonicMiracle}
\\[2mm]
& & 
\hspace{17mm} 
+
\delta \{\delta+H^{\rho}_1(\theta_a)+H^{\rho}_2(\theta_a)\}
\{F(\theta_b)-F(\theta_a)\}
,
\nonumber
\end{eqnarray}
we conclude that~$G^{\rho}$ is monotone non-decreasing.  

Since~$\psi_-$ is monotone non-decreasing, the Monotone Convergence Theorem implies that 
$
\lim_{\theta\to -\infty}
H^\rho_2(\theta)
=
0
$
and
$
\lim_{\theta\to \infty}
H^\rho_1(\theta)
=
0
$.
Since~$H^\rho_1$ is a monotone non-increasing function of~$\theta$, $H^\rho_1(\theta)$ will stay from zero for large negative values of~$\theta$, which implies that
$$
\lim_{\theta\to -\infty}
G^\rho(\theta)
=
\lim_{\theta\to -\infty}
\frac{
H^{\rho}_2(\theta)
+
\delta
F(\theta)
}{
H^{\rho}_1(\theta)
+
H^{\rho}_2(\theta)
+
\delta
}
=
0
.
$$
Similarly, since~$H^\rho_2$ is a monotone non-decreasing function of~$\theta$, $H^\rho_2(\theta)$ will stay from zero for large positive values of~$\theta$, so that
$$
\lim_{\theta\to \infty}
(1-G^\rho(\theta))
=
\lim_{\theta\to \infty}
\frac{
H^{\rho}_1(\theta)
+
\delta
(1-F(\theta))
}{
H^{\rho}_1(\theta)
+
H^{\rho}_2(\theta)
+
\delta
}
=
0
.
$$

It remains to prove that~$G^\rho$ is right-continuous. Since the cumulative distribution function~$F$ is right-continuous, it is sufficient to show that~$H^\rho_1$ and~$H^\rho_2$ also are. To do so, fix an arbitrary monotone increasing sequence~$(c_\ell)$, $\ell=1,2,\ldots$, of positive real numbers such that the set~$\{c_0:=0,\pm c_\ell,\,\ell=1,2,\ldots\}$ contains all discontinuity points of~$\psi_-$ (recall that~$\psi_-$ has an at most countable number of discontinuity points). We can write 
\begin{eqnarray}
H^{\rho}_1(\theta)
&=&
\int_{-\infty}^{\infty}
\psi_-(z-\theta)
\mathbb{I}[z> \theta]
\,dP(z)
+
\psi_-(0) P[Z> \theta]
\nonumber
\\[2mm]
& = & 
\sum_{\ell=0}^{\infty}
\int_{-\infty}^{\infty}
\psi_-(z-\theta) 
\mathbb{I}[\theta+c_\ell<z<\theta+c_{\ell+1}]
\,dP(z)
\nonumber
\\[2mm]
& & 
\hspace{10mm} 
+
\sum_{\ell=1}^{\infty} \psi_-(c_{\ell}) P[\{\theta+c_{\ell}\}]
+
\psi_-(0) P[Z> \theta]
\label{exprH1}
\end{eqnarray}
and
\begin{eqnarray}
H^{\rho}_2(\theta)
&=&
-
\int_{-\infty}^{\infty}
\psi_-(z-\theta)
\mathbb{I}[z\leq \theta]
\,dP(z)
+
\psi_-(0) P[Z\leq \theta]
\nonumber
\\[2mm]
& = & 
-
\sum_{\ell=0}^{\infty}
\int_{-\infty}^{\infty}
\psi_-(z-\theta) 
\mathbb{I}[\theta-c_{\ell+1}<z<\theta-c_{\ell}]
\,dP(z)
\nonumber
\\[2mm]
& & 
\hspace{10mm} 
-
\sum_{\ell=0}^{\infty} \psi_-(-c_{\ell}) P[\{\theta-c_{\ell}\}]
+
\psi_-(0) P[Z\leq \theta]
.
\label{exprH2}
\end{eqnarray}
Lemma~\ref{LemDCTpsi} then provides
\begin{eqnarray*}
\lefteqn{
	\lim_{\theta\stackrel{>}{\to} \theta_0} H^\rho_1(\theta)
=
\sum_{\ell=0}^{\infty}
\int_{-\infty}^{\infty}
\psi_-(z-\theta_0) 
\mathbb{I}[\theta_0+c_\ell< z\leq\theta_0+c_{\ell+1}]
\,dP(z)
}
\\[2mm]
& & 
\hspace{25mm}
+
\psi_-(0) P[Z> \theta_0]
\\[2mm]
& & 
\hspace{10mm}
=
\sum_{\ell=0}^{\infty}
\int_{-\infty}^{\infty}
\psi_-(z-\theta_0) 
\mathbb{I}[\theta_0+c_\ell< z<\theta_0+c_{\ell+1}]
\,dP(z)
\\[2mm]
& & 
\hspace{25mm}
+
\sum_{\ell=0}^{\infty} \psi_-(c_{\ell+1}) P[\{\theta_0+c_{\ell+1}\}]
+
\psi_-(0) P[Z> \theta_0]
\\[2mm]
& & 
\hspace{10mm}
=
H^{\rho}_1(\theta_0)
\end{eqnarray*}
and
%Turning to~$H^{\rho}_2(\theta_0)$, we have
\begin{eqnarray*}
\lefteqn{
	\lim_{\theta\stackrel{>}{\to} \theta_0} H^\rho_2(\theta)
=
-
\sum_{\ell=0}^{\infty}
\int_{-\infty}^{\infty}
\psi_-(z-\theta_0) 
\mathbb{I}[\theta_0-c_{\ell+1}< z\leq\theta_0-c_{\ell}]
 \,dP(z)
}
\\[2mm]
& & 
\hspace{25mm}
+
\psi_-(0) P[Z\leq \theta_0]
\\[2mm]
& & 
\hspace{10mm}
=
-
\sum_{\ell=0}^{\infty}
\int_{-\infty}^{\infty}
\psi_-(z-\theta_0) 
\mathbb{I}[\theta_0-c_{\ell+1}< z<\theta_0-c_{\ell}]
\,dP(z)
\\[2mm]
& & 
\hspace{25mm}
-
\sum_{\ell=0}^{\infty} \psi_-(c_{\ell}) P[\{\theta_0-c_{\ell}\}]
+
\psi_-(0) P[Z\leq \theta_0]
\\[2mm]
& & 
\hspace{10mm}
=
H^{\rho}_2(\theta_0)
,
\end{eqnarray*}
which confirms that~$H^{\rho}_1$ and~$H^{\rho}_2$ are right-continuous functions.

(v) Right-continuity of~$G^{\rho}$  entails that~$G^{\rho}(\theta^\rho_\alpha(P))\geq \alpha$. Hence, monotonicity of~$G^{\rho}$ provides~$G^{\rho}(\theta)\geq \alpha$ for any~$\theta\geq \theta^\rho_\alpha(P)$, which implies that~$O^{\rho\prime}_+(\theta)\geq 0$  for any~$\theta\geq \theta^\rho_\alpha(P)$. Lemma~\ref{LemOneSidedMeanValueTheor} therefore entails that~$O^{\rho}(\theta)\geq O^{\rho}(\theta^\rho_\alpha(P))$  for any~$\theta\geq \theta^\rho_\alpha(P)$. Now, the definition of~$\theta^\rho_\alpha(P)$ ensures that~$G^{\rho}(\theta)< \alpha$ for any~$\theta< \theta^\rho_\alpha(P)$, which implies that~$O^{\rho\prime}_+(\theta)<0$ for any such~$\theta$. Lemma~\ref{LemOneSidedMeanValueTheor} therefore also yields that~$O^{\rho}(\theta)> O^{\rho}(\theta^\rho_\alpha(P))$ for any~$\theta< \theta^\rho_\alpha(P)$, which establishes the result. 

(vi) In Part~(iv) of the result, we showed that 
$$
G^{\rho}(\theta)
=
\frac{
H^{\rho}_2(\theta)
+
\delta
F(\theta)
}{
H^{\rho}_1(\theta)
+
H^{\rho}_2(\theta)
+
\delta
}
$$
is right-continuous. Assume now that~$\psi_-$ is continuous over~$\R$ or that~$P$ is non-atomic. Then~$\theta\mapsto \delta
F(\theta)$ is left-continuous (note that if~$\psi_-$ is continuous at zero, then~$\delta=0$). Now, applying Lemma~\ref{LemDCTpsi} to~(\ref{exprH1}) %and~(\ref{exprH2}) 
yields
\begin{eqnarray*}
\lefteqn{
	\lim_{\theta\stackrel{<}{\to} \theta_0} H^\rho_1(\theta_0)
=
\sum_{\ell=0}^{\infty}
\int_{-\infty}^{\infty}
\psi_-(z-\theta_0+0) 
\mathbb{I}[\theta_0+c_\ell\leq z<\theta_0+c_{\ell+1}]
\,dP(z)
}
\\[2mm]
& & 
\hspace{25mm}
+
\psi_-(0) P[Z\geq \theta_0]
\\[2mm]
& & 
\hspace{10mm}
=
\sum_{\ell=0}^{\infty}
\int_{-\infty}^{\infty}
\psi_-(z-\theta_0) 
\mathbb{I}[\theta_0+c_\ell< z<\theta_0+c_{\ell+1}]
\,dP(z)
\\[2mm]
& & 
\hspace{25mm}
+
\sum_{\ell=0}^{\infty} \psi_-(c_{\ell}+0) P[\{\theta_0+c_{\ell}\}]
+
\psi_-(0) P[Z\geq \theta_0]
%\\[2mm]
%& & 
%\hspace{10mm}
,
\end{eqnarray*}
which provides
\begin{eqnarray*}
\lefteqn{
	\lim_{\theta\stackrel{<}{\to} \theta_0} H^\rho_1(\theta_0)
=
\sum_{\ell=0}^{\infty}
\int_{-\infty}^{\infty}
\psi_-(z-\theta_0) 
\mathbb{I}[\theta_0+c_\ell< z<\theta_0+c_{\ell+1}]
\,dP(z)
}
\\[2mm]
& & 
\hspace{25mm}
+
\sum_{\ell=1}^{\infty} \psi_-(c_{\ell}+0) P[\{\theta_0+c_{\ell}\}]
+
\psi_-(0) P[Z> \theta_0]
\\[2mm]
& & 
\hspace{10mm}
=
H^{\rho}_1(\theta_0)
+
\sum_{\ell=1}^{\infty} (\psi_-(c_{\ell}+0)-\psi_-(c_{\ell})) P[\{\theta_0+c_{\ell}\}]
=
H^{\rho}_1(\theta_0)
.
\end{eqnarray*}
Similarly, applying Lemma~\ref{LemDCTpsi} to~(\ref{exprH2}) 
gives
\begin{eqnarray*}
\lefteqn{
	\lim_{\theta\stackrel{<}{\to} \theta_0} H^\rho_2(\theta_0)
=
-
\sum_{\ell=0}^{\infty}
\int_{-\infty}^{\infty}
\psi_-(z-\theta_0+0) 
\mathbb{I}[\theta_0-c_{\ell+1}\leq z<\theta_0-c_{\ell}]
 \,dP(z)
}
\\[2mm]
& & 
\hspace{25mm}
+
\psi_-(0) P[Z< \theta_0]
\\[2mm]
& & 
\hspace{0mm}
=
-
\sum_{\ell=0}^{\infty}
\int_{-\infty}^{\infty}
\psi_-(z-\theta_0) 
\mathbb{I}[\theta_0-c_{\ell+1}< z<\theta_0-c_{\ell}]
\,dP(z)
\\[2mm]
& & 
\hspace{8mm}
-
\sum_{\ell=0}^{\infty} \psi_-(-c_{\ell+1}+0) P[\{\theta_0-c_{\ell+1}\}]
+
\psi_-(0) P[Z< \theta_0]
%\\[2mm]
%& & 
%\hspace{0mm}
,
\end{eqnarray*}
which entails that
\begin{eqnarray*}
\lefteqn{
	\lim_{\theta\stackrel{<}{\to} \theta_0} H^\rho_2(\theta_0)
=
-
\sum_{\ell=0}^{\infty}
\int_{-\infty}^{\infty}
\psi_-(z-\theta_0) 
\mathbb{I}[\theta_0-c_{\ell+1}< z<\theta_0-c_{\ell}]
\,dP(z)
}
\\[2mm]
& & 
\hspace{8mm}
-
\sum_{\ell=1}^{\infty} \psi_-(-c_{\ell}+0) P[\{\theta_0-c_{\ell}\}]
-
\psi_-(0) P[\{\theta_0\}]
+
\psi_-(0) P[Z\leq \theta_0]
\\[2mm]
& & 
\hspace{0mm}
=
H^{\rho}_2(\theta_0)
+
\sum_{\ell=1}^{\infty} (\psi_-(-c_{\ell})-\psi_-(-c_{\ell}+0)) P[\{\theta_0-c_{\ell}\}]
=
H^{\rho}_2(\theta_0)
.
\end{eqnarray*}
We conclude that~$H_1^{\rho}$ and~$H_2^{\rho}$ are also left-continuous, which implies that~$G^{\rho}$ is continuous. 
\cqfd 
\vspace{3mm}

Lemmas~\ref{LemPourLemMDepthStrictMonotonicity}--\ref{LemPourMonotoneExpect} below will be needed in subsequent proofs.

%
%
%
%\newpage
%
%
%

\begin{Lem}
\label{LemPourLemMDepthStrictMonotonicity}
Fix~$\rho\in\mathcal{C}$ and~$P\in\mathcal{P}^\rho_1$. Denote as~$F$ the cumulative distribution function associated with~$P$. Let~$\theta_a<\theta_b$ with~$F(\theta_a)>0$ and $F(\theta_b-0)-F(\theta_a)>0$, where~$F(\theta_b-0)$ is the limit of~$F(\theta)$ when~$\theta$ converges to~$\theta_b$ from below. 	Then~$G^\rho(\theta_b)>G^\rho(\theta_a)$. 
\end{Lem}

{\sc Proof of Lemma~\ref{LemPourLemMDepthStrictMonotonicity}.}
It directly follows from~(\ref{MonotonicMiracle}) that
\begin{eqnarray}
\nonumber
\lefteqn{
\hspace{0mm} 
\{\delta+H^{\rho}_1(\theta_a)+H^{\rho}_2(\theta_a)\}
\{\delta+H^{\rho}_1(\theta_b)+H^{\rho}_2(\theta_b)\}
\{G^\rho(\theta_b)-G^\rho(\theta_a)\}
}
\\[2mm]
%& & 
%\hspace{3mm} 
%=
%\{\delta(1-F(\theta_a))+H^{\rho}_1(\theta_a)\}
%\{H^{\rho}_2(\theta_b)-H^{\rho}_2(\theta_a)\}
%\\[2mm]
%& & 
%\hspace{17mm} 
%-
%\{\delta F(\theta_a)+H^{\rho}_2(\theta_a)\}
%\{H^{\rho}_1(\theta_b)-H^{\rho}_1(\theta_a)\}
%\\[2mm]
%& & 
%\hspace{17mm} 
%+
%\delta \{\delta+H^{\rho}_1(\theta_a)+H^{\rho}_2(\theta_a)\}
%\{F(\theta_b)-F(\theta_a)\}
%\\[2mm]
& & 
\hspace{6mm} 
\geq 
\{\delta(1-F(\theta_a))+H^{\rho}_1(\theta_a)\}
\{H^{\rho}_2(\theta_b)-H^{\rho}_2(\theta_a)\}
+
\delta^2 \{F(\theta_b)-F(\theta_a)\}
.
\label{basiseq}
\end{eqnarray}
If~$\delta=-\psi_-(0)>0$, then this readily yields
\begin{eqnarray*}
\lefteqn{
\hspace{-44mm} 
\{\delta+H^{\rho}_1(\theta_a)+H^{\rho}_2(\theta_a)\}
\{\delta+H^{\rho}_1(\theta_b)+H^{\rho}_2(\theta_b)\}
\{G^\rho(\theta_b)-G^\rho(\theta_a)\}
}
\\[2mm]
& & 
\hspace{-23mm} 
\geq 
\delta^2 \{F(\theta_b)-F(\theta_a)\}
>0
,
\end{eqnarray*}
which implies that~$G^\rho(\theta_b)>G^\rho(\theta_a)$. We may therefore restrict to the case~$\psi_-(0)=0$ (recall that~$\delta=-\psi_-(0)\geq 0$), under which~(\ref{basiseq}) rewrites
\begin{eqnarray*}
	\lefteqn{
\hspace{-13mm} 
\{H^{\rho}_1(\theta_a)+H^{\rho}_2(\theta_a)\}
\{H^{\rho}_1(\theta_b)+H^{\rho}_2(\theta_b)\}
\{G^\rho(\theta_b)-G^\rho(\theta_a)\}
}
\\[2mm]
& & 
\hspace{13mm} 
\geq 
H^{\rho}_1(\theta_a)
\{H^{\rho}_2(\theta_b)-H^{\rho}_2(\theta_a)\}
.
\end{eqnarray*}
Since~$1-F(\theta_a)\geq F(\theta_b-0)-F(\theta_a)>0$ and~$\psi_-(z)>\psi_-(0)=0$ for any~$z>0$ (recall that~$\psi_-(z)\neq 0$ for~$z\neq 0$), we have 
$$
H^{\rho}_1(\theta_a)
=
\int_{-\infty}^{\infty}
\{
\psi_-(z-\theta_a) + \psi_-(0) 
\}
\mathbb{I}[z> \theta_a]
\,dP(z)
>
0
.
$$
Now, for any~$z<\theta_b$, we have~$\psi_-(z-\theta_b)<\psi_-(0)=0$, whereas for any~$z$, we have~$\psi_-(z-\theta_a)\geq \psi_-(z-\theta_b)$. Therefore, using the fact that~$F(\theta_b-0)-F(\theta_a)>0$, we have
\begin{eqnarray*}	
\lefteqn{
H^{\rho}_2(\theta_b)-H^{\rho}_2(\theta_a)
=
\int_{-\infty}^{\infty}
\big\{
\{
\psi_-(0) - \psi_-(z-\theta_b) 
\}
\mathbb{I}[z< \theta_b]
}
\\[2mm]
& & 
\hspace{43mm} 
-
\{
\psi_-(0) - \psi_-(z-\theta_a) 
\}
\mathbb{I}[z\leq \theta_a]
\big\}
\,dP(z)
> 0
\\[2mm]
& & 
\hspace{23mm} 
\geq 
\int_{-\infty}^{\infty}
\{
\psi_-(0) - \psi_-(z-\theta_b) 
\}
\mathbb{I}[\theta_a< z< \theta_b]
\,dP(z)
%> 0
.
\end{eqnarray*}
We conclude that~$G^\rho(\theta_b)>G^\rho(\theta_a)$. 
\cqfd
\vspace{3mm}

\begin{Lem}
\label{LemPourMonotoneExpect}
Fix~$\rho\in\mathcal{C}$ and~$P\in\mathcal{P}^\rho_1$. Assume that~$\psi_-$ is monotone strictly increasing over~$\R$. Let~$\theta'<\theta''$ with~$\theta',\theta''\in C_P=\{ \theta \in \R : \min( P[Z\leq \theta] , P[Z\geq \theta] ) > 0\}$. 	Then~$G^\rho$ is monotone strictly increasing over~$[\theta',\theta'']$. 
\end{Lem}

{\sc Proof of Lemma~\ref{LemPourMonotoneExpect}.}
Fix~$\theta_a,\theta_b$ with~$\theta'<\theta_a<\theta_b<\theta''$, so that~$P[Z>\theta_a]>0$ and~$P[Z< \theta_b]>0$. Recalling that all terms in the right-hand side of~(\ref{MonotonicMiracle}) are non-negative, we obtain
\begin{eqnarray*}
\lefteqn{
\hspace{-10mm} 
\{\delta+H^{\rho}_1(\theta_a)+H^{\rho}_2(\theta_a)\}
\{\delta+H^{\rho}_1(\theta_b)+H^{\rho}_2(\theta_b)\}
\{G^\rho(\theta_b)-G^\rho(\theta_a)\}
}
\\[2mm]
%& & 
%\hspace{3mm} 
%=
%\{\delta(1-F(\theta_a))+H^{\rho}_1(\theta_a)\}
%\{H^{\rho}_2(\theta_b)-H^{\rho}_2(\theta_a)\}
%\\[2mm]
%& & 
%\hspace{17mm} 
%-
%\{\delta F(\theta_a)+H^{\rho}_2(\theta_a)\}
%\{H^{\rho}_1(\theta_b)-H^{\rho}_1(\theta_a)\}
%\\[2mm]
%& & 
%\hspace{17mm} 
%+
%\delta \{\delta+H^{\rho}_1(\theta_a)+H^{\rho}_2(\theta_a)\}
%\{F(\theta_b)-F(\theta_a)\}
%\\[2mm]
& & 
\hspace{3mm} 
\geq 
\{\delta(1-F(\theta_a))+H^{\rho}_1(\theta_a)\}
\{H^{\rho}_2(\theta_b)-H^{\rho}_2(\theta_a)\}
.
\end{eqnarray*}
Since~$P[Z>\theta_a]>0$ and~$\psi_-(z)>-\psi_-(0)$ for any~$z>0$, we have 
\begin{eqnarray*}
	\lefteqn{
\delta(1-F(\theta_a))+H^{\rho}_1(\theta_a)
}
\\[2mm]
& & 
\hspace{5mm} 
\geq 
H^{\rho}_1(\theta_a)
=
\int_{-\infty}^{\infty}
\{
\psi_-(z-\theta_a) + \psi_-(0) 
\}
\mathbb{I}[z> \theta_a]
\,dP(z)
>
0
.
\end{eqnarray*}
Now, for any~$z<\theta_b$, we have~$\psi_-(0)>\psi_-(z-\theta_b)$, whereas for any~$z$, we have~$\psi_-(z-\theta_a)>\psi_-(z-\theta_b)$. Therefore, using the fact that~$P[Z< \theta_b]>0$, we have
\begin{eqnarray*}	
\lefteqn{
H^{\rho}_2(\theta_b)-H^{\rho}_2(\theta_a)
=
\int_{-\infty}^{\infty}
\big\{
\{
\psi_-(0) - \psi_-(z-\theta_b) 
\}
\mathbb{I}[z< \theta_b]
}
\\[2mm]
& & 
\hspace{40mm} 
-
\{
\psi_-(0) - \psi_-(z-\theta_a) 
\}
\mathbb{I}[z\leq \theta_a]
\big\}
\,dP(z)
> 0
.
\end{eqnarray*}
We conclude that~$G^\rho(\theta_b)>G^\rho(\theta_a)$. This shows that~$G^\rho$ is monotone strictly increasing over~$(\theta',\theta'')$, hence over~$[\theta',\theta'']$. 
%
%Now, assume there exists~$h\in(0,\theta''-\theta')$ such that~$G^\rho(\theta')=G^\rho(\theta'+h)$. Then we have~$G^\rho(\theta'+h/2)<G^\rho(\theta'+h)=G^\rho(\theta')$, which violates the fact that~$G^\rho$ is monotone non-decreasing. Therefore, $G^\rho$ is monotone strictly increasing in~$[\theta',\theta'')$. The same argument shows that~$G^\rho$ must be monotone strictly increasing in~$[\theta',\theta'']$.
\cqfd
\vspace{3mm}

%%%%%%%%%

%\subsection{Proofs of Section~\ref{secproposedmultiquantile}}

In the rest of this appendix, $\theta_\alpha(Z)$, $H_{\alpha,u}^\rho(Z)$,\,\ldots, will respectively stand for~$\theta_\alpha(P)$, $H_{\alpha,u}^\rho(P)$,\,\ldots, where~$P$ is the distribution of~$Z$. 
\vspace{2mm}

{\sc Proof of Theorem~\ref{Theoraffequiv}}. Let~$\Zb$ be a random $d$-vector with distribution~$P$. First note that the assumption on~$\rho$ entails that 
%$$
%\psi_-(\lambda a)
%=
%\lim_{t\to (\lambda a)^-} \frac{\rho(t)-\rho(\lambda a)}{t-\lambda a}
%=
%\lim_{s\to a^-} \frac{\rho(\lambda s)-\rho(\lambda a)}{\lambda (s-a)}
%=
%\lambda^{r-1} 
%\lim_{s\to a^-} \frac{\rho(s)-\rho(a)}{s-a}
%=
%\lambda^{r-1} 
%\psi_-(a)
%.
%$$
$\psi_-(\lambda t)=\lambda^{r-1}\psi_-(t)$ for any~$t\in\R$ and~$\lambda>0$, which provides
\begin{eqnarray}
\label{idabove}
	\lefteqn{
	\hspace{18mm} 
\frac{{\rm E}[|\psi_-(\ub_\Ab'(\Ab\Zb+\bb)-\theta)| \mathbb{I}[\ub_\Ab'(\Ab\Zb+\bb)\leq \theta]]}{{\rm E}[|\psi_-(\ub_\Ab'(\Ab\Zb+\bb)-\theta)|]}
=
}
%$$
%=
%\frac{{\rm E}[|\psi_-(\frac{\ub'Z- \{ \|(\Ab^{-1})'u\| \theta - \ub' \Ab^{-1} b \} }{\|(\Ab^{-1})'u\|})| \mathbb{I}[\frac{\ub'Z- \{ \|(\Ab^{-1})'u\| \theta - \ub' \Ab^{-1} b \} }{\|(\Ab^{-1})'u\|}\leq 0]]}{{\rm E}[|\psi_-(\frac{\ub'Z- \{ \|(\Ab^{-1})'u\| \theta - \ub' \Ab^{-1} b \} }{\|(\Ab^{-1})'u\|})|]}
%$$
\\[2mm]
 & & 
 \hspace{-7mm} 
\frac{{\rm E}[|\psi_-(\ub'\Zb- \{ \|(\Ab^{-1})'\ub\| \theta - \ub' \Ab^{-1} \bb \} )| \mathbb{I}[\ub'\Zb- \{ \|(\Ab^{-1})'\ub\| \theta - \ub' \Ab^{-1} \bb \} \leq 0]]}{{\rm E}[|\psi_-(\ub'\Zb- \{ \|(\Ab^{-1})'\ub\| \theta - \ub' \Ab^{-1} \bb \})|]}
\cdot
\nonumber
\end{eqnarray}
Denote as~$S_\alpha$ the set of real numbers~$\phi$ such that 
$$
\frac{{\rm E}[|\psi_-(\ub'\Zb- \phi )| \mathbb{I}[\ub'\Zb- \phi \leq 0]]}{{\rm E}[|\psi_-(\ub'\Zb-\phi)|]}
\geq \alpha
$$
and as~$T_\alpha$ the set of real numbers~$\theta$ such that 
$$
\frac{{\rm E}[|\psi_-(\ub_\Ab'(\Ab\Zb+\bb)-\theta)| \mathbb{I}[\ub_\Ab'(\Ab\Zb+\bb)\leq \theta]]}{{\rm E}[|\psi_-(\ub_\Ab'(\Ab\Zb+\bb)-\theta)|]}
\geq \alpha
.
$$
It follows from~(\ref{idabove}) that~$\theta\in T_\alpha$ if and only if~$\|(\Ab^{-1})'\ub\| \theta - \ub' \Ab^{-1} \bb \in S_\alpha$.
%, which is if and only if~$\|(\Ab^{-1})'\ub\| \theta - \ub' \Ab^{-1} \bb \in \|(\Ab^{-1})'\ub\| T_\alpha - \ub' \Ab^{-1} \bb$. 
Thus, $S_\alpha=\|(\Ab^{-1})'\ub\| T_\alpha - \ub' \Ab^{-1} \bb$. Since 
$
\theta^\rho_{\alpha}(\ub'\Zb)
=
\inf S_\alpha
$
and
$
\theta^\rho_{\alpha}(\ub_\Ab'(\Ab\Zb+\bb))
=
\inf T_\alpha
$
by definition, this implies that
$$
\theta^\rho_{\alpha}(\ub'\Zb)
%=
%\inf S_\alpha
%=
%\|(\Ab^{-1})'u\| (\inf T_\alpha) - \ub' \Ab^{-1} b
=
\|(\Ab^{-1})'\ub\| \theta^\rho_{\alpha}(\ub_\Ab'(\Ab\Zb+\bb)) - \ub' \Ab^{-1} \bb
.
$$ 
We conclude that $H^{\rho}_{\alpha,\ub_\Ab}(\Ab\Zb+\bb)$ collects the $d$-vectors~$z$ satisfying
$$
\ub_\Ab' \zb 
\geq 
\frac{\theta^\rho_{\alpha}(\ub'\Zb)}{\|(\Ab^{-1})'\ub\|} + \frac{\ub' \Ab^{-1} \bb}{\|(\Ab^{-1})'\ub\|}
,
$$
or equivalently, 
$
\ub' (\Ab^{-1} (\zb-\bb)) 
\geq 
\theta^\rho_{\alpha}(\ub'\Zb)
$. This establishes the result. 
\cqfd
\vspace{3mm}

%%%%%%%%%%%%%%%%%%%%%%%%%%

Let~$\Zb$ be a random $d$-vector with distribution~$P\in\mathcal{P}^\rho_d$. In the subsequent proofs, we let 
\begin{equation}
\label{definGu}
G^\rho_\ub(\theta)
:=
\frac{{\rm E}[|\psi_-(\ub'\Zb-\theta)| \mathbb{I}[\ub'\Zb-\theta\leq 0]]}{{\rm E}[|\psi_-(\ub'\Zb-\theta)|]}
\cdot
\end{equation}
With this notation, $\theta^\rho_\alpha(P_\ub)=\theta^\rho_\alpha(\ub'\Zb)$ is given by~$\inf\big\{\theta\in\R: G^\rho_\ub(\theta)\geq\alpha\big\}$; see Theorem~\ref{Theorlemjones}(v). 
\vspace{2mm}

{\sc Proof of Theorem~\ref{TheorSupportAndcompactness}}.
Fix~$\alpha\in(0,1)$. By definition,~$R_{\alpha}^{\rho}(P)$ is an intersection of closed convex subsets of~$\R^d$, so that it is itself closed and convex. Since the affine-equivariance relation~$
R^{\rho}_{\alpha}(P_{\Ab,\bb})
=
\Ab R^{\rho}_{\alpha}(P)+\bb
$ 
(for~$\rho\in\mathcal{C}_{\rm aff}$) is a direct corollary of Theorem~\ref{Theoraffequiv}, it only remains to show that (i)~$R_{\alpha}^{\rho}(P)\subset C_P$ and that (ii) $R_{\alpha}^{\rho}(P)$ is bounded.  
Let us start with~(i). Fix~$\zb\notin C_P$ and let~$\Zb$ be a random \mbox{$d$-vector} with distribution~$P$. Then there exists~$\ub_0\in\mathcal{S}^{d-1}$ such that $P[\ub_0'\Zb\leq \ub_0'\zb]=0$, so that~$G_{\ub_0}^\rho(\ub_0'\zb)=0$. Right-continuity of~$G_{\ub_0}^\rho$ then entails that~$\ub_0'\zb<\theta_{\alpha}^\rho(\ub_0'\Zb)$. In other words, $\zb\notin H_{\alpha,\ub_0}^{\rho}(P)$, which implies that~$\zb\notin R_{\alpha}^{\rho}(P)$. 
(ii) Fix~$\zb\in R_{\alpha}^{\rho}(P)$. For any~$j\in\{1,\ldots,d\}$, we must have
\vspace{-.5mm} 
 $\zb\in H_{\alpha,\eb_j}^{\rho}(P)\cap H_{\alpha,-\eb_j}^{\rho}(P)$, where~$\eb_j$ denotes the $j$th vector of the canonical basis of~$\R^d$. This implies that, for any~$j$, we have $z_j\leq  \theta_{\alpha}^\rho(Z_j)$ and~$-z_j\leq  \theta_{\alpha}^\rho(-Z_j)$, that is $z_j\in [-\theta_{\alpha}^\rho(-Z_j),\theta_{\alpha}^\rho(Z_j)]$. Since the definition in~(\ref{Mquantdef}) entails that~$\theta_{\alpha}^\rho(Y)$ is finite for any random variable~$Y$, it follows that~$R_{\alpha}^{\rho}(P)$ is a subset of the bounded hyperrectangle~$\bigtimes_{j=1}^d [-\theta_{\alpha}^\rho(-Z_j),\theta_{\alpha}^\rho(Z_j)]$, hence is itself bounded.    
\cqfd
\vspace{3mm}

%Let~$\theta_0$ be such that~$P[Z\leq\theta_0]=0$. Then,~$G^\rho(\theta_0)=0$, so that (from right-continuity), one has~$\theta^\rho_{\alpha}>\theta_0$. \\     
%Let~$\theta_0$ be such that~$P[Z\geq \theta_0]=0$. Then,~$G^\rho(\theta_0)=1$ and~$G^\rho$ is continuous at~$\theta_0$ (since~$P[\{\theta_0\}=0$), so that~$\theta^\rho_{\alpha}<\theta_0$. \\    
%Hence, $\{\theta:\min(P[Z\leq\theta],P[Z\geq\theta])=0\}\subset\R\setminus\{\theta^\rho_{\alpha}:\alpha\in(0,1)\}$.\\
%Therefore, $\{\theta^\rho_{\alpha}:\alpha\in(0,1)\}\subset\{\theta:\min(P[Z\leq\theta],P[Z\geq\theta])>0\}$.

%%%%%%%%%%%%%%%%%%%%%%%%%%

%\subsection{Proofs of Section~\ref{secMdepth}}
 
{\sc Proof of Theorem~\ref{TheorRegionsOrig}.} In this proof, we let~$\tilde{R}^\rho_\alpha(P):=\{\yb\in\R^d:M\!D^\rho(\yb,P)\geq \alpha\}$. Assume first that $\zb\in R^{\rho}_{\alpha}(P)$. Then 
$
M\!D^\rho(\zb,P)
=
\sup
\{
\beta>0 : \zb\in R^{\rho}_{\beta}(P)
\}
\geq \alpha
,
$
so that~$\zb\in \tilde{R}^\rho_\alpha(P)$. 
\vspace{-.4mm}
Now, assume that~$\zb\notin R^{\rho}_{\alpha}(P)$. Then there exists~$\ub$ such that~$\ub'\zb< \theta^\rho_{\alpha}(P_\ub)$. By definition of~$\theta^\rho_{\alpha}(P_\ub)$, we must have~$G^\rho_\ub(\ub'\zb)<\alpha$. Fix then~$\alpha'\in(G^\rho_\ub(\ub'z),\alpha)$. Since~$G^\rho_\ub$ is right-continuous, there exists~$\delta\in(0,\theta^\rho_{\alpha}(P_\ub)-\ub'\zb)$ such that~$G^\rho_\ub(t)<\alpha'$ for any~$r\in[ \ub'\zb,\ub'\zb+\delta]$. The monotonicity of~$G^\rho_\ub$ then implies that~$\ub'\zb<\ub'\zb+\delta\leq \theta^\rho_{\alpha'}(P_\ub)$, which entails that~$\zb\notin R^{\rho}_{\alpha'}(P)$. We conclude that~$M\!D^\rho(\zb,P)\leq \alpha'<\alpha$, so that~$\zb\notin \tilde{R}^\rho_\alpha(P)$.   
\cqfd
\vspace{3mm}

%%%%%%%%%%%%%%%%%%%%%%%%%%

%{\sc Proof of Corollary~\ref{Corolquasiconc}.}
%With~$\alpha:=\min(M\!D^{\rho}(z_0,P),M\!D^{\rho}(z_1,P))$, both~$z_0$ and~$z_1$ belong to~$R^{\rho}_{\alpha}(P)$. Convexity of~$R^{\rho}_{\alpha}(P)$ then entails that~$M\!D^{\rho}((1-\lambda)z_0+\lambda z_1,P)\geq \alpha$. 
%\cqfd
%\vspace{3mm}

%%%%%%%%%%%%%%%%%%%%%%%%%%

{\sc Proof of Theorem~\ref{Theordeepestpoint}}.
Let~$\alpha_*:=\sup_{\zb\in\R^d} M\!D^{\rho}(\zb,P)$. Since the result trivially holds if~$\alpha_*=0$, we may restrict to the case~$\alpha_*>0$. For any positive integer~$n$, there exists~$\zb_n\in\R^d$ such that~$M\!D^{\rho}(\zb_n,P)\geq \alpha_*-(1/n)$. By dropping finitely many terms in the resulting sequence~$(\zb_n)$, one obtains a sequence with values in~$R^{\rho}_{\alpha_*/2}(P)$ (recall that Theorem~\ref{TheorRegionsOrig} states that~$R^{\rho}_{\alpha}(P)$ collects the~$\zb$'s such that~$M\!D^\rho(\zb,P)\geq \alpha$). Thus, from the compactness of~$R^{\rho}_{\alpha_*/2}(P)$ (Theorem~\ref{TheorSupportAndcompactness}), 
\vspace{-.4mm}
there exists a subsequence~$(\zb_{n_\ell})$ that converges in~$R^{\rho}_{\alpha_*/2}(P)$, to~$\zb_*$ say. 
\vspace{-.4mm}
Fix then an arbitrary~$\varepsilon>0$. For~$\ell$ large enough,~$\zb_{n_\ell}$ belongs to the closed set~$R^{\rho}_{\alpha_*-\varepsilon}(P)$, so that~$\zb_*$ also belongs to~$R^{\rho}_{\alpha_*-\varepsilon}(P)$. We therefore proved that $M\!D^\rho(\zb_*,P)\geq \alpha_*-\varepsilon$ for any~$\varepsilon>0$, which establishes the result.      
\cqfd
\vspace{3mm}

%%%%%%%%%%%%%%%%%%%%%%%%%%

The proof of Theorem~\ref{Theordepthaxioms} requires the following lemma.

\begin{Lem}
\label{Lemmaxsymm}	
Fix~$\alpha\in(0,1)$ and let~$Z$ be a random variable with a distribution~$P(\in\mathcal{P}_1^d)$ that is
%centrally 
symmetric about~$\theta_*$. 
%$($in the sense that~$-(Z-\theta_*)$ and $Z-\theta_*$ are equal in distribution$)$. 
Then, unless it is void, the interval~$[\theta^\rho_{\alpha}(Z),-\theta^\rho_{\alpha}(-Z)]$ contains~$\theta_*$.
\end{Lem}

{\sc Proof of Lemma~\ref{Lemmaxsymm}}.
Since the symmetry assumption ensures that~$-Z$ and $Z-2\theta_*$ are equal in distribution, we have 
$$
G^\rho_{-1}(\theta)
=
\frac{{\rm E}[|\psi_-(-Z-\theta)| \mathbb{I}[-Z\leq \theta]]}{{\rm E}[|\psi_-(-Z-\theta)|]}
=
G^\rho_{1}(2\theta_*+\theta)
%=
%\frac{{\rm E}[|\psi_-(Z-(2\theta_*+\theta))| \mathbb{I}[Z-(2\theta_*+\theta)\leq 0]]}{{\rm E}[|\psi_-(Z-(2\theta_*+\theta))|]}
,
$$ 
where we used the notation defined in~(\ref{definGu}). 
Consequently, 
\begin{eqnarray*}
\theta^\rho_{\alpha}(-Z)
&\!\!=\!\!&
\inf\big\{\theta\in\R: G^\rho_{-1}(\theta)\geq\alpha\big\}
\\[1mm]
&\!\!=\!\!&
- \sup \big\{\theta\in\R: G^\rho_{-1}(-\theta)\geq\alpha\big\}
\\[1mm]
&\!\!=\!\!&
- \sup \big\{\theta\in\R: G^\rho_{1}(2\theta_*-\theta)\geq\alpha\big\}
.
\end{eqnarray*}
Therefore, the interval~$[\theta^\rho_{\alpha}(Z),-\theta^\rho_{\alpha}(-Z)]$ rewrites~$[\inf \mathcal{I}_+,\sup \mathcal{I}_-]$, with 
$$
\mathcal{I}_\pm
:=
\big\{\theta\in\R: G^\rho_{1}(\theta_*\pm(\theta-\theta_*))\geq\alpha\big\}
.
$$
The result then follows from the fact that~$\mathcal{I}_-$ is obtained from~$\mathcal{I}_+$ by applying a reflection with respect to~$\theta_*$.
\cqfd
\vspace{3mm}

%%%%%%%%%%%%%%%%%%%%%%%%%%

{\sc Proof of Theorem~\ref{Theordepthaxioms}}.
(i) The claim directly follows from the affine-equivariance result in Theorem~\ref{TheorSupportAndcompactness}. 
(ii) Letting~$\Zb$ be a random $d$-vector with distribution~$P$, assume by contradiction that~$M\!D^{\rho}(\zb,P)>M\!D^{\rho}(\thetab_*,P)$ for some~$\zb\in\R^d$. Thus, there exists~$\alpha\in(0,1)$ such that~$z\in R_\alpha^\rho(P)$ and~$\theta_*\notin R_\alpha^\rho(P)$. In other words, while~$\ub'\zb\geq \theta^\rho_{\alpha}(\ub'\Zb)$ for any~$\ub\in\mathcal{S}^{d-1}$, there exists~$\ub_0$ such that~$\ub_0'\thetab_*< \theta^\rho_{\alpha}(\ub_0'\Zb)$. This implies that~$[\theta^\rho_{\alpha}(\ub_0'\Zb),-\theta^\rho_{\alpha}(-\ub_0'\Zb)]$ is non-empty (it contains~$\ub'\zb$) but does not contain~$\ub_0'\thetab_*$. Since, by assumption, $\ub_0'\Zb$ is centrally symmetric about~$\ub_0'\thetab_*$, this contradicts Lemma~\ref{Lemmaxsymm}. 
(iii) Fix~$r_2>r_1\geq 0$. The quasi-concavity property stated in the paragraph below Theorem~\ref{TheorRegionsOrig} readily provides~$M\!D^{\rho}(\thetab+r_1 \ub,P)\geq \min(M\!D^{\rho}(\thetab,P),M\!D^{\rho}(\thetab+r_2 \ub,P))=M\!D^{\rho}(\thetab+r_2 \ub,P)$, so that~$r\mapsto M\!D^{\rho}(\thetab+r \ub,P)$ is indeed monotone non-increasing over~$\R^+$. Now, fix~$r>r_\ub(P)$. By definition of~$r_\ub(P)$, $\thetab+r \ub$ does not belong to~$C_P$. Theorem~\ref{TheorSupportAndcompactness} then ensures that there is no~$\alpha\in(0,1)$ for which~$\thetab+r \ub\in R_{\alpha}^{\rho}(P)$. Thus, by definition, $M\!D^{\rho}(\thetab+r \ub,P)=0$.    
(iv) Fix~$\varepsilon>0$. In view of Theorem~\ref{TheorSupportAndcompactness}, there exists~$C>0$ such that~$R_{\alpha}^{\varepsilon}(P)$ is included in~$\{\zb\in\R^d:\|\zb\|\leq C\}$. Consequently, Theorem~\ref{TheorRegionsOrig} entails that, as soon as~$\|\zb\|>C$, one has~$M\!D^{\rho}(\zb,P)<\varepsilon$, as was to be shown.    
\cqfd
\vspace{3mm}

%%%%%%%%%%%%%%%%%%%%%%%%%%

{\sc Proof of Theorem~\ref{reginter}.}
Fix~$z\in\R^d$ and let~$\alpha=\inf_{\ub\in\mathcal{S}^{d-1}} G^\rho_\ub(\ub'\zb)$ be the right-hand side of~(\ref{mequ}). Then, 
$
G^\rho_\ub(\ub'\zb)
%=
%\frac{{\rm E}[\psi_-(\ub'Z-\ub'z))]}{{\rm E}[|\psi(\ub'Z-\ub'z))|]}
\geq 
\alpha 
$
for any~$\ub\in\mathcal{S}^{d-1}$. By definition, we must then have~$\ub'\zb\geq \theta_{\alpha}^{\rho}(\ub'\Zb)$ for any~$\ub$, that is, $\zb\in H^{\rho}_{\alpha,\ub}(P)$ for any~$\ub$. This implies that~$\zb\in R^\rho_{\alpha}(P)$, hence that~$M\!D^\rho(\zb,P)\geq \alpha$.
%=\inf_{\ub\in\mathcal{S}^{d-1}} G^\rho_\ub(\ub'z)$. 
By contradiction, assume now that~$M\!D^\rho(\zb,P)>\alpha$.  Then there exists~$\alpha'>\alpha$ such that~$\zb\in R^\rho_{\alpha'}(P)$. This means
\vspace{-.4mm}
  that~$\zb\in H^{\rho}_{\alpha',\ub}(P)$ for any~$\ub$, i.e., that~$\ub'\zb\geq \theta_{\alpha'}^{\rho}(\ub'\Zb)$ for any~$\ub$. Since~$G^\rho_\ub$ is right-continuous and montone non-decreasing, this entails that~$G^\rho_\ub(\ub'\zb)\geq G^\rho_\ub(\theta_{\alpha'}^{\rho}(\ub'\Zb))\geq \alpha'$ for any~$\ub$. Consequently, we must have~$\alpha=\inf_{\ub\in\mathcal{S}^{d-1}} G^\rho_\ub(\ub'\zb)\geq \alpha'$, a contradiction.   
\cqfd
\vspace{3mm}

The proof of Theorem~\ref{TheorMDepthConsistency} requires the following result.

 \begin{Lem}
 \label{LemMDepthConsistency}
Fix~$\rho\in\mathcal{C}$ and let~$\Zb$ be a random $d$-vector with distribution~$P\in\mathcal{P}_d^{\rho}$. Then, (i)
there exist~$c>0$ and~$\varepsilon>0$ such that
$
P[|\ub'(\Zb-\zb)|< c]\leq 1-\varepsilon
$
for any~$\ub\in\mathcal{S}^{d-1}$ and~$\zb\in\R^d$; 
(ii)
$\inf_{(\ub,\zb)\in\mathcal{S}^{d-1}\times \R^d} {\rm E}[|\psi_-(\ub'(\Zb-\zb))|]>0$. 
 \end{Lem}

{\sc Proof of Lemma~\ref{LemMDepthConsistency}.}
(i) First pick~$r$ so large that~$P[\|\Zb\|\geq r/2]\leq 1/2$. For any~$\ub\in\mathcal{S}^{d-1}$ and~$a>r$, we then have  
\begin{equation}
\label{kfhz1}
P[|\ub'\Zb-a|< r/2] 
%= P[ a-(r/2) < \ub'Z < a+(r/2)] 
%\leq P[ \ub'\Zb > a-(r/2)] 
\leq P[ \|\Zb\| > r/2] 
\leq 1/2
.
\end{equation}
It is thus sufficient to show that there exist~$c>0$ and~$\varepsilon>0$ such that $P[|\ub'\Zb-a|< c]\leq 1-\varepsilon$ for any~$\ub\in\mathcal{S}^{d-1}$ and~$a\in[0,r]$. 

By contradiction, assume that for any~$c>0$ and~$\varepsilon>0$, there exist~$\ub\in\mathcal{S}^{d-1}$ and~$a\in[0,r]$ such that~$P[|\ub'\Zb-a|< c]> 1-\varepsilon$. We can thus construct a sequence~$((\ub_n,a_n))$ in~$K=\mathcal{S}^{d-1}\times[0,r]$ such that
    $$
    P[|\ub_n'\Zb-a_n|< 1/n]> 1-(1/n)
    .
    $$
Compactness of~$K$ entails that there exists a subsequence~$((\ub_{n_\ell},a_{n_\ell}))$ that converges in~$K$, to~$(\ub_0,a_0)$ say. Clearly, we may assume that~$(\ub_{n_\ell}' \ub_0)$ is a monotone non-decreasing sequence and that~$(|a_{n_\ell}-a_0|)$ is a monotone non-increasing sequence (if that is not the case, one can always extract a further subsequence meeting these monotonicity properties). Let then~$I_\ell:=[a_0-|a_{n_\ell}-a_0|,a_0+|a_{n_\ell}-a_0|]$ and~$C_\ell:=\{ \ub\in\mathcal{S}^{d-1} : \ub'\ub_0 \geq \ub_{n_\ell}' \ub_0 \}$. Note that the sequences of sets~$(I_\ell)$ and~$(C_\ell)$ are monotone non-increasing with respect to inclusion, with~$\cap_\ell I_\ell=\{a_0\}$ and~$\cap_\ell C_\ell=\{\ub_0\}$. Therefore,
$$
	\lim_{\ell\to\infty} 
s_{\ell}
:=
\lim_{\ell\to\infty} 
P[ \Zb\in \cup_{a\in I_\ell}\cup_{\ub\in C_\ell} \{ \yb:|\ub'\yb-a|\leq 1/n_\ell\} ] 
=
P[ \ub_0' \Zb-a_0=0 ] 
.
$$
But, for any~$\ell$, 
$
s_{\ell}
\geq  
P[ |\ub_{n_\ell}'\Zb-a_{n_\ell}|\leq 1/n_\ell ]  
\geq 
1-(1/n_\ell)
$, which implies that~$(s_{\ell})$ converges to one as~$\ell$ diverges to infinity. Therefore, $P[ \ub_0' \Zb -a_0=0 ]=1$, which, since~$P\in\mathcal{P}_d^{\rho}$, is a contradiction.

(ii) Fix~$c$ and~$\varepsilon>0$ as in Part~(i) of the lemma, that is, such that, letting~$A_{\ub,\zb}:=\{\yb\in\R^d:|\ub'(\yb-\zb)|\geq c\}$, we have~$P[A_{\ub,\zb}]\geq \varepsilon$ for any~$\ub\in\mathcal{S}^{d-1}$ and~$\zb\in\R^d$. Then, for any such~$\ub$ and~$\zb$,
$$
{\rm E}[|\psi_-(\ub'(\Zb-\zb))|]
\geq
\int_{A_{\ub,\zb}} |\psi_-(\ub'(\yb-\zb))| \,dP(\yb)
\geq
\varepsilon \psi_-(c) >0
,
$$ 
since, for~$\rho\in\mathcal{C}$,~$\psi_-(t)$ may be zero at~$t=0$ only. The result follows.  
\cqfd
\vspace{3mm}

%%%%%%%%%%%%%%%%

{\sc Proof of Theorem~\ref{TheorMDepthConsistency}.} 
Let
$$
m_{\zb,\ub}(P)
=
m_{\zb,\ub}^\rho(P)
:=
\frac{{\rm E}[|\psi_-(\ub'(\Zb-\zb))| \mathbb{I}[\ub'(\Zb-\zb)\leq 0]]}{{\rm E}[|\psi_-(\ub'(\Zb-\zb))|]}
\cdot
$$
Fix~$\varepsilon>0$. For any~$Q\in\mathcal{P}_d^\rho$, let~$\ub_\zb(Q)$ be such that~$m_{\zb,\ub_\zb(Q)}(Q)\leq M\!D^\rho(\zb,Q)
\linebreak
+\varepsilon$ (existence follows from Theorem~\ref{reginter}). Then
\begin{eqnarray*}
\lefteqn{
| M\!D^\rho(\zb,P_n) - M\!D^\rho(\zb,P) | \, \mathbb{I}[ M\!D^\rho(\zb,P_n) \geq M\!D^\rho(\zb,P) ]
}
\\[2mm]
& &
\hspace{2mm}
\leq 
m_{\zb,\ub_\zb(P)}(P_n) - m_{\zb,\ub_\zb(P)}(P)+\varepsilon
\leq 
\sup_{\ub} | m_{\zb,\ub}(P_n) - m_{\zb,\ub}(P) | +\varepsilon
\end{eqnarray*}
and
\begin{eqnarray*}
\lefteqn{
| M\!D^\rho(\zb,P_n) - M\!D^\rho(\zb,P) | \, \mathbb{I}[ M\!D^\rho(\zb,P_n) < M\!D^\rho(\zb,P) ]
}
\\[2mm]
& &
\hspace{2mm}
\leq 
m_{\zb,\ub_\zb(P_n)}(P) - m_{\zb,\ub_\zb(P_n)}(P_n) +\varepsilon
\leq 
\sup_{\ub} | m_{\zb,\ub}(P) - m_{\zb,\ub}(P_n) | +\varepsilon
\end{eqnarray*}
(in this proof, all infima/suprema in~$u$ are over~$\mathcal{S}^{d-1}$, whereas those in~$z$ are over~$\R^d$). Therefore,
$$
| M\!D^\rho(\zb,P_n) - M\!D^\rho(\zb,P) |
\leq 
\sup_{\ub} | m_{\zb,\ub}(P_n) - m_{\zb,\ub}(P) | +\varepsilon
,
$$
and since this holds for any~$\varepsilon>0$, this yields
$$
| M\!D^\rho(\zb,P_n) - M\!D^\rho(\zb,P) |
\leq 
\sup_{\ub} | m_{\zb,\ub}(P_n) - m_{\zb,\ub}(P) | 
,
$$
which provides
$$
\sup_{\zb}
 | M\!D^\rho(\zb,P_n) - M\!D^\rho(\zb,P) |
\leq 
\sup_{\zb,\ub} | m_{\zb,\ub}(P_n) - m_{\zb,\ub}(P) | 
.
$$
Now, writing~$
q^-_{\zb,\ub}(P)
:=
{\rm E}[ |\psi_-(\ub'(\Zb-\zb))| \mathbb{I}[\ub'(\Zb-\zb)\leq 0]]
,
$
$
q^+_{\zb,\ub}(P)
:=
{\rm E}[ |\psi_-(\ub'(\Zb-\zb))| \mathbb{I}[\ub'(\Zb-\zb)> 0]]
$
and 
$q_{\zb,\ub}(P)
:=
q^-_{\zb,\ub}(P)+q^+_{\zb,\ub}(P)
=
{\rm E}[ |\psi_-(\ub'(\Zb-\zb))| ]
$,
we have 
\begin{eqnarray*}
\lefteqn{
\hspace{0mm} 
| m_{\zb,\ub}(P_n) - m_{\zb,\ub}(P) | 
=
\Big| 
\frac{q^-_{\zb,\ub}(P_n)}{q_{\zb,\ub}(P_n)}  - \frac{q^-_{\zb,\ub}(P)}{q_{\zb,\ub}(P)}  
\Big|
}
\\[5mm]
& & 
\hspace{3mm} 
\leq
\frac{|q^-_{\zb,\ub}(P_n)  - q^-_{\zb,\ub}(P)|}{q_{\zb,\ub}(P_n)}
+
q^-_{\zb,\ub}(P)
\Big|
\frac{1}{q_{\zb,\ub}(P_n)}  - \frac{1}{q_{\zb,\ub}(P)}  
\Big|
\\[5mm]
& & 
\hspace{3mm} 
\leq
\frac{|q^-_{\zb,\ub}(P_n)  - q^-_{\zb,\ub}(P)|+|q_{\zb,\ub}(P_n)-q_{\zb,\ub}(P)|}{q_{\zb,\ub}(P_n)}
%\leq
%\frac{2|q^-_{\zb,\ub}(P_n)  - q^-_{\zb,\ub}(P)|+|q^+_{\zb,\ub}(P_n)-q^+_{\zb,\ub}(P)|}{q_{\zb,\ub}(P_n)}
\\[4mm]
& & 
\hspace{3mm} 
\leq
\frac{2\sup_{\zb,\ub} |q^-_{\zb,\ub}(P_n)  - q^-_{\zb,\ub}(P)|+ \sup_{\zb,\ub} |q^+_{\zb,\ub}(P_n)-q^+_{\zb,\ub}(P)|}{\inf_{\zb,\ub} q_{\zb,\ub}(P) - \sup_{\zb,\ub} |q_{\zb,\ub}(P_n)-q_{\zb,\ub}(P)|}
\end{eqnarray*}
for any~$z\in\R^d$ and~$u\in\mathcal{S}^{d-1}$. Since~$\inf_{\zb,\ub} q_{\zb,\ub}(P)>0$ (Lemma~\ref{LemMDepthConsistency}), it only remains to prove that 
\begin{equation}
\label{GC2a} 
\sup_{\zb,\ub} |q^-_{\zb,\ub}(P_n)  - q^-_{\zb,\ub}(P)|
\stackrel{a.s.}{\to} 0 
\quad\textrm{ and }\quad
\sup_{\zb,\ub} |q^+_{\zb,\ub}(P_n)  - q^+_{\zb,\ub}(P)|
\stackrel{a.s.}{\to} 0 
\end{equation}
as~$n\to\infty$. Let us focus on the first convergence in~(\ref{GC2a}). Clearly, we are after a Glivenko-Cantelli theorem for the classes of functions
$$
\mathcal{F}^{\rho}
:=
\Big\{
\yb\mapsto f^{\rho}_{\zb,\ub}(\yb) := -\psi_-(\ub'(\yb-\zb)) \mathbb{I}[ \ub'(\yb-\zb)\leq 0 ]
:
\zb\in\R^d
,
\ub\in\mathcal{S}^{d-1}
\Big\}
.
$$
The collection~$\mathcal{H}$ of all halfspaces in~$\R^{d+1}$ is a Vapnik-Chervonenkis class; see, e.g., Page~152 of~\cite{WelVan1996}. %Hence, Lemma~2.6.17 of the same implies that~$\mathcal{H}\sqcap \mathcal{H}:=\{ H_1 \cap H_2: H_1,H_2\in \mathcal{H}\}$ is also a Vapnik-Chervonenkis class. 
Consequently, defining the subgraph of a function~$f:\R^d\to\R$ as~$s_f:=\{(\yb,t)\in\R^{d+1}: t<f(\yb)\}$ and letting
$$
\mathcal{F}
:= 
\Big\{
\yb\mapsto f_{\zb,\ub}(y) := \ub'(\yb-\zb) %\mathbb{I}[ \ub'(\yb-\zb)\leq 0 ]
:
\zb\in\R^d
,
\ub\in\mathcal{S}^{d-1}
\Big\}
,
$$
the collection of subgraphs~$\{s_f:f\in\mathcal{F}\}$, as a subset of~$\mathcal{H}$, is a Vapnik-Chervonenkis class. In other words,~$\mathcal{F}$ is a VC- subgraph class (see, e.g., Section~2.6.2 of~\citealp{WelVan1996}). Now, since~$t\mapsto -\psi_-(t)\mathbb{I}[t\leq 0]$ is a monotone function, Lemma~2.6.18(viii) of~\cite{WelVan1996} implies that~$\mathcal{F}^{\rho}$ is itself a VC-subgraph class, hence is Glivenko-Cantelli, which implies the first convergence in~(\ref{GC2a}). Since the same reasoning establishes the second convergence in~(\ref{GC2a}), the result is proved.   
\cqfd
\vspace{3mm}

{\sc Proof of Theorem~\ref{TheorMmedianConsistency}.} 
The mapping~$\zb\mapsto M\!D^\rho(\zb,P)$ is upper semicontinuous (see the paragraph below Theorem~\ref{TheorRegionsOrig}) and constant over~$R_{\alpha_*}^\rho(P)$, with~$\alpha_*:=\max_{\zb\in\R^d} M\!D^{\rho}(\zb,P)$. Clearly, it is easy to define a mapping~$\zb\mapsto M\!D_0^\rho(\zb,P)$ that is upper semicontinuous, agrees with~$\zb\mapsto M\!D^\rho(\zb,P)$ in the complement of~$R_{\alpha_*}^\rho(P)$, and for which~$\zb_*(P)$ is the unique maximizer. The result then follows from Theorem~2.12 and Lemma~14.3 in \cite{Kos2008} and from Theorem~\ref{TheorMDepthConsistency}.
% that $z_*(P_n) \to z_*(P_n)$ almost surely as~$n\to\infty$.
\cqfd
\vspace{3mm}

The proof of Theorem~\ref{TheorMRegionsConsistency} requires Lemmas~\ref{LemGandMDepthContinuous}--\ref{LemPourMDepthStrictMonotonicity} below.

\begin{Lem}
\label{LemGandMDepthContinuous}
Fix~$\rho\in\mathcal{C}$ and assume that~$P\in\mathcal{P}_d^{\rho}$ assigns probability zero to all hyperplanes in~$\R^d$. Then, (i) 
$
(\ub,\zb)
\mapsto 
G_\ub^\rho(\ub'\zb)
%:=
%\frac{{\rm E}[|\psi_-(\ub'(\Zb-\zb))| \mathbb{I}[\ub'(\Zb-\zb)\leq 0]]}{{\rm E}[|\psi_-(\ub'(\Zb-\zb))|]}
$
is continuous over~$\mathcal{S}^{d-1}\times \R^d$; (ii) $\zb\mapsto M\!D^\rho(\zb,P)$ is uniformly continuous on~$\R^d$. 
\end{Lem}

{\sc Proof of Lemma~\ref{LemGandMDepthContinuous}.}
(i) Letting~$\ell_\pm(\ub,\zb)
:=
{\rm E}[|\psi_-(\ub'(\Zb-\zb))| \mathbb{I}[\pm \ub'(\Zb-\zb)> 0]]$, first note that, under the assumptions considered on~$P$, we have
$$
G_\ub^\rho(\ub'\zb)
=
\frac{\ell_-(\ub,\zb)}{\ell_-(\ub,\zb)+\ell_+(\ub,\zb)}
=:
\ell(\ub,\zb)
\cdot 
$$
We only show that~$(\ub,\zb)\mapsto \ell_-(\ub,\zb)$ is continuous over~$\mathcal{S}^{d-1}\times \R^d$ (continuity of~$\ell_+$ can indeed be proved along the exact same lines). To that end, fix~$(\ub_0,\zb_0)\in\mathcal{S}^{d-1}\times\R^d$ and let~$((\ub_n,\zb_n))$ be a sequence in~$\mathcal{S}^{d-1}\times\R^d$ that converges to~$(\ub_0,\zb_0)$. Let~$\mathcal{D}$ be the collection of discontinuity points of~$\psi_-$ (recall that~$\mathcal{D}$ is at most countable). For any~$\zb\in\R^d$ such that~$\ub_0'(\zb-\zb_0)\notin \mathcal{D}$, 
\begin{eqnarray*}
	\lefteqn{
h_{n}(\zb)
:=
|\psi_-(\ub_n'(\zb-\zb_n))| \mathbb{I}[\ub_n'(\zb-\zb_n)< 0]
}
\\[2mm]
& & 
\hspace{23mm} 
-
|\psi_-(\ub_0'(\zb-\zb_0))| \mathbb{I}[\ub_0'(\zb-\zb_0)< 0]
\to
0
\end{eqnarray*}
as~$n\to \infty$. Since~$P$ assigns probability mass zero to all hyperplanes, this convergence holds $P$-almost everywhere. Assuming that~$h_{n}(\zb)$ is upper-bounded by a $P$-integrable function that does not depend on~$n$, the Lebesgue Dominated Convergence Theorem then yields that
$$
\ell_-(\ub_n,\zb_n)-\ell_-(\ub_0,\zb_0)
=
\int_{\R^d}
h_n(\zb)
\,dP(\zb)
\to 0
$$
as~$n\to\infty$, which proves continuity of~$\ell_-$. 

It therefore remains to show that~$h_{n}(\zb)$ is upper-bounded by a $P$-integrable function that does not depend on~$n$. To do so, let~$v_1=e_1,\ldots,v_d=e_d,v_{d+1}=-e_1,\ldots,v_{2d}=-e_{d}$, where~$e_j$ denotes the $j$th vector of the canonical basis of~$\R^d$. With this notation, note that, for any~$u\in\mathcal{S}^{d-1}$ and~$z\in\R^d$, one has
\begin{eqnarray*}
u'z 
=
\sum_{j=1}^d
u_j z_j
&\!\!\!\!\geq \!\!\!\!& 
\bigg(\sum_{j=1}^d |u_j|\bigg) 
\min\!\big\{\!-|z_j|,\ j=1,\ldots,d\big\}
\\[2mm]
&\!\!\!\!\geq \!\!\!\!& 
\min\!\big\{v_j'z,\ j=1,\ldots,2d\big\}
.
\end{eqnarray*}
For large~$n$, we have~$|\ub_n'\zb_n|\leq \|z_0\|+1=:c$, which then yields
\begin{eqnarray*}
\lefteqn{
|\psi_-(\ub_n'(\zb-\zb_n))| \mathbb{I}[\ub_n'(\zb-\zb_n)< 0]
\leq
|\psi_-(\ub_n'\zb-c)| \mathbb{I}[\ub_n'\zb-c< 0]
}
\\[2mm]
& & 
\hspace{3mm} 
\leq
|\psi_-(\min\!\big\{v_j'z,\ j=1,\ldots,2d\big\}-c)| 
\leq
\sum_{j=1}^{2d}
|\psi_-(v_j'z-c)|
. 
\end{eqnarray*}
We thus conclude that, still for large~$n$, the function~$z\mapsto h_{n}(\zb)$ is upper-bounded by the function
$$
z
\mapsto
|\psi_-(\ub_0'\zb-\ub_0'\zb_0)|
+
\sum_{j=1}^{2d} 
|\psi_-(v_j'z-c)|
,
$$
which is $P$-integrable and does not depend on~$n$. This completes the proof of Part~(i) of the result. 

(ii) Fix then an arbitrary compact set~$K\subset \R^d$ and~$\varepsilon>0$. Since~$\ell$ is continuous over the compact set~$\mathcal{S}^{d-1}\times K$, it is also uniformly continuous on that set. Hence, there exists~$\delta>0$ such that for any~$\ub_1,\ub_2\in\mathcal{S}^{d-1}$ and~$\zb_1,\zb_2\in K$ satisfying $\max(\|\ub_1-\ub_2\|,\|\zb_1-\zb_2\|)<\delta$, we have $|\ell(\ub_1,\zb_1)-\ell(\ub_2,\zb_2)|<\varepsilon$. For any~$\zb\in\R^d$, pick arbitrarily~$u_z\in\mathcal{S}^{d-1}$ such that~$M\!D^\rho(\zb,P)=\ell(\ub_\zb,\zb)$; existence follows by using Part~(i) of the lemma and the compactness of~$\mathcal{S}^{d-1}$ in Theorem~\ref{reginter}. Then, for any~$\zb_1,\zb_2\in K$ with~$\|\zb_1-\zb_2\|<\delta$, we have
$$
M\!D^\rho(\zb_1,P)
=
\ell(\ub_{\zb_1},\zb_1)
>
\ell(\ub_{\zb_1},\zb_2)-\varepsilon
\geq
M\!D^\rho(\zb_2,P)-\varepsilon.
$$  
By symmetry, we also have~$M\!D^\rho(\zb_2,P)
>
M\!D^\rho(\zb_1,P)-\varepsilon$, which yields $|M\!D^\rho(\zb_2,P)-M\!D^\rho(\zb_1,P)|<\varepsilon$. Consequently, $\zb\mapsto M\!D^\rho(\zb,P)$ is uniformly continuous over~$K$. 

We now show that uniform continuity extends to~$\R^d$. To do so, fix~$\varepsilon>0$ and pick~$C$ large enough to have~$M\!D^\rho(\zb,P)<\varepsilon/2$ as soon as~$\zb\notin B_0(C):=\{\zb\in\R^d:\|\zb\|\leq C\}$ (existence of~$C$ follows from Theorem~\ref{Theordepthaxioms}(iv)). Since $\zb\mapsto M\!D^\rho(\zb,P)$ is uniformly continuous over any~$B_0(r)$, there exists~$\delta>0$ such that for any~$\zb_1,\zb_2\in B_0(C+1)$ such that~$\|\zb_1-\zb_2\|<\delta$, one has~$|M\!D^\rho(\zb_2,P)-M\!D^\rho(\zb_1,P)|<\varepsilon$. Letting~$\tilde{\delta}:=\min(\delta,1)$, it is then easy to check that for any~$\zb_1,\zb_2\in\R^d$ such that~$\|\zb_1-\zb_2\|<\tilde{\delta}$, we must have~$|M\!D^\rho(\zb_2,P)-M\!D^\rho(\zb_1,P)|<\varepsilon$ (note that as soon as one of such~$\zb_1,\zb_2$ belongs to~$B_0(C)$, then they both belong to~$B_0(C+1)$).   
\cqfd
\vspace{3mm}

%%%%%%%%%%%%%%%%%

\begin{Lem}
\label{LemPourMDepthStrictMonotonicity}
Fix~$\rho\in\mathcal{C}$ and assume that~$P\in\mathcal{P}_d^{\rho}$ assigns probability zero to all hyperplanes in~$\R^d$ and has a connected support (in the sense defined above Theorem~\ref{TheorMRegionsConsistency}). Fix~$\zb_0,\zb_1\in \R^d$ with $M\!D^{\rho}(\zb_0,P)=\max_{z\in\R^d}M\!D^{\rho}(\zb,P)
\linebreak
>M\!D^{\rho}(\zb_1,P)>0$. Then,
\begin{equation}
\label{qc}
M\!D^{\rho}((1-\lambda)\zb_0+\lambda \zb_1,P)
> 
M\!D^{\rho}(\zb_1,P)
\end{equation}
for any~$\lambda\in(0,1)$. 
\end{Lem}

{\sc Proof of Lemma~\ref{LemPourMDepthStrictMonotonicity}.}
Write~$\zb_\lambda:=(1-\lambda)\zb_0+\lambda \zb_1$ and let~$\ub_\lambda\in\mathcal{S}^{d-1}$ be such that~$M\!D^{\rho}(\zb_\lambda,P)=G_{\ub_\lambda}^\rho(\ub_\lambda'\zb_\lambda)$; existence follows from Lemma~\ref{LemGandMDepthContinuous}(i) and the compactness of~$\mathcal{S}^{d-1}$. 
Assume first that~$\ub_\lambda'\zb_0=\ub_\lambda'\zb_1(=\ub_\lambda'\zb_\lambda)$. Then, since~$z_0$ has maximal depth, we have 
$$
M\!D^{\rho}(\zb_\lambda,P)
%=G_{\ub_\lambda}^\rho(\ub_\lambda'\zb_\lambda)
=G_{\ub_\lambda}^\rho(\ub_\lambda'\zb_0)
= 
M\!D^{\rho}(\zb_0,P)
>
M\!D^{\rho}(\zb_1,P)
.
$$
We can thus restrict to the case~$\ub_\lambda'\zb_0\neq \ub_\lambda'\zb_1$, under which~$\theta_b:=\ub_\lambda'\zb_\lambda=(1-\lambda) \ub_\lambda'\zb_0+\lambda \ub_\lambda'\zb_1> \min(\ub_\lambda'\zb_0,\ub_\lambda'\zb_1)=:\theta_a$. Now, 
$$
G_{\ub_\lambda}^\rho(\theta_a)=\min(G_{\ub_\lambda}^\rho(\ub_\lambda'\zb_0),G_{\ub_\lambda}^\rho(\ub_\lambda'\zb_1))
\geq \min(M\!D^{\rho}(\zb_0,P),M\!D^{\rho}(\zb_1,P))>0
.
$$    
Writing~$F(x):=P[\ub_\lambda'\Zb\leq x]$, this implies that~$F(\theta_a)>0$ (indeed, $F(\theta_a)=0$ would yield~$H^{\rho}_2(\theta)+\delta F(\theta_a)=0$, hence~$G_{\ub_\lambda}^\rho(\theta_a)=0$; see~(\ref{expressionG}) based on~$P=P_{\ub_\lambda}$, the distribution of~$\ub_\lambda'\Zb$). 
Quasi-concavity of~$M\!D^{\rho}(\,\cdot\,,P)$ ensures that
$$
G_{-\ub_\lambda}^\rho(-\ub_\lambda'\zb_\lambda)\geq M\!D^{\rho}(\zb_\lambda,P)\geq \min(M\!D^{\rho}(\zb_0,P),M\!D^{\rho}(\zb_1,P))>0
,
$$
which similarly implies that~$P[-\ub_\lambda'\Zb\leq -\ub_\lambda'\zb_\lambda]>0$, that is, $1-F(\theta_b-0)>0$, where~$F(\theta_b-0)$ still denotes the limit of~$F(\theta)$ when~$\theta$ converges to~$\theta_b$ from below. The assumption of connected support therefore implies that~$F(\theta_b-0)-F(\theta_a)>0$. Lemma~\ref{LemPourLemMDepthStrictMonotonicity} thus entails that
$$
M\!D^{\rho}(\zb_\lambda,P)
=G_{\ub_\lambda}^\rho(\theta_b)
>G_{\ub_\lambda}^\rho(\theta_a)
\geq \min(M\!D^{\rho}(\zb_0,P),M\!D^{\rho}(\zb_1,P))
,
$$ 
which establishes the result.
\cqfd
\vspace{3mm}

%%%%%%%%%%%%%%%%%

We can now prove Theorem~\ref{TheorMRegionsConsistency}. 
\vspace{2mm}

{\sc Proof of Theorem~\ref{TheorMRegionsConsistency}.} 
In view of Theorem~\ref{TheorMDepthConsistency} and the result in Theorem~4.5 from \cite{Dyc2016} (more precisely, its corollary in a random sampling scheme as discussed in page~13 of that paper), it is sufficient to prove that~$M\!D^{\rho}(\,\cdot\,,P)$ is \emph{strictly monotone}, in the sense that, for any~$\alpha\in(0,\alpha_*)$, with~$\alpha_*:=\max_{\yb\in\R^d} M\!D^{\rho}(\yb,P)$, the region~$R^\rho_\alpha(P)$ is the closure~$\bar{R}^\rho_{\alpha,>}(P)$ of~$R^\rho_{\alpha,>}(P):=\{\zb\in\R^d:M\!D^{\rho}(\zb,P)>\alpha\}$. 

Now, since~$R^\rho_\alpha(P)$ is closed and contains~$R^\rho_{\alpha,>}(P)$, we have that~$\bar{R}^\rho_{\alpha,>}(P)\subset R^\rho_\alpha(P)$. To show that~$R^\rho_\alpha(P)\subset \bar{R}^\rho_{\alpha,>}(P)$, fix~$\zb\in R^\rho_\alpha(P)$. If~$M\!D^{\rho}(\zb,P)>\alpha$, then~$\zb$ trivially belongs to~$\bar{R}^\rho_{\alpha,>}(P)$, so that we may assume that~$M\!D^{\rho}(\zb,P)=\alpha$. Consider then the line segment associated with~$\zb_\lambda:=(1-\lambda)\zb_*+\lambda \zb$, $\lambda\in(0,1)$, from an arbitrary deepest point~$\zb_*$ of~$M\!D^{\rho}(\,\cdot\,,P)$ to~$\zb$. Under the assumptions considered, Lemma~\ref{LemPourMDepthStrictMonotonicity} guarantees that~$(\zb_{1-(1/n)})$ is a sequence in~$R^\rho_{\alpha,>}(P)$ that converges to~$\zb$, so that~$\zb\in\bar{R}^\rho_{\alpha,>}(P)$. We conclude that we also have~$R^\rho_\alpha(P)\subset \bar{R}^\rho_{\alpha,>}(P)$, hence that~$M\!D^{\rho}(\,\cdot\,,P)$ is strictly monotone. This establishes the result.  
\cqfd
\vspace{3mm}

%%%%%%%%%%

%%%%%%%%%%%%%%%%%%%%%%%%%%

%\subsection{Proofs of Section~\ref{secExpectiledepth}}
%\label{secAppSec5}

{\sc Proof of Theorem~\ref{TheorExpectDeepestPoint}}. 
Let~$Y$ be a random variable with a distribution in~$\mathcal{P}_1$. It follows from Theorem~\ref{Theorlemjones} and Lemma~\ref{LemPourMonotoneExpect} that the order-$1/2$ expectile of~$Y$ is~${\rm E}[Y]$ and that this expectile is the only value of~$\theta$ such that 
$
{\rm E}[|Y-\theta| \mathbb{I}[Y-\theta\leq 0]]/{\rm E}[|Y-\theta|]
=
1/2
$.
For any~$\ub\in\mathcal{S}^{d-1}$, we thus have that
$$
\frac{{\rm E}[ |\ub'(\Zb-{\rm E}[\Zb])| \mathbb{I}[\ub'(\Zb-{\rm E}[\Zb])\leq 0 ]]}{{\rm E}[|\ub'(\Zb-{\rm E}[\Zb])|]}
=
\frac{1}{2}
$$
($\Zb$ is a random $d$-vector with distribution~$P$), so that, writing~$\mub_P={\rm E}[\Zb]$,~(\ref{DefinExpectileDepth2}) yields~$E\!D(\mub_P,P)=1/2$. Assume now that there exists~$\zb\in\R^d$ such that $E\!D(\zb,P)=e>1/2$. Then, for an arbitrary fixed~$\ub\in\mathcal{S}^{d-1}$, one has
\begin{equation}
\label{ineq1}
\frac{{\rm E}[ |\ub'(\Zb-\zb)| \mathbb{I}[\ub'(\Zb-\zb)\leq 0 ]]}{{\rm E}[|\ub'(\Zb-\zb)|]}
\geq e
\end{equation}
and
\begin{equation}
\label{ineq2}
\frac{{\rm E}[ |-\ub'(\Zb-\zb)| \mathbb{I}[-\ub'(\Zb-\zb)\leq 0 ]]}{{\rm E}[|-\ub'(\Zb-\zb)|]}
\geq e.
\end{equation}
Adding up these two inequalities yields
%$$
%\frac{{\rm E}[ |\ub'(\Zb-\zb)| (1+\mathbb{I}[\ub'(\Zb-\zb)=0 ])]}{{\rm E}[|\ub'(\Zb-\zb)|]}
%\geq 2e,
%$$
%that is, 
$1\geq 2e$, a contradiction. We conclude that $E\!D(\mub_P,P)\geq E\!D(\zb,P)$ for any~$\zb\in\R^d$, and it only remains to show that~$\mub_P$ is the only maximizer of expectile depth. For that purpose, assume that~$\zb$ is such that~$E\!D(\zb,P)=1/2$. Then for any~$\ub\in\mathcal{S}^{d-1}$, the inequalities in~(\ref{ineq1})--(\ref{ineq2}) hold with~$e=1/2$ and are actually equalities (indeed, would there be a direction~$\ub$ for which at least one of these inequalities would be strict, then adding up both inequalities as above would provide~$1>1$, a contradiction). Thus, for any~$\ub\in\mathcal{S}^{d-1}$, $\ub'\zb$ is the order-$1/2$ expectile of~$\ub'\Zb$, that is, $\ub'\zb={\rm E}[\ub'\Zb]$ (see above). This means that~$\ub'(\zb-\mub_P)=0$ for any~$\ub\in\mathcal{S}^{d-1}$, which shows that~$\zb=\mub_P$. 
\cqfd  
\vspace{3mm}

%%%%%%%%%%%%%%%%%%

The proof of Theorem~\ref{TheorMonotExpectDepth} requires the following \emph{strict} quasi-concavity property.

\begin{Lem}
\label{Lemstrictquasiconc}
Let~$P$ be a probability measure in~$\mathcal{P}_d$ and denote as~$\mub_P$ the corresponding mean vector.
Then, 
\begin{equation}
\label{qc}
E\!D((1-\lambda)\mub_P+\lambda \zb,P)
> 
E\!D(\zb,P)
\end{equation}
for any~$\lambda\in[0,1)$ and~$z(\neq \mub_P)$ in the c-support~$C_P$ of~$P$. 
\end{Lem}

{\sc Proof of Lemma~\ref{Lemstrictquasiconc}}. 
Fix~$\zb_\lambda:=(1-\lambda)\mub_P+\lambda \zb$, with~$\zb(\neq \mub_P)\in C_P$ and~$\lambda\in(0,1)$ (for~$\lambda=0$, the result directly follows from Theorem~\ref{TheorExpectDeepestPoint}). Let~$A:=\{u\in\mathcal{S}^{d-1}: u'(z-\mub_P)=0 \}$. First note that the proof of Theorem~\ref{TheorExpectDeepestPoint} entails that, for any~$u\in A$, 
$$
G_{\ub}^\rho(\ub'\zb_\lambda)
=
G_{\ub}^\rho(\ub'\mub_P)
=
\frac{1}{2}
=
E\!D(\mub_P,P)
>
E\!D(\zb_\lambda,P)
$$
(throughout the proof, $\rho(t)=t^2$ is the quadratic loss function), so that~(\ref{DefinExpectileDepth2}) yields
\begin{equation}
\label{toprostqu}
E\!D(\zb_\lambda,P)
=
\min_{u\in\mathcal{S}^{d-1}\setminus A} G_{\ub}^\rho(\ub'\zb_\lambda)
.
\end{equation}
Fix then~$u\in \mathcal{S}^{d-1}\setminus A$. Since both~$z$ and~$\mub_P$ belong to~$C_P$, we have~$P[\ub'\Zb\leq\min(\ub'\mub_P,\ub'\zb)]>0$ and~$P[\ub'\Zb\geq\max(\ub'\mub_P,\ub'\zb)]>0$. Recalling that~$u\notin A$, we have $\ub'\zb_\lambda \in(\min(\ub'\mub_P,\ub'\zb),\max(\ub'\mub_P,\ub'\zb))$. Lemma~\ref{LemPourMonotoneExpect} then yields
\begin{eqnarray*}
\lefteqn{
G_{\ub}^\rho(\ub'\zb_\lambda) 
> 
G_{\ub}^\rho(\min(\ub'\mub_P,\ub'\zb)) 
= 
\min(G_{\ub}^\rho(\ub'\mub_P),
G_{\ub}^\rho(\ub'\zb))
}
\\[2mm]
& & 
\hspace{16mm} 
\geq 
\min(E\!D(\mub_P,P),E\!D(\zb,P))
=
E\!D(\zb,P)
\end{eqnarray*}
for any~$\ub\in\mathcal{S}^{d-1}\setminus A$. The result thus follows from~(\ref{toprostqu}).  
\cqfd
\vspace{3mm}

%%%%%%%%%%

{\sc Proof of Theorem~\ref{TheorMonotExpectDepth}.}
Fix~$0\leq r_1<r_2<r_\ub(P)$. Then, $\mub_P+r \ub\in C_P$ for any~$r\in [0,r_2]$. Therefore, Lemma~\ref{Lemstrictquasiconc} yields~$E\!D(\mub_P+r_1 \ub,P)>E\!D^{\rho}(\mub_P+r_2 \ub,P)$, so that~$r\mapsto E\!D(\mub_P+r \ub,P)$ is monotone strictly decreasing in~$[0,r_\ub(P))$. Since Theorem~\ref{Theordepthaxioms}(iii) ensures that~$E\!D(\mub_P+r \ub,P)=0$ for any~$r>r_\ub(P)$, the result follows from Theorem~\ref{TheorexpectileDepthContinuous} (whose proof is independent of Theorem~\ref{TheorMonotExpectDepth}). 
\cqfd
\vspace{3mm}
 
%%%%%%%%%%%%%%%%%%%

The proof of Theorem~\ref{TheorexpectileDepthContinuous} requires the following preliminary result.

\begin{Lem}
\label{LemExpectileDepthSmoothness}
Let~$\Zb$ be a random $d$-vector with distribution~$P\in\mathcal{P}_d$. 
Then, (i) 
$
(\ub,\zb) \mapsto 
%\min_{\ub\in\mathcal{S}^{d-1}} 
\ell_-(\ub,\zb):={\rm E}[ |\ub'(\Zb-\zb)| \mathbb{I}[\ub'(\Zb-\zb)<0]]
$
and
$
(\ub,\zb) \mapsto 
%\min_{\ub\in\mathcal{S}^{d-1}} 
\ell_+(\ub,\zb):={\rm E}[ |\ub'(\Zb-\zb)| \mathbb{I}[\ub'(\Zb-\zb)> 0]]
$
are continuous over~$\mathcal{S}^{d-1}\times \R^d$;
(ii) for any~$\ub\in\mathcal{S}^{d-1}$, the functions~$\zb \mapsto \ell_\pm(\ub,\zb)$ admit, at any~$\zb\in\R^d$, directional derivatives in any direction;
(iii) if, moreover, $P$ is smooth in a neighbourhood of~$\zb_0$ (in the sense defined in Theorem~\ref{TheorexpectileDepthContinuous}), then,  for any~$\ub\in\mathcal{S}^{d-1}$, the functions~$\zb \mapsto \ell_\pm(\ub,\zb)$
%$
%(\ub,\zb) \mapsto 
%\min_{\ub\in\mathcal{S}^{d-1}} 
%\ell_-(\ub,\zb)
%:={\rm E}[ |\ub'(\Zb-\zb)| \mathbb{I}[\ub'(\Zb-\zb)<0]]
%$
%and
%$
%(\ub,\zb) \mapsto 
%%\min_{\ub\in\mathcal{S}^{d-1}} 
%\ell_+(\ub,\zb):={\rm E}[ |\ub'(\Zb-\zb)| \mathbb{I}[\ub'(\Zb-\zb)> 0]]
%$
are continuously differentiable in a neighbourhood of~$\zb_0$.
\end{Lem}

{\sc Proof of Lemma~\ref{LemExpectileDepthSmoothness}}. 
(i) We only prove the result for~$\ell_-$, as the proof for~$\ell_+$ is entirely similar. 
Fix~$(\ub_0,\zb_0)\in\mathcal{S}^{d-1}\times \R^d$ and write~$B_{\zb_0}(r):=\{ \zb\in\R^d:\|\zb-\zb_0\|\leq r\}$. For any~$\yb\in\R^d$, we have that 
$
(\ub,\zb) \mapsto \ub'(\zb-\yb) \mathbb{I}[\ub' \yb< \ub'\zb]
$
is continuous at~$(\ub_0,\zb_0)$. Moreover, for any~$(\ub,\zb)\in\mathcal{S}^{d-1}\times B_{\zb_0}(1)$,
the function~$\yb\mapsto \ub'(\zb-\yb) \mathbb{I}[\ub' \yb< \ub'\zb]$ is upper-bounded by the function~$\yb\mapsto
\|\zb_0\|+1+\|\yb\|
$ that is~$P$-integrable and does not depend on~$(\ub,\zb)$. The Lebesgue Dominated Convergence Theorem therefore yields the result. 

(ii) %For the same reason as above, we only prove the result for~$\ell_-$. Fix~$u,v\in\mathcal{S}^{d-1}$ and~$z_0\in\R^d$. 
We will show that 
\begin{equation}
\label{partialdetlmoins}
\frac{\partial\ell_-}{\partial \vb}(z_0)
=
(\ub'\vb) {\rm P}[\ub'\Zb<\ub'\zb_0] \mathbb{I}[\ub'\vb<0]
+
(\ub'\vb) {\rm P}[\ub'\Zb\leq \ub'\zb_0] \mathbb{I}[\ub'\vb>0]
.
\end{equation}
To do so, note that, for any~$h>0$, 
\begin{eqnarray*}
\lefteqn{
m_{\zb_0,\ub,\vb}(h,\yb)
}
\\[2mm]
& & 
\hspace{1mm}
:=
\frac{1}{h}
\big\{
\ub'(\zb_0+h\vb-\yb)\mathbb{I}[\ub'\yb<\ub'(\zb_0+h\vb)] - \ub'(\zb_0-\yb)\mathbb{I}[\ub'\yb<\ub'\zb_0]
\big\}
\\[2mm]
& & 
\hspace{12mm}
-
\big\{
\ub'\vb\mathbb{I}[\ub'\yb< \ub'\zb_0] \mathbb{I}[\ub'\vb<0]
+
\ub'\vb\mathbb{I}[\ub'\yb\leq \ub'\zb_0] \mathbb{I}[\ub'\vb>0]
\big\}
\\[2mm]
& & 
\hspace{1mm} 
=
\frac{1}{h}\, \big( \ub'(\zb_0+h\vb) - \ub'\yb \big)
S(\ub'\vb)
\mathbb{I}[\ub'\yb \in \mathcal{I}_{\zb_0,\ub,\vb}(h)]
,
\end{eqnarray*}
where the sign function~$S$ was defined on page~\pageref{pagedefsign} and where~$\mathcal{I}_{\zb_0,\ub,\vb}(h)$ denotes the open interval with endpoints~$\ub'\zb_0$ and~$\ub'(\zb_0+h\vb)$. This shows that, for any~$\yb\in\R^d$,~$m_{\zb_0,\ub,\vb}(h,\yb)$ converges to zero as~$h$ goes to zero from above and that the function~$\yb\mapsto |m_{\zb_0,\ub,\vb}(h,\yb)|$ is upper-bounded by the function~$\yb\mapsto |\ub'\vb|$ that is $P$-integrable and does not depend on~$h$. Consequently, the Lebesgue Dominated Convergence Theorem entails that, as~$h$ goes to zero from above, 
\begin{eqnarray*}
\lefteqn{
\frac{\ell_-(\ub,\zb_0+h\vb)-\ell_-(\ub,\zb_0)}{h} - 
\Big\{
(\ub'\vb) {\rm P}[\ub'\Zb<\ub'\zb_0] \mathbb{I}[\ub'\vb<0]
}
\\[2mm]
 & & 
\hspace{13mm}
+
(\ub'\vb) {\rm P}[\ub'\Zb\leq \ub'\zb_0] \mathbb{I}[\ub'\vb>0]
\Big\}
=
\int_{\R^d} m_{\zb_0,\ub,\vb}(h,\yb) \,dP(\yb)
\to 0,
\end{eqnarray*}
which establishes~(\ref{partialdetlmoins}). The exact same reasoning allows to show that 
\begin{equation}
\label{partialdetlplus}
\frac{\partial\ell_+}{\partial \vb}(\zb_0)
=
- (\ub'\vb) {\rm P}[\ub'\Zb\geq \ub'\zb_0] \mathbb{I}[\ub'\vb<0]
- (\ub'\vb) {\rm P}[\ub'\Zb> \ub'\zb_0] \mathbb{I}[\ub'\vb>0]
.
\end{equation}

(iii) It trivially follows from the Lebesgue Dominated Convergence Theorem that, under the smoothness assumption considered, the functions~$\zb\mapsto \frac{\partial\ell_\pm}{\partial \vb}(\zb)$ in~(\ref{partialdetlmoins})-(\ref{partialdetlplus}) are continuous in a neighborhood of~$\zb_0$.
\cqfd
\vspace{3mm}
 
%%%%%%%%%

{\sc Proof of Theorem~\ref{TheorexpectileDepthContinuous}}.
Throughout the proof, $\Zb$ denotes a random $d$-vector with distribution~$P$. (i) It directly follows from Lemma~\ref{LemExpectileDepthSmoothness}(i) that 
\begin{equation}
\label{ellfunction}
(\ub,\zb) \mapsto 
%\min_{\ub\in\mathcal{S}^{d-1}} 
\ell(\ub,\zb):=\frac{{\rm E}[ |\ub'(\Zb-\zb)| \mathbb{I}[\ub'\Zb\leq \ub'\zb]]}{{\rm E}[|\ub'(\Zb-\zb)|]}
=\frac{\ell_-(\ub,\zb)}{\ell_-(\ub,\zb)+\ell_+(\ub,\zb)}
\end{equation}
($\ell_-$ and~$\ell_+$ refer to the same functions as in Lemma~\ref{LemExpectileDepthSmoothness}) is continuous over~$\mathcal{S}^{d-1}\times \R^d$. The proof of uniform continuity then proceeds exactly as in the proof of Lemma~\ref{LemGandMDepthContinuous}. 
(ii) Fix~$\zb_0\in\R^d$ and~$\ub,\vb\in\mathcal{S}^{d-1}$.  Lemma~\ref{LemExpectileDepthSmoothness}(ii) implies that~$\zb\mapsto \ell(\ub,\zb)$ admits a directional derivative at~$\zb_0$ in direction~$\vb$. In dimension~$d=1$, there are finitely many~$\ub$'s that are to be considered in~(\ref{DefinExpectileDepth2}), so that the aforementioned differentiability readily entails equidifferentiability in the sense of \cite{MilSeg2002}. The result then follows from Theorem~3 of \cite{MilSeg2002}.
(iii) By assumption, $P$ is smooth over a neighbourhood~$\mathcal{N}$ of~$\zb_0$. Lemma~\ref{LemExpectileDepthSmoothness}(iii) then yields that, for any~$\ub\in\mathcal{S}^{d-1}$, $\zb\mapsto \ell(\ub,\zb)$ is continuously differentiable over~$\mathcal{N}$. The result then follows from Theorem~1 in \cite{Dan1966} or Proposition~1 in \cite{Dem2004}. 
\cqfd  
\vspace{3mm}

%%%%%%%%%%%%%%%% 
 
{\sc Proof of Theorem~\ref{Theordepthuniqueu}}.
From affine invariance, there is no loss of generality in assuming that~$\zb=0$, $\ub_0=(0,\ldots,0,1)'\in\R^d$ and that the path~$\ub_t$ is of the form~$\ub_t=(0,\ldots,0,\cos (t+\frac{\pi}{2}),\sin (t+\frac{\pi}{2}))'$, $t\in[0,\pi]$. Since~$P[\Pi\setminus \{0\}]=0$ for any hyperplane~$\Pi$ containing~$0$, we have, for any~$t\in[0,\pi]$,
%, write then (see~(\ref{EDcomput}))
$$
\gamma(t)
:=
e_0(\ub_t)
%=
%\frac{{\rm E}[ |u_t'(Z-z)| \mathbb{I}[u_t'(Z-z)\leq 0]]}{{\rm E}[|u_t'(Z-z)|]}
%=
%\frac{{\rm E}[ |u_t'(Z-z)| \mathbb{I}[u_t'(Z-z)< 0]]}{{\rm E}[ |u_t'(Z-z)| \mathbb{I}[u_t'(Z-z)< 0]]+{\rm E}[ |u_t'(Z-z)| \mathbb{I}[u_t'(Z-z)\geq 0]]}
%$$ 
%$$
%=
%\frac{-{\rm E}[ u_t'(Z-z) \mathbb{I}[u_t'(Z-z)< 0]]}{-{\rm E}[ u_t'(Z-z) \mathbb{I}[u_t'(Z-z)< 0]]+{\rm E}[ u_t'(Z-z) \mathbb{I}[u_t'(Z-z)\geq 0]]}
%$$
%$$
=
\frac{-h_{<}(t) }{- h_{<}(t) + h_{>}(t)} 
=
\frac{h_{<}(t) }{h_{<}(t)-h_{>}(t)}
,
$$
with
$$
h_{<}(t)
:=
\int_{\R^d} \ub_t' \yb \mathbb{I}[\ub_t' \yb < 0 ] \, dP(\yb)
%=:u_t' \mu_{t,<}
\
\textrm{ and }
\
h_{>}(t)
:=
\int_{\R^d} \ub_t' \yb \mathbb{I}[\ub_t' \yb > 0 ] \, dP(\yb)
%=:u_t' \mu_{t,\geq}
.
$$
Throughout the proof, we will use the notation
$
\mub={\rm E}[ \Zb ] 
$,
$
\mub_{t,<}
:=
{\rm E}[ \Zb \mathbb{I}[ \ub_t' \Zb <  0]]
$,
and
$
\mub_{t,>}
:=
{\rm E}[ \Zb \mathbb{I}[ \ub_t' \Zb >  0]]
$.
%and
%$
%\mu_{t,\geq}
%:=
%{\rm E}[ Z \mathbb{I}[ u_t' Z \geq  0]]
%$. 
%
Note that under the assumptions of the theorem, we have~$\mub=\mub_{0,>}+\mub_{0,<}$. 

We start by considering differentiability of (a)~$h_{<}(t)=\ub_t' \mub_{t,<}$ and~(b) $h_{>}(t)=\ub_t' \mub_{t,>}$.
(a) Since~$P[\Pi\setminus \{0\}]=0$ for any hyperplane~$\Pi$ containing~$0$, the mapping~$t\mapsto \ub_t' \yb \mathbb{I}[\ub_t' \yb  < 0 ]$ is $P$-almost everywhere differentiable at any~$t\in[0,\pi]$, with derivative~$ 
 t\mapsto   \dot{\ub}_t' \yb  \mathbb{I}[\ub_t' \yb<0]$, where we let~$\dot{\ub}_t:=(0,\ldots,0,-\sin (t+\frac{\pi}{2}),\cos (t+\frac{\pi}{2}))'$.
Since the function~$(t,\yb) \mapsto   \dot{\ub}_t' \yb \mathbb{I}[\ub_t' \yb<0]$ is upper-bounded by the $t$-independent $P$-integrable function~$\yb\mapsto \|\yb\|$, the mapping~$t\mapsto h_{<}(t)$ is differentiable at any~$t\in[0,\pi]$, with derivative~$\dot{h}_{<}(t)
:=\dot{\ub}_t' \mub_{t,<}$.
(b) Similarly,  for any~$\yb\in\R^d$, the mapping~$t\mapsto \ub_t' \yb \mathbb{I}[\ub_t' \yb  > 0 ]$ is $P$-almost everywhere differentiable  at any~$t\in[0,\pi]$, with derivative~$
 t\mapsto   \dot{\ub}_t' \yb \mathbb{I}[\ub_t' \yb>0]$. Since the function~$(t,\yb) \mapsto  \dot{\ub}_t' \yb \mathbb{I}[\ub_t' \yb>0]$ is still upper-bounded by the $t$-independent $P$-integrable function~$\yb\mapsto \|\yb\|$, the mapping~$t\mapsto h_{>}(t)$ is differentiable at any~$t\in[0,\pi]$, with derivative~$\dot{h}_{>}(t):=\dot{\ub}_t' \mub_{t,>}$.

We conclude that~$t\mapsto \gamma(t)$ is differentiable at any~$t\in[0,\pi]$, with a derivative~$\dot\gamma(t)$ that satisfies
\begin{eqnarray*}
	\lefteqn{
	(h_{<}(t)-h_{>}(t))^2
\dot{\gamma}(t)
=
\dot{h}_{<}(t) (h_{<}(t)-h_{>}(t))
-
h_{<}(t)(\dot{h}_{<}(t)-\dot{h}_{>}(t))
}
\\[2mm]
& & 
\hspace{7mm} 
=
h_{<}(t) \dot{h}_{>}(t)
-
h_{>}(t) \dot{h}_{<}(t)
=
(\ub_t' \mub_{t,<}) (\dot{\ub}_t' \mub_{t,>})
-
(\ub_t' \mub_{t,>}) (\dot{\ub}_t' \mub_{t,<})
.
\end{eqnarray*}
 
Let us introduce some further notation. For any~$t\in[0,\pi]$, write the projections of~$\mub_{t,<}$, $\mub_{t,>}$, and~$\mub$ onto the plane spanned by the last two vectors of the canonical basis of~$\R^d$  
as
$(0,\ldots,0,r_{t,<}\cos \alpha_{t,<},r_{t,<}\sin \alpha_{t,<})'$,
$(0,\ldots,0,r_{t,>}\cos \alpha_{t,>},r_{t,>}\sin \alpha_{t,>})'$,
and $(0,\ldots,0,r\cos \alpha,r\sin \alpha)'$, respecti\-vely, where all $r$'s are nonnegative and all~$\alpha$'s belong to~$[0,2\pi)$. Since it is assumed that~$E\!D(\zb,P)>0$, we must have
\begin{equation}
\label{spanp}
r_{t,>}>0
\quad\textrm{ and }\quad
\alpha_{t,>}\in (t,t+\pi)
\end{equation}
for any~$t\in[0,\pi]$
and
\begin{equation}
\label{spanm}
r_{t,<}>0
\quad\textrm{ and }\quad
\alpha_{t,<}\in [0,t) \cup (t+\pi,2\pi)
\end{equation}
for any~$t\in[0,\pi]$. Note that~$t\mapsto \alpha_{t,>}$ is monotone non-decreasing over~$[0,\pi]$ and that~$t\mapsto \alpha_{t,<}$ is monotone non-decreasing ``modulo~$2\pi$" over the same range. Finally, note also that
\begin{equation}
\label{g0mean}	
e_0(\ub_0)
=
\frac{{\rm E}[ |\ub_0'\Zb| \mathbb{I}[\ub_0'\Zb\leq 0]]}{{\rm E}[|\ub_0'\Zb|]}
\leq 
\frac{1}{2}
\end{equation}
(if not, then~$e_0(\ub_\pi)=1-e_0(\ub_0)<e_0(\ub_0)$, which contradicts the definition of~$\ub_0$). If~$e_0(\ub_0)=1/2$, then~$e_0(\ub_t)=1/2$ for any~$t\in[0,\pi]$ (if~$e_0(\ub)>1/2$ for some~$\ub\in\mathcal{C}$, then~$e_0(-\ub)=1-e_0(\ub)<1/2=e_0(\ub_0)$, a contradiction), so that the result holds with~$t_a=t_b=\pi$. We may thus assume that the inequality in~(\ref{g0mean}) is strict, which implies that~$\ub_0'\mub={\rm E}[\ub_0'\Zb]> 0$, hence that $r>0$ and~$\alpha\in [0,\pi]$. 

With the notation introduced above, we have
\begin{eqnarray*}
	\lefteqn{
\hspace{-1mm} 
 	(h_{<}(t)-h_{>}(t))^2
\dot{\gamma}(t)
=
(\ub_t' \mub_{t,<}) (\dot{\ub}_t' \mub_{t,>})  
-
(\ub_t' \mub_{t,>})  
(\dot{\ub}_t' \mub_{t,<}) 
}
\\[2mm]
& & 
\hspace{6mm} 
=
r_{t,<}
r_{t,>}
(\cos (t+{\textstyle\frac{\pi}{2}}) \cos \alpha_{t,<} + \sin (t+{\textstyle\frac{\pi}{2}}) \sin \alpha_{t,<} )
\\[2mm]
& & 
\hspace{32mm} 
\times (-\sin (t+{\textstyle\frac{\pi}{2}}) \cos \alpha_{t,>} + \cos (t+{\textstyle\frac{\pi}{2}}) \sin \alpha_{t,>})
\\[2mm]
& & 
\hspace{18mm} 
-
r_{t,<}
r_{t,>}
(\cos (t+{\textstyle\frac{\pi}{2}}) \cos \alpha_{t,>} + \sin (t+{\textstyle\frac{\pi}{2}}) \sin \alpha_{t,>} )
\\[2mm]
& & 
\hspace{38mm} 
\times (-\sin (t+{\textstyle\frac{\pi}{2}}) \cos \alpha_{t,<} + \cos (t+{\textstyle\frac{\pi}{2}}) \sin \alpha_{t,<})
\\[2mm]
& & 
\hspace{6mm} 
=
r_{t,<}
r_{t,>}
\cos (\alpha_{t,<}-(t+{\textstyle\frac{\pi}{2}}))
\sin (\alpha_{t,>}-(t+{\textstyle\frac{\pi}{2}}))
\\[2mm]
& & 
\hspace{18mm} 
-
r_{t,<}
r_{t,>}
\cos (\alpha_{t,>}-(t+{\textstyle\frac{\pi}{2}}))
\sin (\alpha_{t,<}-(t+{\textstyle\frac{\pi}{2}}))
\\[2mm]
& & 
\hspace{6mm} 
=
r_{t,<}
r_{t,>}
\sin (\alpha_{t,>}-\alpha_{t,<})
=:
\ell(t)
.
\end{eqnarray*}
 
 %%%%%%%%%%%%%%%%%%%%%%%%%%%%%%%%

Now, since~$\ub_0$ is a minimizer of~$e_0(\cdot)$ on~$\mathcal{C}$, 
$
\ell(0)
=
- r_{0,<}r_{0,>}\sin (\alpha_{0,<}-\alpha_{0,>})
= 0
.
$
Since~(\ref{spanp})-(\ref{spanm}) entail that~$\alpha_{0,<}>\alpha_{0,>}$, this yields~$\alpha_{0,<}=\alpha_{0,>}+\pi$. By using the identity~$\mub=\mub_{0,>}+\mub_{0,<}$ and the fact that~$\alpha\in[0,\pi]$, we conclude that~$\alpha=\alpha_{0,>}\in (0,\pi)$.       
Similarly, using the fact that~$\ub_\pi=-\ub_0$ is a maximizer of~$e_0(\cdot)$ on~$\mathcal{C}$ (this follows from the fact that~$e_0(-u)=1-e_0(u)$ for any~$u$), we have
$
\ell(\pi)
=
- r_{\pi,<}r_{\pi,>}\sin (\alpha_{\pi,<}-\alpha_{\pi,>})
= 0
,
$
which implies that~$\alpha_{\pi,>}=\alpha_{\pi,<}+\pi$ (recall that we cannot have~$\alpha_{\pi,<}=\pi$, nor~$0$). Thus~$\pi<\alpha_{\pi,>}<2\pi$. By using the identity~$\mub=\mub_{\pi,>}+\mub_{\pi,<}$ and the fact that~$\alpha\in(0,\pi)$, we conclude that~$\alpha_{\pi,<}=\alpha$. 
 
Now, fix~$t_0\in[0,\pi]$ with~$\ell(t_0)=- r_{t_0,<}r_{t_0,>}\sin (\alpha_{t_0,<}-\alpha_{t_0,>})\neq 0$ (if there is no such~$t_0$, then the result holds with~$t_a=t_b=\pi$). 
Monotonicity of~$t\mapsto \alpha_{t,>}$ yields $\alpha_{t_0,>}\geq \alpha_{0,>}= \alpha$. Since~$\alpha_{t_0,>}=\alpha$ would lead to~$\alpha_{t_0,<}=\alpha_{t_0,>}+\pi$ (due to~$\mu=\mu_{t_0,>}+\mu_{t_0,<}$), hence to~$\ell(t_0)=0$, we must actually have~$\alpha_{t_0,>}>\alpha$.  

Pick then an arbitrary~$t\in[t_0,\pi)$. Monotonicity of~$t\mapsto \alpha_{t,>}$ then yields $\alpha_{t,>}\geq \alpha_{t_0,>}>\alpha$, 
hence~$\alpha_{t,>}\in (\alpha,t+\pi)$. Since~$\mub=\mub_{t,>}+\mub_{t,<}$, we must have~$\alpha_{t,<}\in (\alpha_{t,>}+\pi,t+2\pi]\!\mod 2\pi$, that is, (a) $\alpha_{t,<}\in(\alpha_{t,>}+\pi ,2\pi)$ or (b) $\alpha_{t,<}\in[0,t)$. In case~(a), we have~$\alpha_{t,<}-\alpha_{t,>}\in(\pi,2\pi)$, so that $\ell(t)=
-r_{t,<}r_{t,>}\sin (\alpha_{t,<}-\alpha_{t,>})>0$. In case~(b), in view of the monotonicity (modulo~$2\pi$) of~$t\mapsto \alpha_{t,<}$, we have~$\alpha_{t,<}\in[0,\alpha_{\pi,<}=\alpha]$. Therefore, the identity~$\mub=\mub_{t,>}+\mub_{t,<}$ implies that~$\alpha_{t,<}<t<\alpha_{t,>} \leq \alpha_{t,<}+\pi$. Therefore, $\ell(t)=r_{t,<}r_{t,>}\sin (\alpha_{t,>}-\alpha_{t,<})\geq 0$. 

Assume that~$\ell(t)=0$. As shown above, we must thus be in case~(b). Then~$\alpha_{t,>} = \alpha_{t,<}+\pi$, so that (still due to~$\mub=\mub_{t,>}+\mub_{t,<}$) we must have~$\alpha_{t,<}=\alpha=\alpha_{\pi,<}$. Monotonicity then implies that for any~$t'\in [t,\pi)$, we have~$\alpha_{t',<}=\alpha$, which, in view of~$\mub=\mub_{t',>}+\mub_{t',<}$ yields~$\alpha_{t',>}=\alpha+\pi$. Consequently, we have~$\ell(t')=0$, which establishes the result.
\cqfd  
\vspace{3mm}

%%%%%%%%%%%%%%%%

% \subsection{Proofs of Section~\ref{secrisks}}
%\label{secAppSec6}

{\sc Proof of Theorem~\ref{Theoraxioms}}. 
Throughout the proof, $\rho(t)=t^2$ is the quadratic loss function. 
(i)
First note that since~$\ub\in \R^d_+$, we have~$\ub'\Xb\leq \ub'\Yb$ almost surely, so that the monotonicity of (univariate) expectiles entails that $\theta^\rho_{\alpha}(\ub'\Xb)\leq \theta^\rho_{\alpha}(\ub'\Yb)$. It trivially follows that~$H^{\rho}_{\alpha,\ub}(\Yb)\subset H^{\rho}_{\alpha,\ub}(\Xb)\subset H^{\rho}_{\alpha,\ub}(\Xb) \oplus \R^d_+$. To establish the other inclusion, fix~$\zb\in H^{\rho}_{\alpha,\ub}(\Xb)$. We may assume that~$\zb\notin H^{\rho}_{\alpha,\ub}(\Yb)$ (if~$\zb\in H^{\rho}_{\alpha,\ub}(\Yb)$, then~$\zb=\zb+0\in H^{\rho}_{\alpha,\ub}(\Yb) \oplus \R^d_-$). We then have
$$
\zb
=
(\theta^\rho_{\alpha}(\ub'\Yb)\ub+(I_d-\ub\ub')\zb)+(\ub'\zb-\theta^\rho_{\alpha}(\ub'\Yb))\ub
=: \zb_0+\zb_1
.
$$ 
Since~$\ub'\zb_0=\theta^\rho_{\alpha}(\ub'\Yb)$, we have that~$\zb_0\in H^{\rho}_{\alpha,\ub}(\Yb)$. Since~$\ub'\zb-\theta^\rho_{\alpha}(\ub'\Yb)<0$ (recall that~$\zb\notin H^{\rho}_{\alpha,\ub}(\Yb)$) and~$\ub\in\R^d_+$, we also have~$\zb_1\in\R^d_-$. This shows that $H^{\rho}_{\alpha,\ub}(\Xb)\subset H^{\rho}_{\alpha,\ub}(\Yb) \oplus \R^d_-$. 
(ii)
Let~$\zb\in H^{\rho}_{\alpha,\ub}(\Xb+\Yb)$ and decompose it into~$
\zb
=
(\theta^\rho_{\alpha}(\ub'\Xb)\ub+(I_d-\ub\ub')\zb)+(\ub'\zb-\theta^\rho_{\alpha}(\ub'\Xb))\ub
=:
\zb_0+\zb_1$. Obviously,~$\zb_0\in H^{\rho}_{\alpha,\ub}(\Xb)$. As for~$\zb_1$, the superadditivity of univariate expectiles for~$\alpha\in(0,\frac{1}{2}]$ implies that $\ub'\zb_1=\ub'\zb-\theta^\rho_{\alpha}(\ub'\Xb)\geq \theta^\rho_{\alpha}(\ub'(\Xb+\Yb))- \theta^\rho_{\alpha}(\ub'\Xb)\geq \theta^\rho_{\alpha}(\ub'\Yb)$, which shows that~$\zb_1\in H^{\rho}_{\alpha,\ub}(\Yb)$. 
(iii) 
If~$\zb_0\in H^{\rho}_{\alpha,\ub}(\Xb)$ and~$\zb_1\in H^{\rho}_{\alpha,\ub}(\Yb)$, then the subadditivity of univariate expectiles for~$\alpha\in[\frac{1}{2},1)$ readily yields~$\ub'(\zb_0+\zb_1)\geq \theta^\rho_{\alpha}(\ub'\Xb)+\theta^\rho_{\alpha}(\ub'\Yb)\geq \theta^\rho_{\alpha}(\ub'(\Xb+\Yb))$, which shows that~$\zb_0+\zb_1\in H^{\rho}_{\alpha,\ub}(\Xb+\Yb)$. 
(iv) The result trivially follows from the monotonicity of univariate expectiles with respect to their order~$\alpha$.
\cqfd

%%%%%%%%%%%%%%%%%%%%%%%%%%%%%%%%%%%%%%%%%%%%%%%%%%

\end{document}